\documentclass[twoside,a4paper,12pt,centertags]{amsart}

\usepackage{amsbsy}
\usepackage{amscd}
\usepackage{amsmath,amssymb,verbatim,vmargin,slashed,graphicx,subcaption}
\usepackage[all]{xy}
\usepackage[bookmarks=true]{hyperref}   
\usepackage{color}
\usepackage{tikz}
\usetikzlibrary{arrows,decorations.markings}

\theoremstyle{plain}
\newtheorem{thm}{Theorem}[section]

\theoremstyle{definition}

\newtheorem{rem}[thm]{Remark}

\newtheorem{ex}[thm]{Example}

\newcommand{\myvec}[1]{{\bf {#1}}}
\newcommand{\mv}[1]{\boldsymbol{#1}} 

\newcommand{\japsigma}{\langle\sigma\rangle}

\mathchardef\semic="303B

\newcommand{\diag}{\text{{\rm diag}}\,}
\newcommand{\dirac}{{\mathbf D}}
\newcommand{\rev}[1]{\overline{#1}}

\newcommand{\R}{{\mathbf R}}
\newcommand{\C}{{\mathbf C}}

\newcommand{\mH}{{\mathcal H}}
\newcommand{\mE}{{\mathcal E}}

\DeclareMathOperator{\re}{Re}
\newcommand{\im}{\text{{\rm Im}}\,}
\newcommand{\sett}[2]{ \{ #1 \, \semic \, #2 \} }

\newcommand{\nul}{\textsf{N}}
\newcommand{\ran}{\textsf{R}}

\newcommand{\clos}[1]{\overline{#1}}
\newcommand{\barint}{\mbox{$ave \int$}}

\newcommand{\curl}{{\text{{\rm curl}}}}

\newcommand{\pd}{\partial}
\newcommand{\pv}{\text{{\rm p.v.\!}}}

\def\barint_#1{\mathchoice
            {\mathop{\vrule width 6pt
height 3 pt depth -2.5pt
                    \kern -8.8pt
\intop}\nolimits_{#1}}%
            {\mathop{\vrule width 5pt height
3 pt depth -2.6pt
                    \kern -6.5pt
\intop}\nolimits_{#1}}%
            {\mathop{\vrule width 5pt height
3 pt depth -2.6pt
                    \kern -6pt
\intop}\nolimits_{#1}}%
            {\mathop{\vrule width 5pt height
3 pt depth -2.6pt
          \kern -6pt \intop}\nolimits_{#1}}}

\usepackage{color}

\definecolor{gr}{rgb}   {0.,   0.8,   0. }
\definecolor{bl}{rgb}   {0.,   0.5,   1. }
\definecolor{mg}{rgb}   {0.7,  0.,    0.7}

\begin{document}

\title[An efficient  full-wave  solver for eddy currents]
{An efficient  full-wave  solver for eddy currents}

\author[Johan Helsing]{Johan Helsing$\,^1$}
\author[Anders Karlsson]{Anders Karlsson$\,^2$}
\author[Andreas Ros\'en]{Andreas Ros\'en$\,^3$}
\thanks{$^1\,$Centre for Mathematical Sciences, Lund
    University, Box 118, 221 00 Lund, Sweden (johan.helsing@math.lth.se).}
\thanks{$^2\,$Electrical and Information Technology, Lund
    University, Box 118, 221 00 Lund, Sweden 
    (anders.karlsson@eit.lth.se).}
\thanks{$^3\,$Mathematical Sciences, Chalmers University of
    Technology and the University of Gothenburg, 412 96 G{\"o}teborg,
    Sweden (andreas.rosen@chalmers.se).}

\begin{abstract}
An integral equation reformulation of the Maxwell transmission
   problem is presented. The reformulation uses techniques such as
   tuning of free parameters and augmentation of  close-to-rank-deficient  operators. It is 
   designed for the eddy current regime and works  both for surfaces of genus $0$
   and $1$. Well-conditioned systems and field representations are
   obtained despite the Maxwell transmission problem being
   ill-conditioned for genus $1$ surfaces due  to  the presence of Neumann
   eigenfields. Furthermore, it is shown that these eigenfields, for
   ordinary conductors in the eddy current regime, are different from the
    classical   Neumann eigenfields for superconductors.  
Numerical  examples,  based on
   the reformulation, give an unprecedented $13$-digit accuracy both
   for transmitted and scattered fields.
\end{abstract}

\keywords{
Maxwell transmission problem,
Eddy current,
Boundary integral equation,
Neumann eigenfield,
Low-frequency breakdown}
 
\subjclass[2010]{15A66, 35Q61, 45E05, 78M99}

\maketitle

\section{Introduction}

This work concerns the Maxwell transmission problem (MTP), 
which is the problem of computing the electromagnetic wave 
transmitted through and scattered from a bounded object $\Omega_+\subset \R^3$, 
given an incident time-harmonic electromagnetic wave in the 
exterior region $\Omega_-=\R^3\setminus\clos\Omega_+$. 
Consider $\Omega_+$ with boundary surface $\Gamma$, 
generalized diameter $L=\sup\sett{|x-y|}{x,y\in \Omega_+}$ 
and $\Omega_-$ being vacuum. The corresponding wavenumbers are 
\begin{align}
   k_+&=\omega\sqrt{(\epsilon_0\epsilon_r+ i\sigma/\omega)\mu_0}\,, \label{eq:kplusformula}\\
   k_-&=\omega\sqrt{\epsilon_0\mu_0}\,,\label{eq:kminusformula}
\end{align}
where $\omega, \epsilon_0, \mu_0, \epsilon_r, \sigma$ denote 
frequency, permittivity and permeability of vacuum, and relative 
permittivity and conductivity of $\Omega_+$. 
In terms of the wavenumbers, the conductivity is  
$\sigma= \re(k_+^2/(i\eta_0 k_-))$ and the skin depth equals
$1/\im(k_+)$. 
Since magnetic materials often have non-linear properties and exhibit 
hysteresis we restrict ourselves to non-magnetic materials, for which the linear MTP \eqref{eq:maxwtranspr} below is an accurate physical model. 

  As $\sigma\to\infty$ for fixed $\omega>0$, the object $\Omega_+$
  approaches a perfect electric conductor (PEC) with zero internal
  fields and with the electric field normal and the magnetic field
  tangential to $\Gamma$. The PEC boundary condition, that the
  electric field is normal to $\Gamma$, is only proper if the skin
  depth is much smaller than the diameter of the object. This is why
  the PEC boundary condition cannot be applied to metals at low
  frequencies. On the other hand it applies, with good accuracy, to a
  superconductor for frequencies ranging from zero up to very high
  values. The reason is that a superconductor is also a perfect
  diamagnet. Therefore we refer to the limit $\sigma\to\infty$ for the
  MTP as the {\em superconducting limit}, despite restricting
  ourselves to non-magnetic materials. 

In scattering theory, the regime $k_-L\ll 1$ is referred to as the Rayleigh regime. Here the scattered far fields are accurately determined by the induced electric and magnetic dipole moments in $\Omega_+$, an approximation widely used in optics and microwave theory~\cite[Sec.~10.1]{Jackson:98} and with application in radar, lidar, and radio communication~\cite{SchneiderETAL:05,HuangETAL:19,
Shamsan:19}. The scattered near fields from sub-wavelength objects are important in non-destructive testing, where they serve as input data to solvers that extract information about the objects'
interior~\cite{GarciaMartinetal:11,EgorovETAL:17}. The design of integrated circuits often requires the determination of inductances, capacitances, and resistances of sub-wavelength components based on both transmitted and scattered fields~\cite{ZhuETAL:05}. 

\begin{figure}[t]
\centering
\includegraphics[height=90mm]{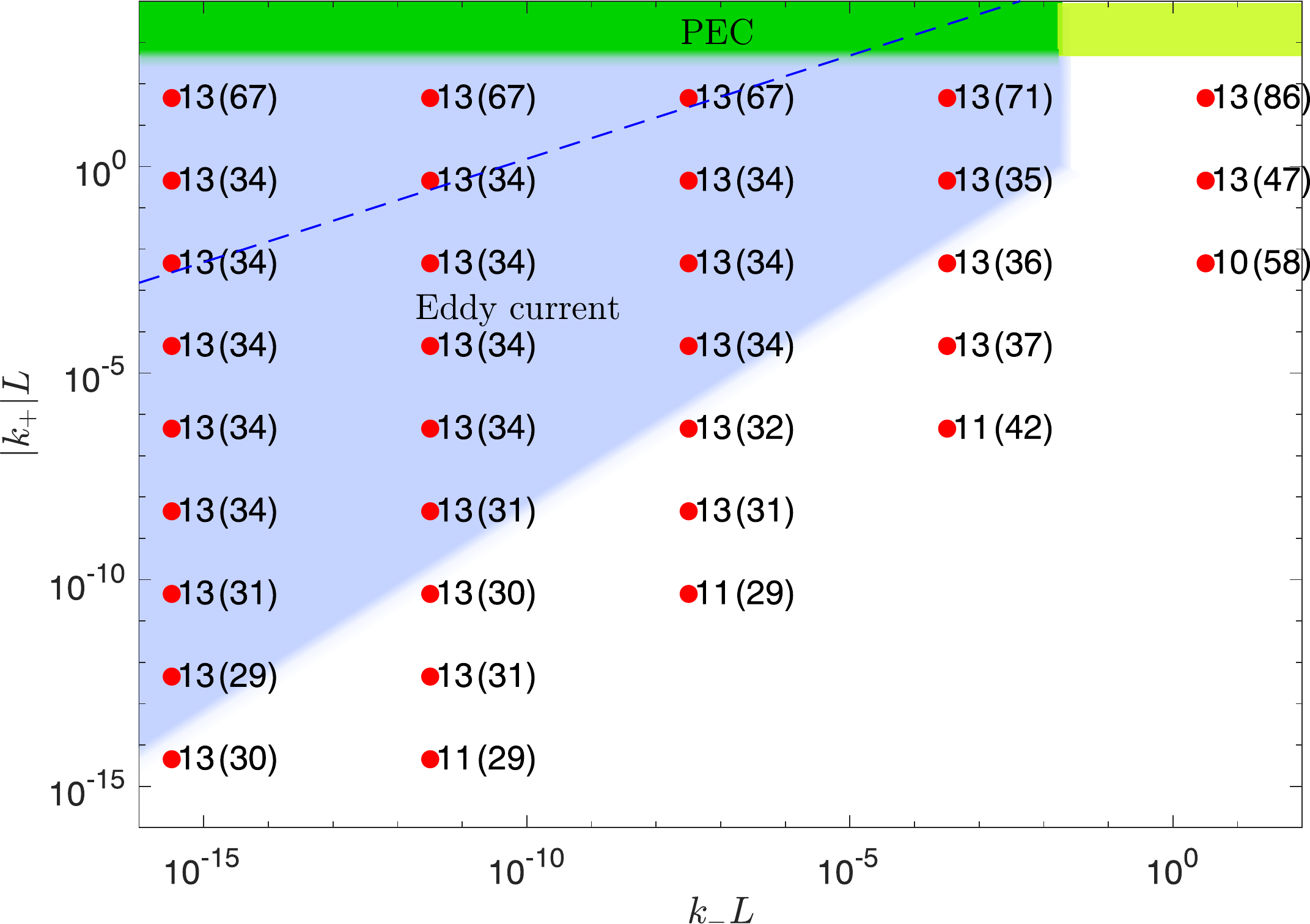}
\caption{\sf 
   Performance of Dirac (B-aug1) on the ``starfish
   torus''~\eqref{eq:starfishtorus} with incident partial
   waves~\eqref{eq:partialwave} and $\arg(k_+)=\pi/4$. The expression
   $Y(X)$ at red points $(k_-L,|k_+|L)$ says that GMRES needs $X$
   iterations and that $Y$-digit  accuracy~\eqref{eq:Ydigits},  
   or better, is achieved in
   each of the fields  $\{E^+, E^-, H^+, H^-\}$  at all
   $90,\!000$ field points in the computational domain. The PEC regime
   is in dark green and light green. The eddy current regime is in blue
   and dark green (vacuum in $\Omega_-$). The dashed line is the upper
   limit of realizable $|k_+|L$ for $L=1$ m (silver in $\Omega_+$,
   vacuum in $\Omega_-$). Red points outside the eddy current regime,
   with $\arg(k_+)=\pi/4$, can be realized if $\Omega_-$ is a
   dielectric.}
\label{fig:quadrant0genus1}
\end{figure}

When $\sigma/\omega \gg\epsilon_0\epsilon_r$, then
$\arg(k_+)\approx \pi/4$ in \eqref{eq:kplusformula} and
\begin{equation}
  |k_+|L\approx \sqrt{\eta_0\sigma L}\sqrt{k_-L},
\end{equation}
where also  $|k_+|L\gg k_-L$.  
We refer to $\sigma L$ as the {\em scaled conductivity} of $\Omega_+$. 
Here $\eta_0\sigma L$ and $k_-L$ are dimensionless, and 
$\eta_0=\sqrt{\mu_0/\epsilon_0}\approx 377$ Ohm is the wave impedance of vacuum.
The pair of wavenumbers $(k_-, k_+)$ is said to be in the {\em eddy current regime} if 
$L \ll1/k_-$ and $|k_+| \gg k_-$.  
The eddy current regime is in blue
and dark green in Figure~\ref{fig:quadrant0genus1}.
There is an upper limit on conductivity $\sigma\lesssim 6\cdot 10^7$ S/m in
ordinary conductors, that is, non-superconducting materials. 
Thus $|k_+|L\lesssim C\sqrt{k_-L}$, where $C\approx 1.5\cdot 10^5\sqrt{L}$
with $L$ measured in meters,
is the physical part of the eddy current regime for objects of given size $L$.
When $L=1$ m this limit is the dashed line in Figure~\ref{fig:quadrant0genus1}. 
The low-frequency asymptote for any ordinary conductor with $\sigma>0$
 and $L=1$ m,  is
a line parallel to, and below, this dashed line. 
The dark green and light green areas  in Figure~\ref{fig:quadrant0genus1} is the regime where the 
PEC boundary condition is applicable with  a reasonably small 
relative error. 
To stay in this area when $\omega\to 0$, 
it is necessary that $\sigma\to \infty$. 
     It is seen that for $L=1$ m and $k_-L< 10^{-5}$ the PEC approximation is invalid for
     ordinary conductors, but it holds for superconductors. This is so since  the  surface resistance of a superconductor is  low enough to be 
     considered to be zero in the entire eddy current regime, see~\cite{Halbritter:74}, and zero surface resistance implies the PEC boundary      
     condition on $\Gamma$. 

We refer to solvers based on boundary integral equations (BIEs) that model the full MTP as {\em
   full-wave solvers}, in contrast to solvers which build on an
   approximation to the MTP.
When $|k_+|L\ll 1$, a standard approximation is to determine 
the dipole moments by solving Laplace's equation, and for 
$|k_+|L\gg 1$ the standard approximation is to use  the PEC boundary 
condition. 
Between these two extremes it is necessary to solve
the MTP without approximations.  
The full-wave solvers for the MTP in the eddy current regime  
that we have found in the literature are \cite{RuckerETAL:95, ZhuETAL:05, Andriulli20}.
Rucker et al.~\cite{RuckerETAL:95} give an overview of full-wave solvers, which give
at best  a relative  error of $1 \%$ in scattering situations comparable to our Figure~\ref{fig:genus0}.
Zhu et al.~\cite{ZhuETAL:05} describe a full-wave solver used for an open source program
FastImp.
Chhim et al.~\cite{Andriulli20} present a full-wave solver based on the PMCHWT BIE. 
We lack data to judge the accuracy of the BIEs in \cite{ZhuETAL:05} 
and \cite{Andriulli20}. 
The analysis in \cite{Andriulli20} appears to be limited to $\omega\to 0$ for fixed $\sigma$, 
leaving a possible gap 
to the green PEC regime in Figure~\ref{fig:quadrant0genus1}.
Also, numerical evaluation of the field representations is missing in \cite{Andriulli20},
cf. Section~\ref{sec:conditioning}.
Other BIEs rather solve the quasi-static approximation of the MTP obtained by neglecting 
the dispacement current in Amp{\`e}re's law, which limits their validity.
For justifications of such eddy current models, see \cite{AmmariBuffaN:00, PPicard:17, BonnetDemaldent19}.
As these rather sparse results in the literature indicate, it is indeed a challenging problem 
to design BIEs for the MTP in the eddy current regime.
As we discuss below, 
reasons for this is that the fields may differ much in size and
for $\Gamma$ of non-zero genus the MTP itself is actually ill-posed 
in the eddy current regime as $k_-\to 0$. 
By the limit being in the eddy current regime, we mean that $|k_+|/k_-\to \infty$.
More precisely, in such limits there is an incident field for which some of the transmitted 
and scattered fields differ drastically from their generic magnitude.
We refer to Sections~\ref{sec:physics}, \ref{sec:hightorus} and \ref{sec:finitetorus}
for details. 
See  Hiptmair~\cite{Hiptmair:07}  for more background on eddy current computations.

For our numerical method to be efficient when $\im(k_+)$ is
   large, we always assume that 
$|k_+|L\lesssim 50$, so our standing assumption  in numerical evaluations   is that
\begin{equation}  \label{eq:goodcond}
0<k_-L\ll |k_+|L\lesssim 50.
\end{equation}  
In this paper, we achieve BIEs that
   compute all the fields to a minimum of $13$ accurate digits 
 in the entire regime \eqref{eq:goodcond}, for $\Gamma$ of genus $0$ as well as
genus $1$.
See Figure~\ref{fig:quadrant0genus0} in Section~\ref{sec:perform} for genus $0$
and Figure~\ref{fig:quadrant0genus1} for genus $1$.  
Our BIEs appear to be the only known full-wave solvers for the MTP
  that compute 
all fields accurately and fast in all  of  \eqref{eq:goodcond}. 
We point out that although $\arg(k_-)=0$ and $\arg(k_+)=\pi/4$ in
all our numerical examples, there is no approximation of the full MTP
\eqref{eq:maxwtranspr} involved, and a known permittivity can be included in $k_+$. 

The Dirac BIE from~\cite{HelsRose20,HelsKarlRos20} 
is the starting point for the present
work and from now on referred to as Dirac (A). This BIE is based on
the embedding of Maxwell's equations into an elliptic Dirac equation
and a Cauchy integral representation for the fields (Eq. \eqref{eq:projdens} below).
Schematically, given the incident wave $g$ on $\Gamma$,
we have
\begin{equation}   \label{eq:abstractIER}
  h\mapsto (F^+,F^-) \mapsto g,
\end{equation}
and solve the Dirac BIE for the density $h$ on $\Gamma$, from which the
ansatz/Cauchy representation yields the transmitted 
electric and magnetic fields $F^+=(E^+, H^+)$ 
in $\Omega_+$ and the scattered fields $F^-=(E^-,H^-)$ in $\Omega_-$.
The Dirac BIE is a size $8\times 8$ block system using $50$, not all distinct, 
integral operators of double and single layer type, which
can be used on any Lipschitz regular boundary $\Gamma$.
It has $12$ free parameters, as recalled in \eqref{eq:12param} below,
and these can be chosen to avoid
false eigenwavenumbers for all passive materials.
As for low-frequency breakdown,
the only regimes found in~\cite{HelsRose20,HelsKarlRos20} where 
the Dirac (A)  exhibits false
eigenwavenumbers is when $k_\pm\to 0$ at the same time as
$\hat k= k_+/k_-\to\infty$ or $\hat k\to 0$.  
(We refer to~\cite[Sec.~3]{HelsKarlRos20}  
for a discussion about the
notions of eigenwavenumbers and low-frequency breakdown. 
See also Section~\ref{sec:chooseaug} for the notions of 
dense-mesh breakdown and topological low-frequency breakdown.)  
This corresponds to the false eigenwavenumber at $x=-1$ 
in~\cite[Fig.~9(b)]{HelsKarlRos20},  and it should be noted
that this single peak contains the whole eddy current regime
shown in Figure~\ref{fig:quadrant0genus1}.  
(The reverse eddy current regime corresponding to 
$x= +1$, where the roles of $\Omega_\pm$ have been swapped, is 
less important  in applications, and we omit the details.)  

The design of all Dirac BIEs starts by {\em tuning}, that is, carefully assigning values to the free
parameters, to avoid false eigenwavenumbers and optimize numerical performance.
The main result of this paper consists of two new parameter choices
   for Dirac BIEs in the eddy current regime, referred to as Dirac
   (A$\infty$) and Dirac (B), which, at a low cost, enable us to
   simultaneously compute the four fields  $\{E^+, E^-, H^+, H^-\}$  to  
  almost full machine precision. A key problem in the
   eddy current regime is
   that the fields  $\{E^+, E^-, H^+, H^-\}$  may differ much in 
size, and tuning the parameters to account for this becomes a non-trivial
matter. In particular, the transmitted electric field $E^+$ is typically much 
smaller than the other fields. It is nevertheless important to compute also 
$E^+$ with a small relative error, since the measurable eddy current
$J=\sigma E^+$ will have the same relative error. 

An important take-home message from the present paper
is that in order to stably solve the MTP through \eqref{eq:abstractIER}, 
it is necessary both (a) to have a well-conditioned 
BIE for computing the density $h$ from the boundary datum $g$, and (b)
to have a well-conditioned representation of the fields, for computing the fields $F^\pm$ from the density $h$. 
For Dirac (A$\infty$), it is in general (b)  that  is problematic since its field
evaluation allows for large fields, see \eqref{eq:Ascale}, 
and this sometimes leads to  cancellation  and loss
of accuracy for small fields like $E^+$.
Dirac (B) is adapted to the small size of $E^+$, 
see \eqref{eq:Bscale}.  It is then rather (a)
        which is challenging, but we still obtain a well-conditioned system. 
To eliminate  null densities  $h$ we use {\em augmentation}, that
is we suitably add a  finite-rank  matrix to the system to be solved, without
changing the physical solutions. We 
explain in Sections~\ref{sec:chooseparam} and \ref{sec:chooseaug} our general process for designing Dirac BIEs.
The augmentation techniques explained in 
Section~\ref{sec:chooseaug} are of 
independent interest beyond the MTP.
After such augmentations,  which for our BIEs are needed only as $k_-\to 0$,  
we obtain BIEs referred to as (A$\infty$-aug)
and (B-aug$0/1$) for the MTP. 
Our augmentations build on a careful analysis of null spaces and ranges of the quasi-static limit operators, given in the Appendix. 
In absense of such detailed knowledge one could try random augmentations
as in \cite{SifuentesGimbutasGreengard:15}. However, most of our operators require a careful choice of
augmentation or else well-conditioning and accuracy will be lost.

Turning to objects $\Omega_+$ with non-zero genus, an additional
difficulty at high scaled conductivities is that Neumann eigenfields, 
similar to those in PEC scattering \cite{AndiulliETAL09, EpsteinETAL13},
appear also in the MTP  in the eddy current regime  as $k_-\to 0$. 
More surprisingly, such Neumann eigenfields are present at low frequencies even 
for finite  non-zero  scaled conductivities,  although the terminology ``eigenfield'' may
not be appropriate to describe  this phenomenon. 
  The same Neumann eigenfield appears in the low-frequency limit for all ordinary conductors, regardless of the value of $\sigma$.
We have not found this Neumann eigenfield for ordinary conductors in
the literature. In particular it is not related to the notion of $k$-Neumann
fields from \cite{Kress:86, EpsGreNei10, EpsGreNei15}, 
since it is a static field and the magnetic field is not 
tangential on the surface.
However, when $\sigma\to \infty$ as $\omega\to 0$,  our 
Neumann eigenfield approaches the classical static PEC Neumann eigenfield. 
See~\eqref{eq:defneigfield}  
and  numerical examples and discussion in Sections~\ref{sec:hightorus} and \ref{sec:finitetorus},  
where it is shown that there exists an incident field for which some of the transmitted and scattered fields differ drastically from their 
generic magnitude. 
This means that for $\Gamma$ of non-zero genus, the MTP itself is 
ill-conditioned in the eddy current regime.
We discuss the Neumann eigenfields in some detail in Section~\ref{sec:physics},
and here only stress one important point:
according to our discussion above, the physical eddy current eigenfields appearing 
in ordinary conductors are those shown in 
Figure~\ref{fig:Neumanneig}(g,h,i).
The Neumann eigenfields computed  with  the PEC boundary condition   appear only 
in superconductors. See Figure~\ref{fig:Neumanneig}(a,b,c).

 At the end of the paper, we give 
in Section~\ref{sec:Maxwess} a proof of the result, which we have not 
found in the literature, that the essential spectrum of the MTP coincides with that of 
the Neumann--Poincar\'e operator. 
 This proof further illustrates the flexibility of the free Dirac
   parameters.
We conclude the paper in Section~\ref{sec:concl} with some remarks
on the usage of (A$\infty$-aug) and (B-aug$0/1$).

\section{The Maxwell and the two Helmholtz problems}

We fix notation for the remainder of the paper.
Let $\Omega_+$ be a bounded, connected
domain in $\R^3$ with Lipschitz regular boundary surface $\Gamma$ and an unbounded, connected, exterior $\Omega_-$. 
Let $L=\sup\sett{|x-y|}{x,y\in \Omega_+}$ be the generalized diameter of $\Omega_+$. Starting from Section~\ref{sec:finitecond}, we use unit
of length so that $L$ is of order $1$, 
which is convenient in numerical computations.
The outward unit normal 
on $\Gamma$ is $\nu$,  surface measure is $d\Gamma$,  and $\{\nu,\tau,\theta\}$
denotes a positive ON-frame on $\Gamma$. 
(Singularities of the frame on a null set does not present a problem.)  
In $\R^3$, $\{\mv\rho, \mv\theta, \mv z\}$ denotes the standard cylindrical
ON-frame. 
We consider time-harmonic
fields with time dependence $e^{-i \omega t}$, and angular
frequency $\omega>0$. 
The domains $\Omega_\pm$ are homogeneous with material properties described by wavenumbers $k_\pm$, and we write 
\begin{equation}
\hat k= k_+/k_-.
\end{equation} 
All our numerical examples use $\arg(k_-)= 0$ and $\arg(k_+)=
   \pi/4$, but our BIEs (A$\infty$) and (B-aug0/1) apply to more
   general  wavenumbers  $\im(k_\pm)\ge 0$
   satisfying~\eqref{eq:goodcond}. 

We consider Maxwell transmission problems MTP($k_-,k_+,\alpha$) 
\begin{equation}  \label{eq:maxwtranspr}
\begin{cases}
  \nu\times E^+= \nu\times (E^0+ E^-), & x\in \Gamma,\\
   \nu\times H^+= (\hat k^2/\alpha) \nu\times (H^0+H^-), & x\in \Gamma,\\
   \nabla\times E^+= ik_+ (\hat k^{-1} H^+), \quad\nabla\times (\hat k^{-1} H^+)= -ik_+ E^+, &x\in\Omega_+,\\
   \nabla\times E^-= ik_- H^-, \quad\nabla\times H^-= -ik_- E^-, &x\in\Omega_-,\\
   x/|x|\times E^- - H^-= o(|x|^{-1}e^{\im(k_-) |x|}), & x\to\infty,\\
      x/|x|\times H^- - E^-= o(|x|^{-1}e^{\im(k_-) |x|}), & x\to\infty,\\
   \end{cases}
\end{equation} 
where $E^0$ and $H^0$ are the incident fields from sources in 
$\Omega_-$,  and we want to solve for $E^\pm$ and $H^\pm$.  
In particular
\begin{equation}  \label{eq:DirEminusaug}
  \int_\Gamma \nu\cdot  E^- d\Gamma=0
\end{equation}
holds by the jump relations and the divergence theorem.   
For a discussion of the $L_2$ topology considered for the fields and the
corresponding trace space, we refer to~\cite[Sec.~5]{HelsRose20}.
The  physical Maxwell  transmission problem 
that  we aim to solve is MTP($k_-,k_+$)= MTP($k_-,k_+,\hat k^2$), with 
$\alpha= \hat k^2$, where the tangential parts of both the 
electric and magnetic fields are continuous across $\Gamma$.
However, we also use auxiliary MTPs with other values of the parameter
$\alpha$.

\begin{rem}  \label{rem:BH}
The $k_\pm$ are related to the total permittivities $\epsilon_\pm$ and permeabilities $\mu_\pm$ by $k_\pm=\omega\sqrt{\epsilon_\pm\mu_\pm}$
in $\Omega_\pm$ respectively. 
  We follow the convention from~\cite{HelsKarlRos20} where in all $\R^3$,
   the $H$ field
  is the magnetic field rescaled by the wave impedance
$\sqrt{\mu_-/\epsilon_-}$.
  The field $\hat k^{-1}H^+$, appearing in \eqref{eq:maxwtranspr}, 
  is the magnetic field rescaled by the wave impedance
$\sqrt{\mu_+/\epsilon_+}$ in $\Omega_+$.
This field is natural when formulating Maxwell's equations as a Dirac equation, 
and was denoted $B^+$ in~\cite{HelsRose20}.
\end{rem}

Besides MTPs, we also consider auxiliary 
Helmholtz transmission problems HTP($k_-,k_+,\beta$) 
\begin{equation}   \label{eq:Helmtransmprobl}
\begin{cases}
   u^+= u^0+u^-, & x\in\Gamma, \\
   \pd_{\nu} u^+={\beta}\pd_{\nu} (u^0+ u^-), & x\in\Gamma,\\
    \Delta U^+ +k_+^2 U^+=0, &x\in\Omega_+,\\
   \Delta U^- +k_-^2 U^-=0, &x\in\Omega_-,\\
    \pd_{x/|x|}U^- -ik_- U^-=o(|x|^{-1}e^{\im(k_-) |x|}), & x\to\infty,
   \end{cases}
\end{equation} 
where $u^0$ is the trace of the incident wave $U^0$, and we
want to solve for $U^\pm$.
We use the elliptic Dirac type equation
\begin{equation}   \label{eq:hodgedirac}
    \begin{bmatrix} 0 & \nabla\cdot & 0 & 0 \\ 
    \nabla & 0 & -\nabla\times & 0 \\ 
    0 & \nabla\times & 0 & \nabla \\
    0 & 0 & \nabla\cdot & 0 \end{bmatrix}
  \begin{bmatrix} F_0 \\ F_1 \\ F_2 \\ F_3 \end{bmatrix}
  = ik
  \begin{bmatrix} F_0 \\ F_1 \\ F_2 \\ F_3 \end{bmatrix},
\end{equation}
for two scalar fields $F_0$ and $F_3$, and two vector fields
$F_1$ and $F_2$,
which embeds one Maxwell and two Helmholtz  equations. 
Indeed, for $F_0=F_3=0$, \eqref{eq:hodgedirac} amounts
to Maxwell's equations for $F_1=E^+$ and 
\begin{equation}   \label{eq:Hplus}
F_2=\hat k^{-1}H^+
\end{equation}
 in $\Omega_+$, and  for  $F_1=E^-$ and $F_2=H^-$ in $\Omega_-$.
Moreover,  the  Helmholtz equation for $U$
amounts to \eqref{eq:hodgedirac} for $F_0=ikU$, 
$F_1= \nabla U$ and $F_2=F_3=0$, as well as for
$F_3=ikU$, 
$F_2= \nabla U$ and $F_0=F_1=0$.
A main point with the Dirac formalism is that it avoids 
divergence- and curl-free constraints, 
by complementing the divergence-free Maxwell vector
fields with the Helmholtz gradient vector fields.

The Dirac transmission problem DTP($k_-,k_+,\alpha,\beta,\gamma$)
which is fundamental in our formalism is 
\begin{equation}   \label{eq:Diractransmissionpr}
\begin{cases}
   F^+= M (F^0+ F^-), & x\in\Gamma, \\
    \dirac F^+=ik_+ F^+, &x\in\Omega_+,\\
    \dirac F^-=ik_- F^-, &x\in\Omega_-,\\
    (x/|x|-1)F^-=o(|x|^{-1}e^{\im(k_-) |x|}), & x\to\infty.
   \end{cases}
\end{equation} 
The Dirac derivative $\dirac F$ is the  left-hand  side in \eqref{eq:hodgedirac},
and replacing $\nabla$ by the vector $x$ in this matrix yields the Clifford product
$xF$ appearing in the Dirac radiation condition.
(For a short explanation of the underlying multivector formalism we refer  to~\cite[Sec.~3]{HelsRose20},  and for the long explanation we refer 
to~\cite{RosenGMA19}.)
On $\Gamma$, we write Dirac fields $F$ as
\begin{equation}  \label{eq:dirfields}
   F= \begin{bmatrix}
      F_0 & \nu\cdot F_2 & \mv F_{2T} & F_3 & \nu\cdot F_{1} &  \mv F_{1T}
   \end{bmatrix},
\end{equation} 
with tangential fields $\mv F_{jT}=(\tau\cdot F_j)\tau+ (\theta\cdot F_j)\theta$,
$j=1,2$.
The jump matrix $M$ is the diagonal matrix
\begin{equation}  
  M= \diag\begin{bmatrix} \hat k/(\alpha\beta) & 1/\hat k  & 
  \mv {\hat k/\alpha}  & 1/\gamma & 1/\alpha 
   & \mv{1}  \end{bmatrix},
\end{equation} 
when acting on $F$ in \eqref{eq:dirfields}. 
This special structure of $M$ ensures that the DTP decouples in a certain way
into one MTP with parameter $\alpha$ and two HTPs with parameters $\beta$ and $\gamma$. See~\cite[Prop.~8.4]{HelsRose20}.

For a given wavenumber $k$, the Dirac equation $\dirac  F= ikF$ comes with a Cauchy reproducing formula for $F$, 
similar to the classical one for analytic functions and $k=0$. 
When acting on $F$, written as in \eqref{eq:dirfields},  the singular Cauchy integral on $\Gamma$ is

\begin{equation}   \label{eq:EkCauchyinte}
E_k=
\begin{bmatrix}
-K^{\nu'}_k  & 0 & \mv K_{1,3:4} &
 0 & S_k & \mv{ 0}  \\
 K_{2,1} & -K^\nu_k & 
  \mv K_{2,3:4} & S_{2,5} & 0 & \mv S_{2,7:8} \\
 \mv K_{3:4,1} & \mv K_{3:4,2} & 
-\mv M_k^* & \mv S_{3:4,5} & \mv{ 0} & \mv S_{3:4,7:8}\\
  0 & S_k & \mv{ 0} & -K^{\nu'}_k & 0 & \mv K_{5,7:8} \\
 S_{6,1} & 0 & \mv S_{6,3:4} & K_{6,5} &
 -K^\nu_k & \mv K_{6,7:8} \\
 \mv S_{7:8,1} & \mv 0 & \mv S_{7:8,3:4} & \mv K_{7:8,5} & 
\mv K_{7:8,6} & -\mv{M}_k^*
\end{bmatrix}.
\end{equation}
In this paper, we enumerate the 
scalar components of Dirac fields \eqref{eq:dirfields} by $1,\ldots, 8$, and use 
boldface vector notation for
 the tangential vectors
parts 3:4 and 7:8. 
We denote the causal fundamental 
solution to the Helmholtz equation, normalized with a factor of $-2$, by
\begin{equation}
   \Phi_k(x)= \frac{e^{ik|x|}}{2\pi |x|}.
\end{equation} 
The operators appearing along the diagonal are the acoustic
double layer
\begin{equation}   \label{eq:acousticNP}
  K^{\nu'}_k f(x) = 
  \pv\int_\Gamma \nabla\Phi_k(y-x) \cdot \nu(y) f(y) d\Gamma(y)
\end{equation}
with real adjoint $-K^\nu_k f$, where $\nu(y)$ is replaced by $\nu(x)$,
and the real adjoint $\mv M^*_k$
of the acoustic magnetic dipole operator
\begin{equation}  \label{eq:magndipoleop}
  \mv M_k f(x) = 
  \nu(x)\times \pv\int_\Gamma \nabla\Phi_k(y-x) \times f(y) d\Gamma(y).
\end{equation}
Also appearing is the acoustic single layer
\begin{equation}
  S_k f(x) = 
 ik\int_\Gamma \Phi_k(y-x) f(y) d\Gamma(y),
\end{equation}
with the scaling factor $ik$.
The remaining operators $K$ and $S$ include various products of the
frame vectors $\{\nu,\tau,\theta\}$ and are detailed in~\cite[Eq.~(27)]{HelsKarlRos20}.

A fundamental algebraic property of $E_k$ is that for each 
 wavenumber  $k\in \C$, we have 
$E_k^2=I$. Using the associated Hardy projections
\begin{equation}
  E_k^\pm= \tfrac 12(I\pm E_k),
\end{equation} 
we can express DTP($k_-,k_+,\alpha,\beta,\gamma$) compactly as 
\begin{equation}   \label{eq:DTPwithprojs}
\begin{cases}
   F^+= M (F^0+ F^-), \\
   E_{k^+}^- F^+=0, \\
   E_{k^-}^+ F^-=0,
   \end{cases}
\end{equation} 
for $F^\pm$ and $F^0$ on $\Gamma$.
The condition $E_{k^+}^- F^+=0$ is equivalent to $F^+$ belonging
to the range of the projection operator $E_{k^+}^+$, that is, $F^+$ is the trace
of a solution to $\dirac F^+=ik_+ F^+$ in $\Omega^+$.
Similarly $E_{k^-}^+ F^-=0$ encodes that $F^-$ is the trace of an
exterior Dirac solution with  wavenumber  $k_-$, which satisfies the Dirac
radiation condition.

\section{The Dirac integral equation}   \label{sec:chooseparam}

In~\cite{HelsRose20}, BIEs for DTP($k_-,k_+,\alpha,\beta,\gamma$)
were derived as follows.
We make a field representation $F^+= rE_{k_+}^+(M'P'h)$, 
$F^-= -E_{k_-}^-(P' h)$ and multiply the 
jump relation $F^+= M (F^-+ F^0)$ by $P$.
In matrix notation, this amounts to the BIE
\begin{equation}   \label{eq:DiracBIEfactors}
  P \begin{bmatrix} E_{k_+}^+ & -M E_{k_-}^-\end{bmatrix}
  \begin{bmatrix} r & 0 \\ 0 & 1 \end{bmatrix}
  \begin{bmatrix} E_{k_+}^+M' \\- E_{k_-}^- \end{bmatrix}P'h=
  PMf^0,
\end{equation}  
for the density $h$ with $8$ scalar components.
We write $f^0= F^0|_\Gamma$ for the incident field on $\Gamma$.
The preconditioning matrices $P, P'$ will be chosen as constant diagonal
matrices, and the scaling parameter $r$ is $r=1/\hat k$ in all our Dirac BIEs.

The well posedness of DTP($k_-,k_+,\alpha,\beta,\gamma$) is equivalent
to invertibility of $\begin{bmatrix} E_{k_+}^+ & -M E_{k_-}^-\end{bmatrix}$, 
as a map from the direct sum of the ranges of $E_{k_+}^+$ and
$E_{k_-}^-$.
Consider also an auxiliary DTP($k_+,k_-,\alpha',\beta',\gamma'$),
with the  wavenumbers  swapped, with an auxiliary Maxwell jump parameter $\alpha'$ and with two 
auxiliary Helmholtz jump parameters $\beta', \gamma'$.
The duality result from~\cite[Prop.~8.5]{HelsRose20} shows that
well posedness of DTP($k_+,k_-,\alpha',\beta',\gamma'$) is equivalent 
to invertibility of
$\begin{bmatrix} E_{k_+}^+M' \\- E_{k_-}^- \end{bmatrix}$,
as a map onto  the direct sum of the ranges of $E_{k_+}^+$ and
$E_{k_-}^-$, with
\begin{equation}
  M'= \diag\begin{bmatrix} 1/\alpha' & 1/\gamma' & \mv{1}  & \hat k & 1/(\hat k\alpha'\beta') &
  \mv{1/(\alpha'\hat k)}   \end{bmatrix}.
\end{equation}

The resulting Dirac BIE \eqref{eq:DiracBIEfactors} has $12$ free parameters
\begin{equation}   \label{eq:12param}
  r, \beta, \gamma, \alpha', \beta', \gamma'\quad\text{and}\quad
  P'=\diag\begin{bmatrix} p'_1 & p'_2 & \mv p'_{3:4} & p'_5 & p'_6 & \mv p'_{7:8}   \end{bmatrix}
\end{equation}
to be chosen.
We recall that $\alpha=\hat k^2$ for the non-magnetic Maxwell transmission problem 
that we want to solve.
For $r$ and $P'$, any non-zero and finite complex numbers are allowed,
although we have always used $r=1/\hat k$ so far.
Given $P'$, we choose $P$ so that 
\begin{equation}
  P(rM'+M)P'=I,
\end{equation} 
and set $N= PM$ and $N'=rM'P'$.
This turns the Dirac BIE \eqref{eq:DiracBIEfactors} into a second kind
integral equation
\begin{equation}  \label{eq:DiracBIEmain}
  h+(P E_{k_+} N'- N E_{k_-} P') h= 2N f^0,
\end{equation}
where we have used that $(E_k^\pm)^2= E_k^\pm= \tfrac 12(I\pm E_k)$.
The operator to invert on $\Gamma$ is $I+G$, where $G$ denotes the 
singular integral operator 
\begin{equation}  \label{eq:G}
  P E_{k_+} N'- N E_{k_-} P'= P(E_{k_+}(rM')-ME_{k_-})P',
\end{equation}
containing $30$ scalar integral operators
of double layer type, and $20$ scalar integral operators
of single layer type, according to \eqref{eq:EkCauchyinte}.
For a given choice of parameters, the operator $G$ is computed 
using~\cite[Eqs.~(131),~(132)]{HelsRose20}.

From the density $h$, obtained by solving \eqref{eq:DiracBIEmain},
we compute the traces on $\Gamma$ of
the Dirac fields $F^\pm$ as
\begin{equation}   \label{eq:projdens}
  F^+|_\Gamma= E_{k_+}^+ N' h\quad\text{and}\quad
   F^-|_\Gamma= -E_{k_-}^- P' h.
\end{equation} 
The fields $F^\pm$ in $\Omega^\pm$ are computed by instead using the corresponding
non-singular (replacing $x\in\Gamma$ by $x\in\Omega^\pm$) versions
of the Cauchy integrals $E_{k_\pm}^\pm$ in \eqref{eq:projdens}.
Discarding the two auxiliary Helmholtz components in $F^\pm$ (which are
$0$ for Maxwell data $f^0$), 
and recalling \eqref{eq:Hplus}, this yields $E^\pm$ and $H^\pm$. 
These field formulas are detailed in~\cite[Eqs.~(28)--(31)]{HelsKarlRos20}.
Note that in computing $H^+$, the factor $\hat k$ from \eqref{eq:Hplus}
cancels the factor $r=1/\hat k$ in $N'$ from \eqref{eq:projdens},
whereas this factor $1/\hat k$ remains in $E^+$. 
Note also that the minus sign in the second equation in \eqref{eq:projdens}
is contained in the field equations~\cite[Eqs.~(28), (30)]{HelsKarlRos20}.

\subsection{Dirac (A)}   
The Dirac parameters 
\begin{equation} \label{eq:betaDir}
\begin{split}
  &\begin{bmatrix} r & \beta & \gamma & \alpha' & \beta' & \gamma'
  \end{bmatrix} = 
  \begin{bmatrix} \frac 1{\hat k} &  \xi & a & 
  \frac 1{\hat k} & \frac 1{\hat k} &  \rev a \end{bmatrix},
   \\
  P' &=\diag\begin{bmatrix} 
1 & \hat k^{1/2}(1+a)^{-1/2}
  & \mv{\hat k^{1/2}} & 1 & 1 
  & \mv{\frac{\hat k}{{\hat k}+1}} 
  \end{bmatrix},
\end{split}
\end{equation}
were proposed in~\cite[Thm.~2.3]{HelsRose20},
where $a=\hat k/|\hat k|$.
Here the overline symbol denotes the complex conjugate and
\begin{equation}  \label{eq:xi}
\xi=1+i\delta \arg(\hat k)
\end{equation}
 is a tuning factor, set to $\xi=1$ in~\cite{HelsRose20} and further discussed below. 
The relations above then give
\begin{equation}  
\begin{split}
  P &= \begin{bmatrix}\frac \xi{{\hat k}^{-1}+\xi} &  \hat k^{1/2}(1+a)^{-1/2}
   & \mv{\frac {\hat k^{1/2}}2 }& \frac 1{1+\rev a} & \frac 1{1+\hat k^{-2}} & \mv 1 
   \end{bmatrix},  \\
  N &=\begin{bmatrix}  \frac 1{1+\xi{\hat k}} &  \frac 1{\hat k^{1/2}}(1+a)^{-1/2}
  &   \mv{\frac 1{2\hat k^{1/2}} } & \frac 1{1+a}  & \frac {1}{1+{\hat k^2}} & \mv 1 
   \end{bmatrix},  \\
  N' &= \begin{bmatrix}1 &   \frac 1{|\hat k|^{1/2}}(1+\rev a)^{-1/2}
   & \mv{\frac 1{\hat k^{1/2}}} & 1  & 1 & \mv{ \frac 1{1+\hat k}}
     \end{bmatrix}.
\end{split}
\end{equation} 
The main operator $G$ from \eqref{eq:G} 
is in general not  compact, not even 
on smooth domains. However, we have that
$G$ is nilpotent modulo compact operators on smooth domains. 
This property is important for the efficiency of GMRES, since it 
implies that the only accumulation point for the spectrum of 
$G$ is zero. 
The above choices of $r,\beta=1, \alpha', \beta'$ ensure this
through cancellations in the (1:2,3:4) and (7:8,5:6) blocks.
The choices of $\gamma, \gamma'$ were made only to avoid false
eigenwavenumbers by chosing suitable complex arguments for them.

In~\cite[Sec.~5.1]{HelsKarlRos20}, we chose $\delta= 0.2/\pi$ in \eqref{eq:xi}. 
This turns the complex argument of $\beta$ slightly towards the argument of $\hat k$, which avoids false
   eigenwavenumbers in plasmonic  scattering. The change  is still small enough for $G$ to be close to a nilpotent
        operator modulo compact operators on smooth $\Gamma$, so that iterative solvers converge rapidly. 
See Section~\ref{sec:chooseaug}.
We refer to this Dirac BIE as Dirac (A), which performs well on any
Lipschitz surface, of any genus, as long as $|\hat k|$ is bounded away
from $0$ and $\infty$.
Throughout~\cite{HelsRose20, HelsKarlRos20}, we only  consider
$|\hat k|$ which are not too large or small, and 
the choice for $P'$ is then less  important.
In the present paper, we allow  $\hat k\to\infty$, and then
the parameters need to be chosen with more care,
since Dirac (A) is no longer performing well.

\section{Tuning and augmentation}  \label{sec:chooseaug}

We describe in this section some ideas for designing 
efficient BIEs in general, and Dirac BIEs in particular.
A first step is to obtain Fredholm operators both (a) for the system
to be solved and (b) for the field representation,  by
tuning the free parameters in the BIE.
Recall that an operator is Fredholm if 
its null space is finite-dimensional and its range is closed
and has finite codimension. 
An operator  that  fails to be Fredholm, for example a differential 
or hypersingular operator, may lead to  ``dense-mesh
breakdown'' in the computation,  in  the 
terminology of~\cite{ValdesETAL:11}.
Even if the operators are Fredholm for each $k_->0$, but not
uniformly as $k_-\to 0$, there may be a low-frequency breakdown 
in the computation. 

The top three considerations in tuning Dirac BIEs are the following.

\begin{itemize}
\item
The singular integral operator $G$ should be close  in norm  
to a nilpotent operator
modulo compact operators on smooth $\Gamma$,
to work well in an iterative solver.
Our experience is that the choice $r=1/\hat k$ is necessary for this.

\item
The complex arguments of $\beta/\hat k$, $\gamma/\hat k$,
$\alpha'\hat k$, $\beta'\hat k$ and $\gamma'\hat k$ 
should be
\begin{equation}  \label{eq:falsesuffi}
  \le \pi-\arg(k_+/i)-\arg(k_-/i)
\end{equation}
to avoid false eigenwavenumbers.
When $k_\pm$ are on the real or imaginary
axis, a strict inequality is required, and sometimes $\arg(-z)$
may replace $\arg(z)$, for $z$ being one of the five complex numbers above.
We refer to \cite[Props. 8.4-5 and Defn. 2.1]{HelsRose20}
for details.
False essential spectrum may appear when one of 
$\beta, \gamma, \alpha', \beta', \gamma'$ is in a compact
subset of the negative real axis $(-\infty,0)$.
See Section~\ref{sec:Maxwess}.

The role of the tuning factor $\xi$ from \eqref{eq:xi} in the Dirac BIEs is to 
slightly adjust parameters, not only $\beta$, but also $\alpha'$
and $\gamma'$ in the Dirac BIEs presented below,
to obtain a strict inequality in \eqref{eq:falsesuffi}.
This typically ensures  
that false eigenwavenumbers are avoided in plasmonic scattering, that is 
when $(\arg(k_-), \arg(k_+))=(0,\pi/2)$ or $(\pi/2,0)$.
\item
The coefficients in $G$, depending on $P, P', N, N'$, should be
uniformly bounded in the set of $k_\pm$ considered.
For this, we note that it is sufficient but not necessary 
to have $P, P', N, N'$ bounded.
That it is needed to have $N$ bounded is clear from the  right-hand  side
in \eqref{eq:DiracBIEmain}.
Most importantly, $N'$ and $P'$ should be chosen so that
\eqref{eq:projdens} computes the fields at the correct scale.
The generic scale for the fields in the eddy current regime 
is discussed in Section~\ref{sec:physics} below. 
\end{itemize}

A second step in the design of efficient BIEs is to remove remaining 
finite-dimensional null spaces in the Fredholm maps.
In an abstract setting, the typical situation is that finite-dimensional 
null spaces open up as a parameter $\lambda\to 0$,
depending on the topology of $\Gamma$, causing a ``topological low-frequency 
breakdown''. 
See \cite{AndiulliETAL09, EpsteinETAL13, EpsGreNei19} 
for examples of this generic problem for BIEs.
We remove such null spaces through augmentation. 
Recall from \eqref{eq:abstractIER} that the typical construction of a BIE is
to compose jump relations $g= B_\lambda F$ and a field representation $F= A_\lambda h$
to obtain a linear system  
\begin{equation}  \label{eq:ABfact}
g= (I+G_\lambda)h= B_\lambda A_\lambda h.
\end{equation} 
Here both $g$ and $h$ belong to a suitable space of functions on
   $\Gamma$, but the domain of $B_\lambda$ and the range of $A_\lambda$
   consist of a space of fields $F=F^\pm$ in $\Omega_\pm$ satisfying the
   differential equation, or equivalently a corresponding space of traces
   $F^\pm|_\Gamma$ on $\Gamma$. We assume that both $A_\lambda$ and
   $B_\lambda$ are Fredholm maps with index $0$. In general both maps
   $A_\lambda$ and $B_\lambda$ may require augmentation. We refer to
   augmentation of the right factor, the field representation, as (R)
   augmentation, and to augmentation of the left factor, the jump
   relations, as (L) augmentation.

   The logic behind the two types of augmentations is quite different.
   We therefore discuss them separately, starting with two elementary but
   illustrative examples of augmentation of BIEs for Helmholtz boundary
   value problems, with no aim for completeness or full proofs. It is
   clear from these examples that (R) augmentation is required when we
   have null densities $h+G_0 h=0$ corresponding to zero fields $F=A_0
   h=0$. Otherwise (L) augmentation is required in order to obtain an
   invertible system.

\begin{ex}[Exterior Dirichlet problem]
  Consider the Dirichlet problem $u|_\Gamma=g$ for $\Delta u+k^2 u=0$ in $\Omega_-$,
with standard Sommerfeld radiation condition at $\infty$, 
for $k$ in a neighbourhood of $0$.
This has a unique solution $u$ for all $k$ in a neighbourhood of $0$,
including $0$, so no (L) augmentation is needed.
Using the standard double layer potential representation of $u$ however, leads to
a BIE $h+ K^{\nu'}_k h=2g$, which at $k=0$ has a null space spanned by constant functions
$h$. We resolve this problem by using an  (R) augmented  field representation
\begin{equation}   \label{eq:extdirRaug}
  u(x)= \int_\Gamma \nabla\Phi_k(y-x)\cdot \nu(y) h(y) d\Gamma(y)+ \Phi_k(x-p)\int_\Gamma h d\Gamma(y), \qquad x\in\Omega_-,
\end{equation}
with some fixed $p\in \Omega_+$. 
In dimension two, the fundamental solution $\Phi_k(x-p)=(i/2)H^{(1)}_0(k|x-p|)$
uses the Hankel function and
needs to be divided by $\log(k)$ in the second term in \eqref{eq:extdirRaug}.  
This leads to the  (R) augmented  BIE $h+ K^{\nu'}_k h+ b_k(ch)=2g$
with the  finite-rank  operator $b_kc$ added, where $b_k=\Phi_k(\cdot-p)|_\Gamma$ and
$ch= \int_\Gamma h d\Gamma$.
This is stably solvable for $h$, for all $k$ near $0$, and the field $u$ is computed from the
augmented representation \eqref{eq:extdirRaug}. 
This can be seen as a
rank-1 version of the combined field integral 
equation~\cite[Eq.~(3.25)]{ColtonKress:inverseAc}, which only
removes the false eigenwavenumber at $k=0$.
It can also be seen as a generalization to $k\ne 0$ of the standard
treatment for $k=0$ in \cite[p.~345]{GuentherLee:96}.
\end{ex}

In an abstract setting, (R) augmentation can be described as follows.
We consider a parameter $\lambda\to 0$, and for each fixed $\lambda \ne 0$ we
compose jump relations $B_\lambda$ and field representations $A_\lambda$ as 
 in \eqref{eq:ABfact}  to obtain a system
\begin{equation}   \label{eq:Glambdaeq}
  h+G_\lambda h= g.
\end{equation}
We assume that free parameters have been tuned so that $I+G_\lambda$ is invertible
for each $\lambda \ne 0$, but that it is  merely  a Fredholm map at $\lambda=0$.
Assume for simplicity that the null space $\nul(I+G_0)$ is one-dimensional.
If $B_0$ is invertible and the null space for $I+G_0$ comes from the field representation $A_0$,
then we identify a functional $c$  that  is non-zero on 
$\nul(A_0)=\nul(I+G_0)$ and a field/solution $F_\lambda$
for each parameter $\lambda$, such that $F_\lambda\to F_0$ where $F_0$ does not belong 
to the range $\ran(A_0)$.
Let $b_\lambda= B_\lambda F_\lambda$ be the corresponding boundary datum.
The obtained  (R) augmented  BIE uses the field representation 
$F= A_\lambda h+ F_\lambda(ch)$, where the density $h$ now solves
the  (R) augmented  system $(I+ G_\lambda+ b_\lambda c) h=g$.

\begin{ex}[Interior Neumann problem]    \label{ex:modelinhomoL}
  Consider the Neumann problem $\pd_\nu u|_\Gamma=g$ for $\Delta u+k^2 u=0$ in $\Omega_+$,
for $k$ in a neighbourhood of $0$.
We have a unique solution $u$, except at $k=0$, and an (L) augmentation is required.
To note is that 
\begin{equation}
\int_{\Gamma} g d\Gamma= -k^2 \int_{\Omega_+} u dx, 
\end{equation}
which forces $\int_{\Omega_+} udx=0$
if $k\ne 0$ and $\int_{\Gamma} gd\Gamma=0$.
The standard single layer potential
$u(x)= \int_\Gamma \Phi_k(x-y) h(y) d\Gamma(y)$, $x\in\Omega_+$,
gives a representation of all solutions $u$ in $\Omega_+$ and needs no  (R) augmentation, 
but leads to the BIE $h-K^\nu_k h=2g$, which  for $k=0$  has a one-dimensional null space.
For $k\ne 0$, define the functional
\begin{equation}  \label{eq:Neumannc}
c_kh= \int_{\Omega_+} u dx= -\frac 1{k^2}\int_{\Gamma} g d\Gamma=
 -\frac 1{2k^2}\int_\Gamma (h-K_k^\nu h)1 d\Gamma=   
 \int_{\Gamma} hw_k d\Gamma, 
\end{equation}
with weight function
$w_k= (K_0^{\nu'}1- K_k^{\nu'}1)/(2k^2)$. 
The last equality in \eqref{eq:Neumannc} follows from duality and
$K_0^{\nu'}1=-1$.  
By Taylor expansion of $\nabla \Phi_k$, this computation of $w_k$ can be stabilized, and for $k\ne 0$ 
our BIE is seen to be equivalent to 
\begin{equation}
  h-K^\nu_k h+ b(c_k h)= 2g-k^{-2}b\int_{\Gamma} g d\Gamma.
\end{equation}
Choosing the constant function $b=1$ on $\Gamma$, the operator $I-K^\nu_k+bc_k$ is well-conditioned for all $k$ near $0$, 
and the previous low-frequency breakdown has
moved to the simpler computation of the second term on the 
 right-hand  side.
\end{ex}

In an abstract setting, (L) augmentation can be described as follows.
Consider $I+G_\lambda= B_\lambda A_\lambda$ as $\lambda\to 0$ as above,
but assume now that the field representation $A_0$ is invertible but that $B_0$ 
has a one-dimensional null space.
We identify a scalar  equation
\begin{equation}  \label{eq:ceqd}
   c_\lambda h = d_\lambda g,
\end{equation} 
which follows from $h+G_\lambda h=g$ for $\lambda\ne 0$, but may fail at $\lambda=0$.
The equation \eqref{eq:ceqd} often appears by rescaling one scalar component of the jump relations so that $c_\lambda$ is normalized.
We assume that $c_\lambda\to c_0$ as $\lambda\to 0$, where  $c_0\ne 0$  on $\nul(I+G_0)$.
Typically  $d_\lambda g$  will not stay bounded as $\lambda\to 0$, unless we assume that data satisfy $d_\lambda g=0$,
in which case we speak of a {\em homogeneous} (L) augmentation.
Choosing an auxiliary function $b\notin \ran(A_0)= \ran(I+G_0)$, we obtain the (L)
augmented BIE with system
$(I+ G_\lambda +b c_\lambda)h= g+ b(d_\lambda g)$ and 
field representation as before.

For Dirac BIEs, the factorization $I+G_\lambda= B_\lambda A_\lambda$ is 
given by \eqref{eq:DiracBIEfactors}.
Following the ideas presented in this section, we obtain  the augmentations for our Dirac BIE
that  are stated in Sections~\ref{sec:finitecond} and \ref{sec:DirB}.
The  derivation  of these augmentations are found in~\ref{app:aug}. 
It should be noted that Dirac (A$\infty$) only requires a single (L) augmentation for the Dirichlet eigenfield, 
regardless  of the genus of $\Gamma$, 
whereas Dirac (B) for $\Gamma$ of genus $g\ge 0$ requires 
$1+g$ (R) augmentations and
$1+g$ (L) augmentations for the Dirichlet and Neumann eigenfields, and  an 
(L) augmentation of a Helmholtz eigenfield. 

\begin{rem}
  An alternative technique for augmenting a BIE in the form $(I+G_\lambda)h=g$ is to add
equations and unknowns. Consider 
\begin{equation}  \label{eq:rowcoladded}
\begin{bmatrix} I+G_\lambda & B_\lambda \\ C_\lambda & D_\lambda \end{bmatrix}
\begin{bmatrix} h \\ a \end{bmatrix}= \begin{bmatrix} g \\ d_\lambda g \end{bmatrix}.
\end{equation} 
With only one extra equation and unknown in \eqref{eq:rowcoladded}, and
with $D_\lambda= -1$, the elimination  of $a$ shows that 
\eqref{eq:rowcoladded} is equivalent to the 
additive (L) augmentation described above. 

Now \eqref{eq:rowcoladded} can be converted to Fredholm second kind block triangular form, which is
better suited for an iterative solver, using 
\begin{equation}  \label{eq:ggprecond}
  \begin{bmatrix} I & B_\lambda \\ C_\lambda & D_\lambda \end{bmatrix}
\end{equation}
as a left- or right preconditioner in \eqref{eq:rowcoladded}.
The inverse of \eqref{eq:ggprecond} can be efficiently
applied via the solution of a system involving the Schur complement of $I$.
 See~\cite[Sec.~4.1]{GreenbaumETAL:93}  for an example.
Nevertheless, we find that our additive
augmentations are easier to work with than this more traditional technique
involving extra equations and preconditioners.
\end{rem}

Augmentations have been frequently used in literature. Most of these
   concern the simpler problem of augmenting static Laplace or biharmonic
   problems, which corresponds to augmenting only $I+G_0$, that is,
   $I+G_\lambda$ at $\lambda=0$. Static homogeneous (L) augmentation
   appears in \cite[p. 257]{Mikhlin:64}. Static (R) augmentation appears
   in \cite[Eqs.~(11),(23)]{GreenbaumETAL:93}. Non-static inhomogeneous
   (L) augmentation appears in \cite[Rem.~1]{EpsteinETAL13}.

\section{Eddy current eigenfields}   \label{sec:physics}

Assume that the incident fields $E^0$ and $H^0$ have magnitude of order $1$ on $\Gamma$.
Then the transmitted magnetic field $H^+$ is of order $1$ since $\Omega_+$ is non-magnetic, and the scattered electric field $E^-$ is also of order $1$, assuming that we are in the eddy current regime.
The transmitted electric field $E^+$ satisfies
$\nabla\times E^+= ik_- H^+$, $\nabla\cdot E^+=0$ and 
\begin{equation}   \label{eq:jumpEnor}
\nu\cdot E^+=
\hat k^{-2} \nu\cdot (E^- + E^0),
\end{equation} 
which shows that $E^+$ is of order
$\max(k_-L, |\hat k^{-2}|)$ in the generic scattering situation.
From this we conclude that
the magnitude of the eddy current 
\begin{equation}
J= \sigma E^+= \re(k_+^2/(i\eta_0 k_-))E^+
\end{equation}
is of order
$\max(|k_+|^2 L, k_-)$, 
and the scattered magnetic field $H^-$ is of order $\max((|k_+| L)^2, k_-L)$.

   The incident field that we use in our numerical examples is a sum of
   the two lowest order axially symmetric spherical vector waves
   \cite[Eq.~(7)]{KarlssonKristensson:83}
   \begin{equation}
     E^0(x)=G(x)+k_-^{-1}\nabla\times G(x)\,,\qquad H^0(x)=-iE^0(x)\,.
   \label{eq:partialwave}
   \end{equation}
   Here $G(x)=\sqrt{3/(8\pi)}j_1(k_-|x|)\rho|x|^{-1}\mv\theta$ and $j_1$
   is the spherical Bessel function of order $1$.

\begin{figure}[t!]
\centering
\includegraphics[height=40mm]{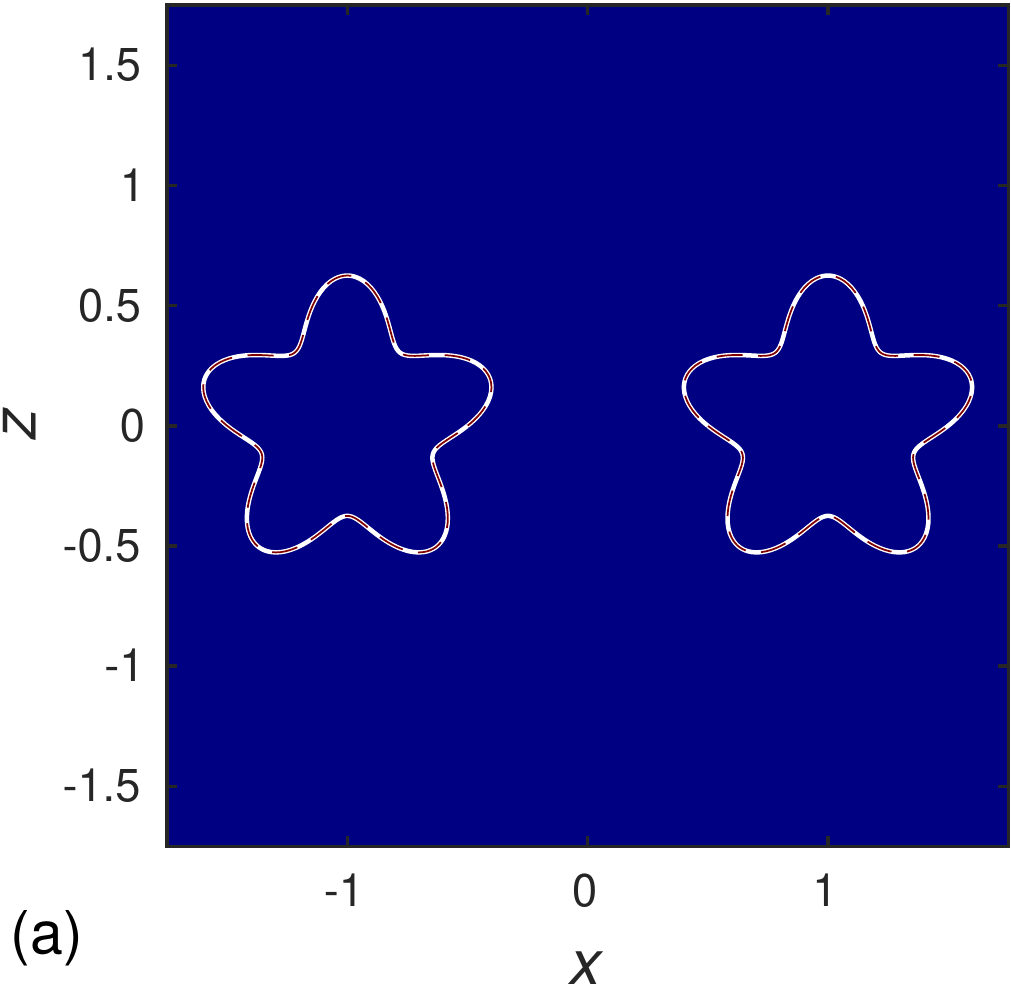}
\hspace{7mm}
\includegraphics[height=40mm]{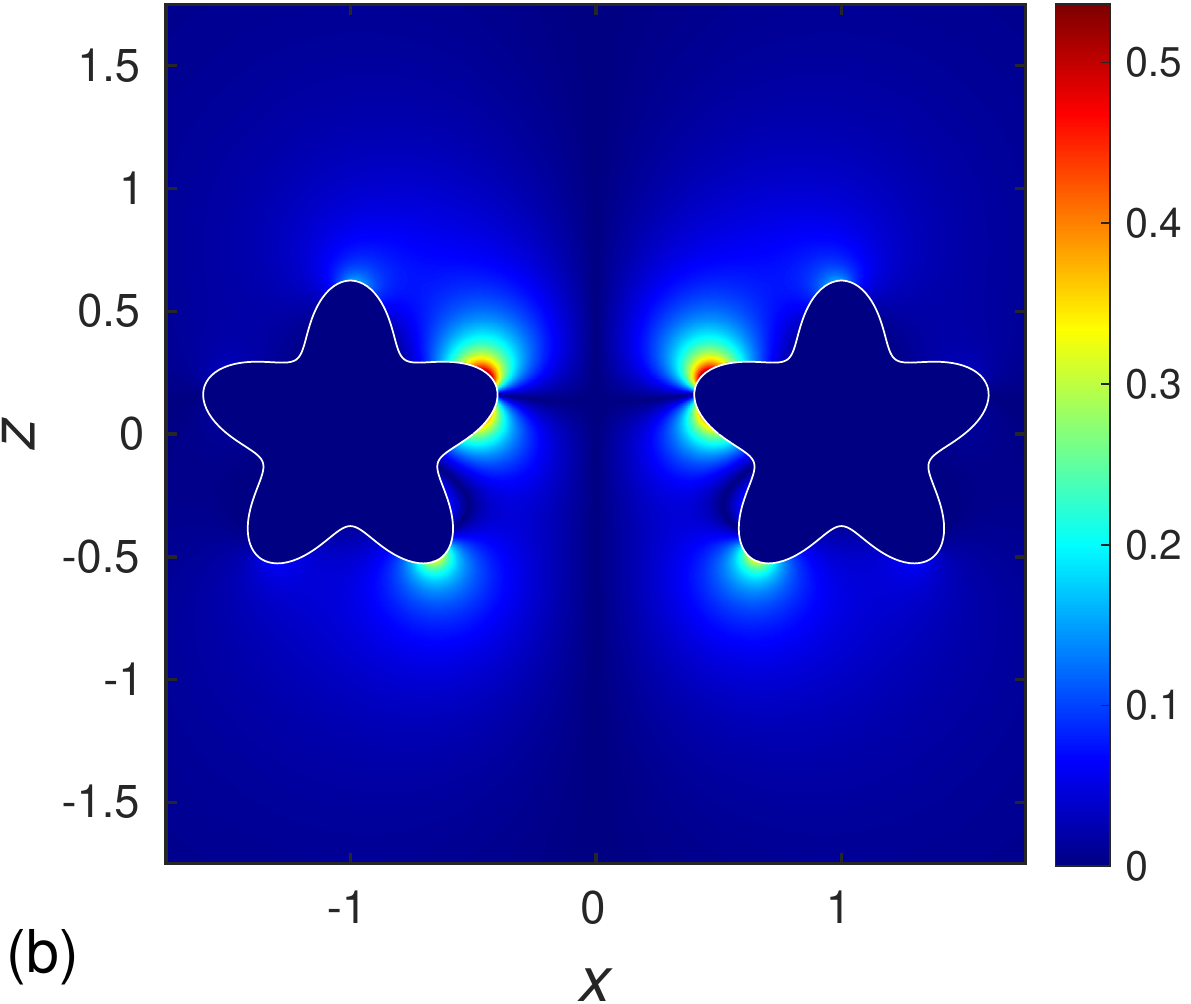}
\includegraphics[height=40mm]{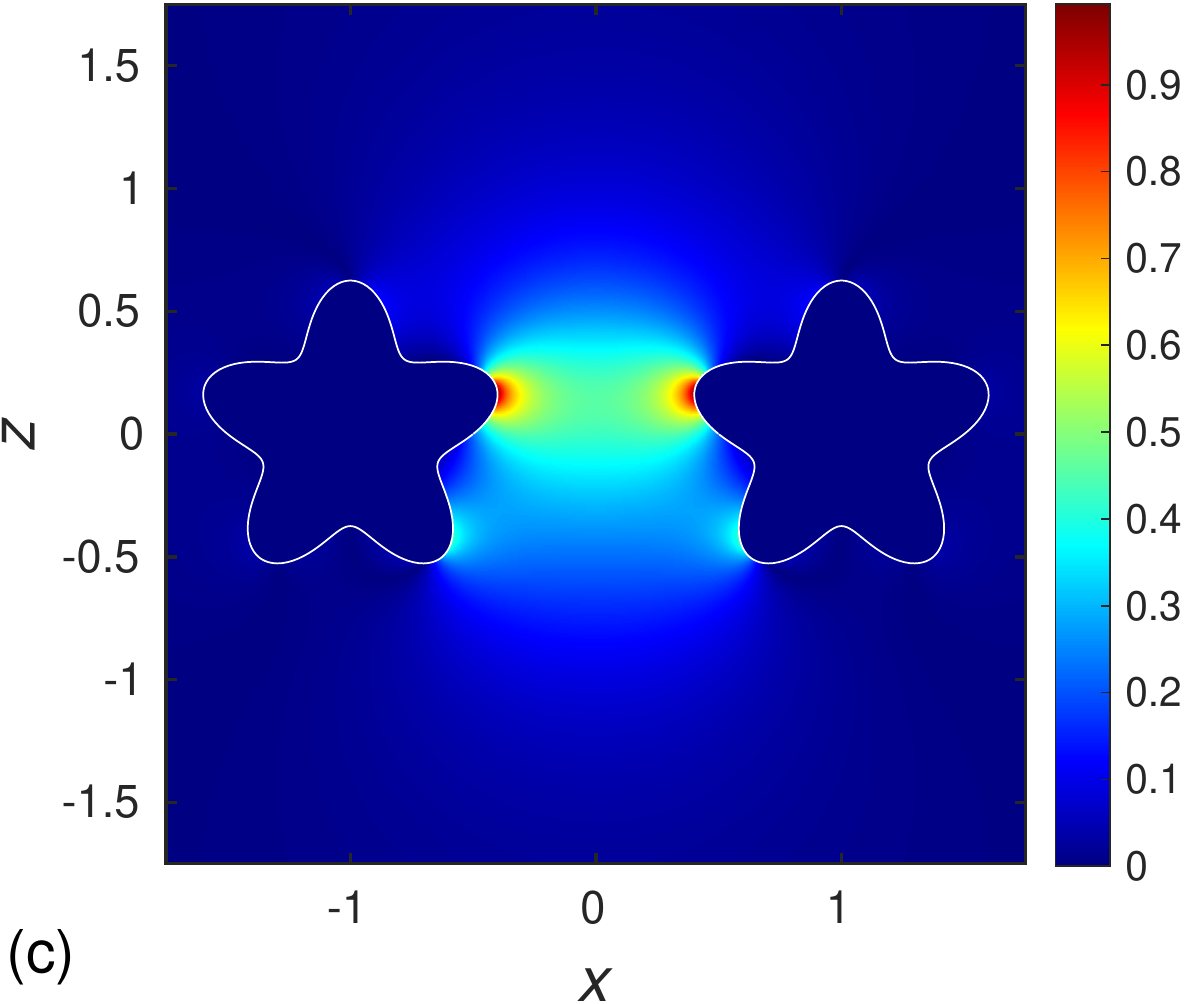}
\vspace{1mm}

\includegraphics[height=40mm]{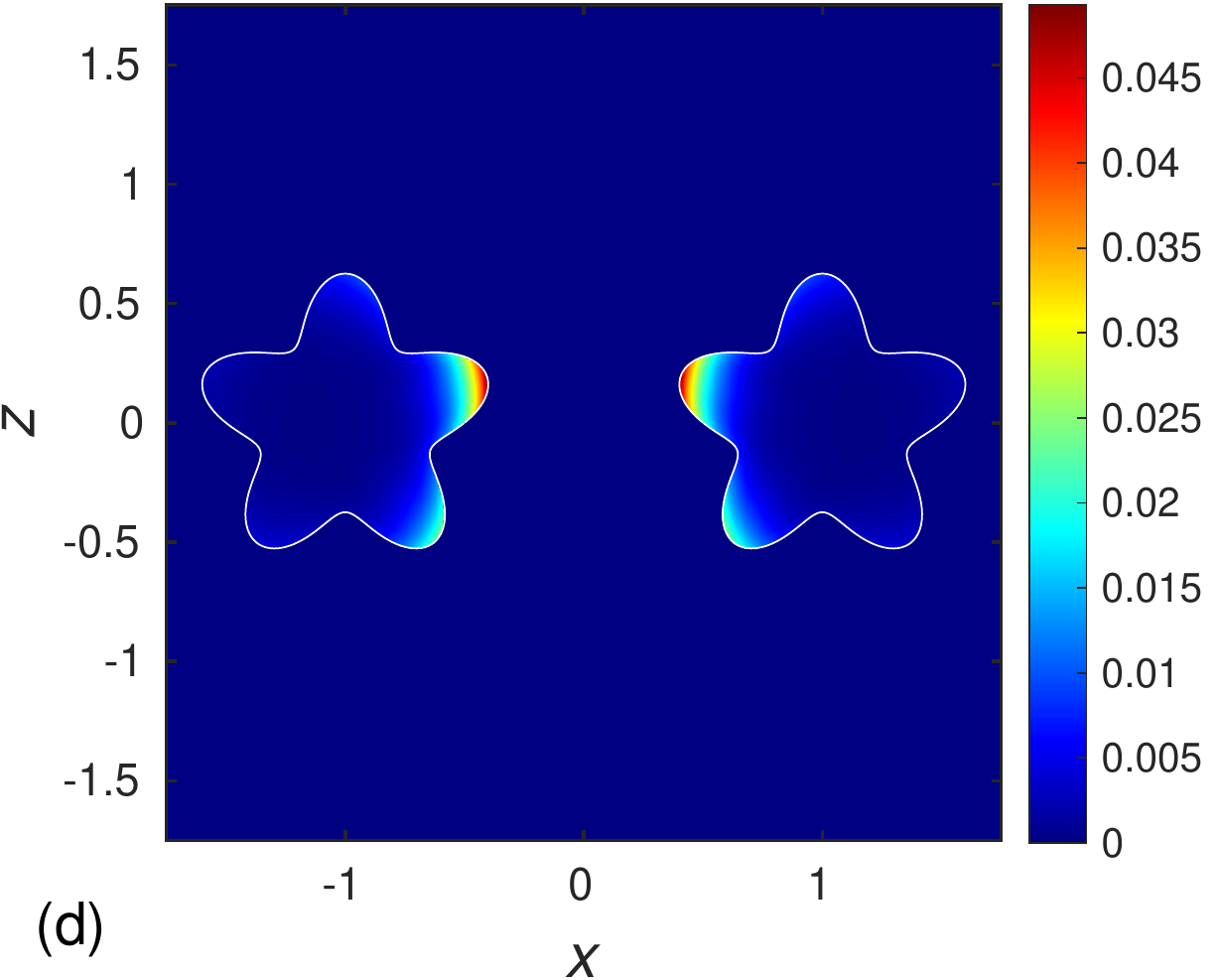}
\includegraphics[height=40mm]{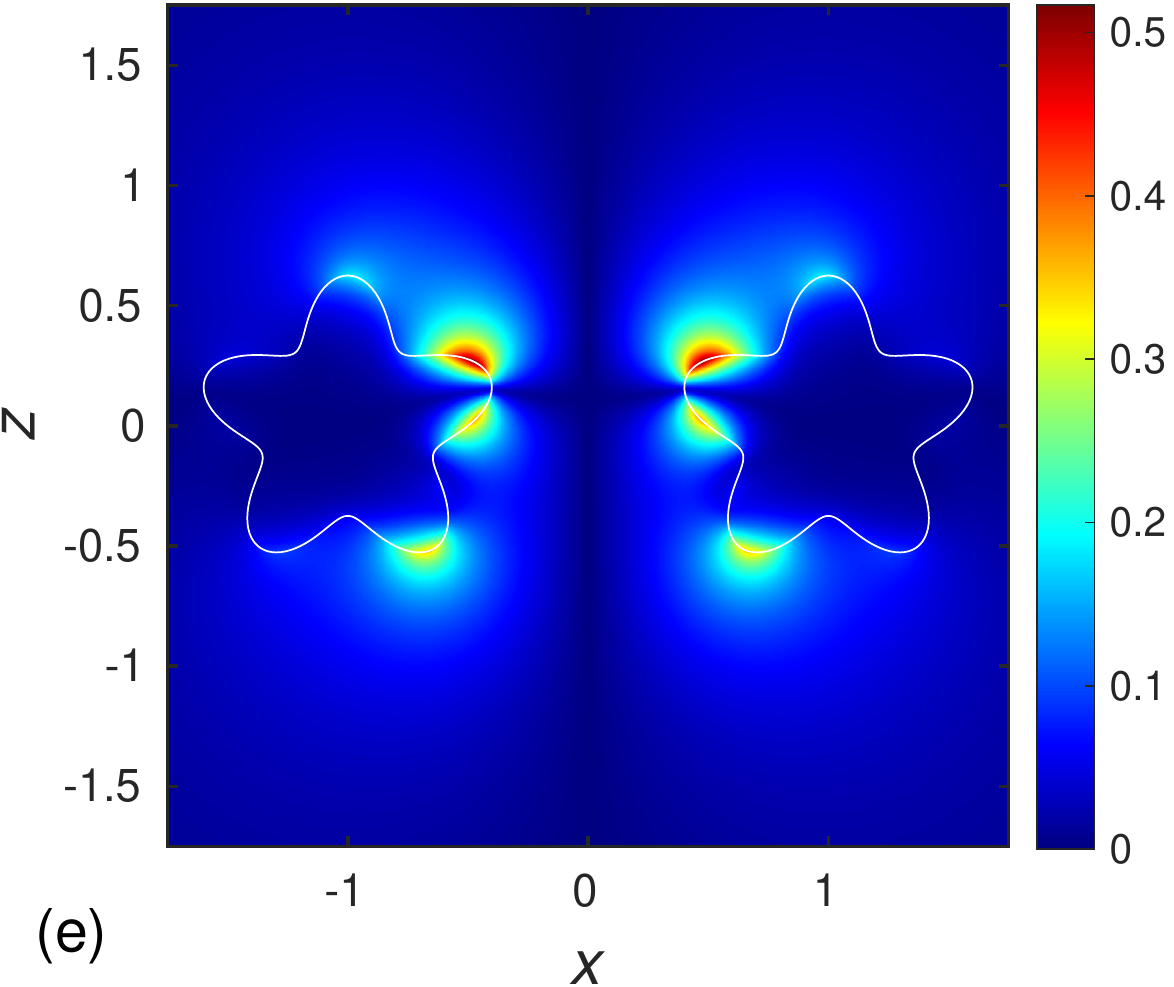}
\includegraphics[height=40mm]{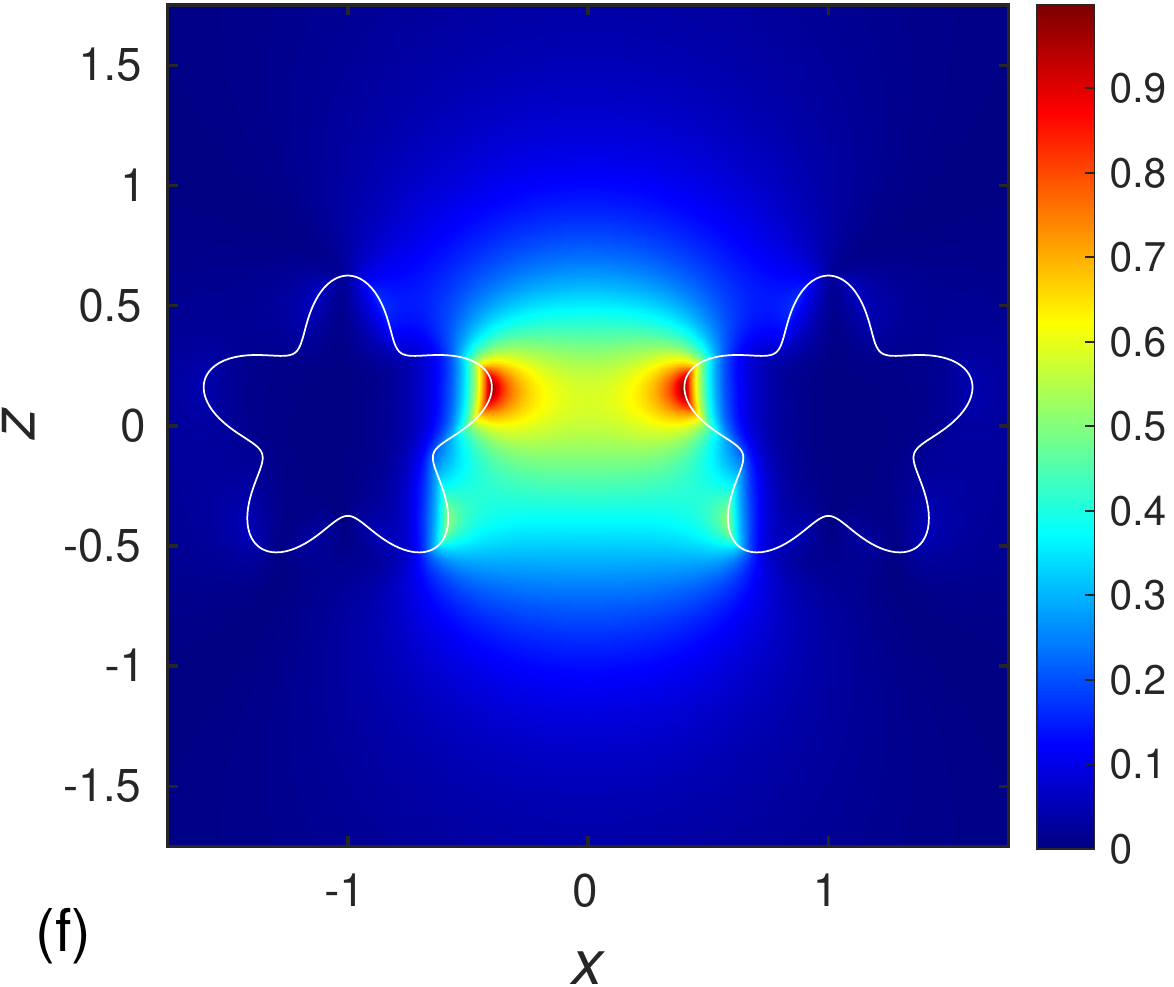}
\vspace{1mm}

\includegraphics[height=40mm]{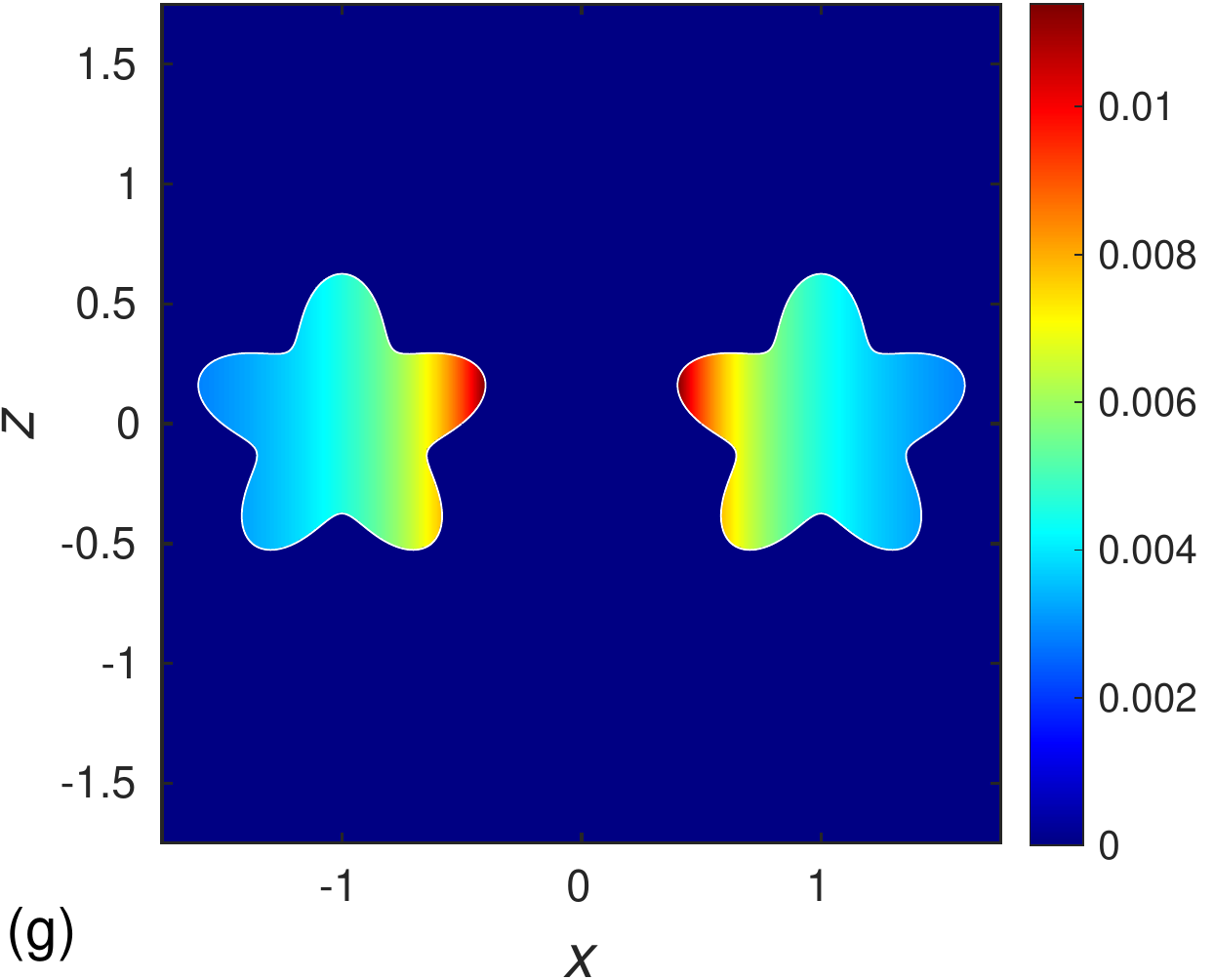}
\includegraphics[height=40mm]{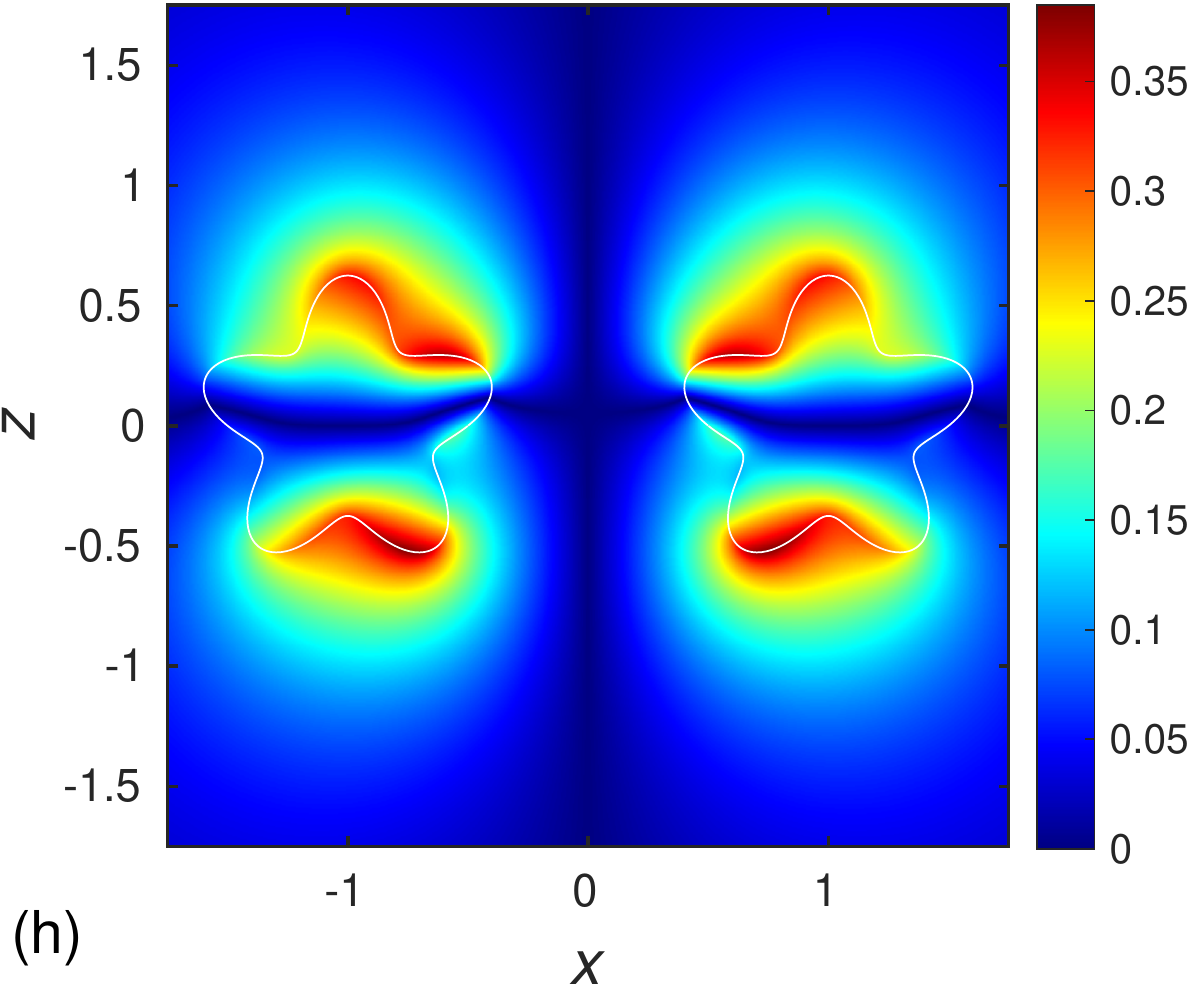}
\includegraphics[height=40mm]{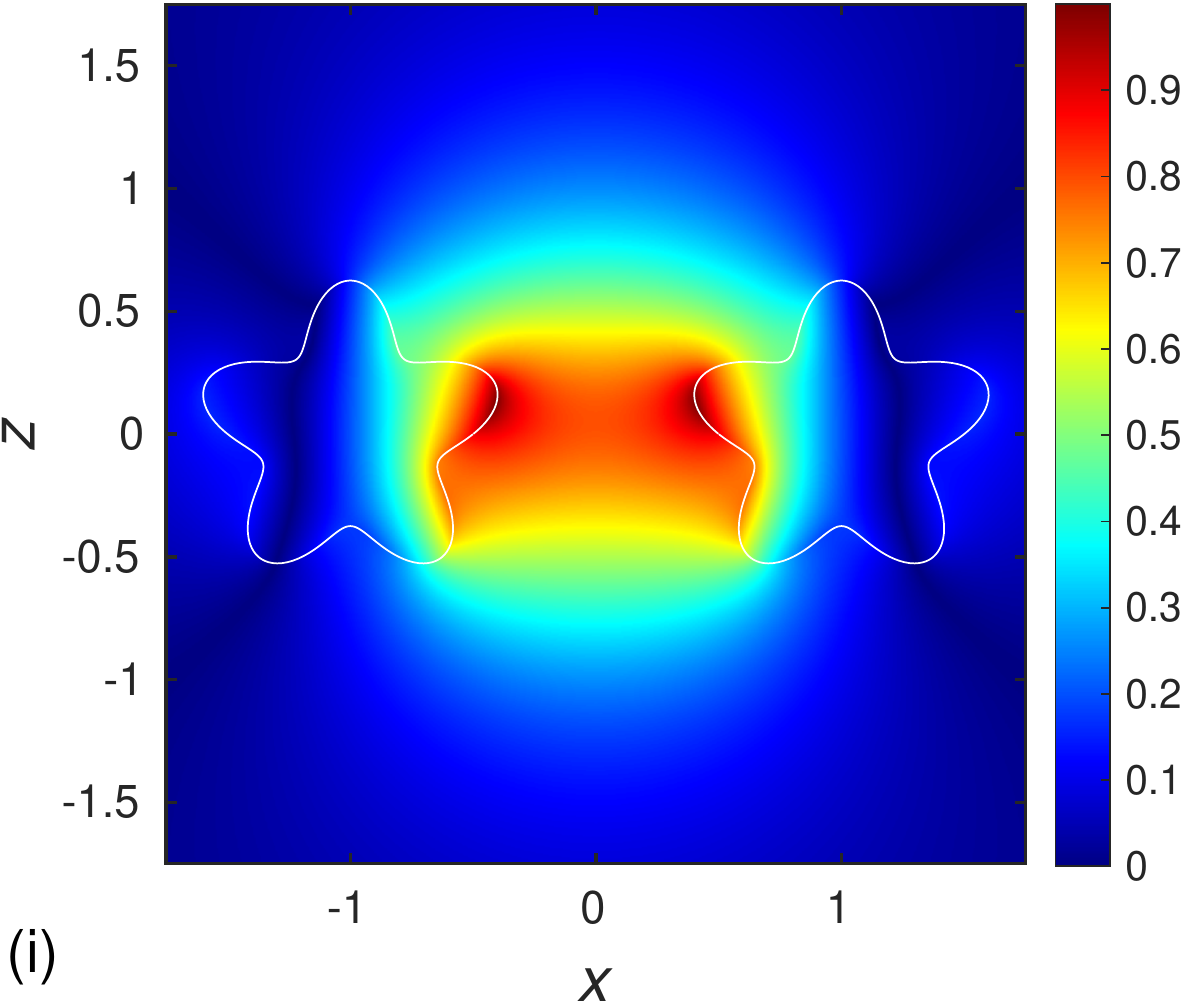}
\caption{\sf Neumann eigenfields of the ``starfish torus''~\eqref{eq:starfishtorus}, normalized so that 
$\max|H|=1$;
(a,d,g) eddy current $|J_\theta|$; 
(b,e,h) $|H_\rho|$; (c,f,i) $|H_z|$;
(a,b,c) eigenfield for a superconductor;
(d,e,f) borderline eddy current-PEC eigenfield at $k_+= 10(1+i)$; 
(g,h,i) eigenfield of an ordinary conductor.}
\label{fig:Neumanneig}
\end{figure}

We see two types of eigenfields appearing in the eddy current regime
as $k_-\to 0$,
which generalize the  classical  (exterior) Dirichlet and Neumann eigenfields
(or Dirichlet and Neumann vector fields in the terminology of \cite{ColtonKress92})
in PEC scattering, and where the magnitude of the fields can 
differ drastically from those described above. 
As discussed below, the Dirichlet eigenfield is non-physical in the sense that
it cannot be excited by sources in $\Omega_-$, whereas the Neumann eigenfield,
see \eqref{eq:Neumanneigeqs} and Figure~\ref{fig:Neumanneig}(g,h,i) 
for ordinary conductors,
indeed can be excited by such sources, and is of physical origin and not an 
artifact of any field representation. 
 We say that there exists an eigenfield  in the eddy current regime
if there are incident fields such that 
\begin{equation}  \label{eq:defneigfield}
   \frac{\|E^+\|/\max(k_-L, |\hat k^{-2}|)+\|E^-\|
+\|H^+\|+ \|H^-\|}{\|E^0\|+\|H^0\|} \to \infty,
\end{equation}
as $k_-\to 0$ along a curve in the $(k_-, k_+)$ plane, 
contained in the eddy current regime. 
Here the norm $\|\cdot\|$ of a vector field on $\Gamma$ is the
sum of the $H^{-1/2}(\Gamma)$ norm of its normal component
and the $H^{-1/2}(\curl,\Gamma)$ norm of its tangential part. 
See~\cite[Eq.~(65)]{HelsRose20}. 
The eigenfield is defined to be the limit of $\{E^+, E^-, H^+, H^-\}$,
normalized suitably.
In the sum in \eqref{eq:defneigfield}, we have scaled
   each of the four fields $\{E^+, E^-, H^+, H^-\}$ by the generic
   size of the corresponding measurable field $\{E^+, E^0+E^-, H^+,
   H^0+H^-\}$.

The classical exterior  Dirichlet eigenfield is a divergence- and curl-free electric field 
in $\Omega_-$ which is normal on $\Gamma$ and decays at $\infty$, 
resulting from a net charge in $\Omega_+$. 
However, $\Omega_+$ is assumed to have zero net charge
and thus such eigenfields cannot be excited by sources located in 
$\Omega_-$. Since we only assume sources in $\Omega_-$, we will
not see any Dirichlet eigenfields appearing. 
Our augmentations of null spaces related to the Dirichlet eigenfields build 
on~\eqref{eq:DirEminusaug}.
This condition excludes Dirichlet eigenfields $E^-$, since the
   maximum principle applied to the electric potential shows that
   $\nu\cdot E^-$ cannot change sign, which forces $E^-$ to be 
zero in $\Omega^-$.

The classical exterior Neumann eigenfield is a divergence- and curl-free
magnetic field in $\Omega_-$ which is tangential on $\Gamma$ and decays at
$\infty$, resulting from an electric current on the 
surface of genus $\ge 1$ of a superconductor. 
For a torus, the current is in the $\theta$ direction, and the Neumann 
eigenfield
is in the $\tau$ direction  on $\Gamma$. 
For an ordinary conductor, we see the Neumann eigenfield
appearing in \eqref{eq:maxwtranspr} as $k_-\to 0$ in the 
eddy current regime, as follows.
Instead of $E^+$ we consider the auxiliary field
$\widetilde E^+= E^+/\max(k_-L, |\hat k^{-2}|)$, since in the generic
scattering situation the field $\tilde E^+$ is of order $1$.
Setting $E^0=H^0=0$, the first equation in  \eqref{eq:maxwtranspr} 
shows that $\nu\times E^-=0$, since $E^+\to 0$ as $k_-\to 0$ in the
eddy current regime. 
This forces $E^-=0$ since $E^-$ is a divergence- and curl-free
  vector field in $\Omega_-$ and satisfies \eqref{eq:DirEminusaug}.
From the jump relation $\nu\cdot E^+= \hat k^{-2} \nu\cdot E^-$, 
which follows from the second equation in  \eqref{eq:maxwtranspr},
we conclude that $\nu \cdot \widetilde E^+=0$.
Therefore $\widetilde E^+$, as well as the eddy current
$J=\sigma E^+$, is a divergence- and curl-free vector field  in $\Omega_+$ with vanishing normal component at $\Gamma$.
Finally the magnetic field $H$ is a divergence-free vector field in $\R^3$ solving
\begin{equation}   \label{eq:Neumanneigeqs}
  \nabla\times H= 
  \begin{cases}
    \eta_0 J, & \text{in } \Omega_+,\\
    0, & \text{in } \Omega_-,
  \end{cases}
\end{equation} 
which follows from \eqref{eq:maxwtranspr} since $-i k_+\hat k=
  \eta_0 \sigma$.
To summarize, the De Rham cohomology space
$H^1(\Omega_+)$, see for example \cite[Sec. 10.6]{RosenGMA19},
parametrizes the Neumann eigenfields. More precisely the eddy current
$J=\sigma E^+$ is a divergence- and curl-free field which is tangential 
on $\Gamma$, and 
the associated magnetic field is obtained by applying the
Biot--Savart operator to $\eta_0 J$.
Note from \eqref{eq:Neumanneigeqs} 
for the Neumann eigenfield that $H$ and $J$ are 
independent of $\sigma$ for all ordinary conductors with $0<\sigma<\infty$.
In contrast to the Dirichlet eigenfields, the Neumann eigenfields can be excited
by sources in $\Omega_-$. Thus some
ill-conditioning is inevitable in any BIE
since the Neumann eigenfield is a physical eigenfield. 
By an inhomogeneous (L) augmentation
similiar to Example~\ref{ex:modelinhomoL}, we shall obtain a 
well-conditioned system and field representation, locating the 
ill-conditioning to the simpler computation of the right hand side
for the system.  
In Section~\ref{sec:numerics}, we demonstrate that the Neumann 
eigenfield  in an ordinary conductor  can be excited by the incident field
\begin{equation}   \label{eq:zcoil}
  E^0(x)= i c_2 H_1^{(1)}(k_-\rho) \mv \theta,\qquad H^0(x)= c_2 H_0^{(1)}(k_-\rho) \mv z,
\end{equation}
normalized with $c_2= 1/|H_1^{(1)}(k_-)|$ and
where $H_n^{(1)}$ is the first kind Hankel function of order $n$. 
This field can in principle be generated
by a thin wire along the $z$-axis made of a material  with
high relative permeability $\mu_r\gg 1$. By closing the wire in a large loop and exciting a
magnetic field inside the wire by a coil around the wire,  a field like \eqref{eq:zcoil} will be incident on $\Gamma$.

An important point  is that the Neumann eigenfield obtained in the PEC approximation is not the eigenfield appearing in ordinary conductors. The difference is illustrated in
 Figure~\ref{fig:Neumanneig} where
(g,h,i)  shows the Neumann eigenfield appearing in ordinary conductors. Note that the $H$ field penetrates into $\Omega_+$
and that the eddy current flows in the interior of $\Omega_+$, 
unlike the Neumann eigenfield for
superconductors shown in Figure~\ref{fig:Neumanneig}(a,b,c),
where the current flows on the surface $\Gamma$, as shown qualitatively in 
Figure~\ref{fig:Neumanneig}(a).  
Figure~\ref{fig:Neumanneig}(d,e,f)  shows an intermediate eigenfield in the case when the scaled conductivity grows inversely proportional to $k_-$ 
as $k_-\to 0$ so that the skin depth
is fixed in the sense that $k_+=10(1+i)$. Here we see how $J$ and $H$ 
begin to be expelled from $\Omega_+$.

Note that the ill-posedness of the MTP in the eddy current regime for $\Gamma$
of genus $1$, due to the existence of Neumann eigenfields, does not contradict conservation of energy.
Since we consider total permittivities $\epsilon_+$  that  are imaginary,
there will be no transmitted electric energy in $\Omega^+$, 
only a large magnetic energy.  
However, to excite the Neumann eigenfield requires an incident
   field like \eqref{eq:zcoil}, which requires a large  power  to be
   produced by the coil.

\section{Dirac (A$\infty$)}  \label{sec:finitecond}

For the remainder of this paper, we choose unit of length so that $L$ is of
order~$1$.
In \cite{HelsKarlRos20}, it was demonstrated that Dirac (A) works well near the quasi-static
limit $k_\pm\to 0$, provided $|k_+|$ and $k_-$ are comparable in size.
In this section, we first formulate a Dirac BIE, referred to as (A$\infty$),  that,
   after augmentation,  is intended to be used in the eddy current regime
   \eqref{eq:goodcond}.
The resulting BIE, referred to as (A$\infty$-aug), is insensitive to the genus of $\Gamma$ and
builds on \eqref{eq:betaDir}, but differs
in the choices of $\alpha', \beta'$ and the preconditioning $P'$.
It is demonstrated in  Section~\ref{sec:numerics} that Dirac (A$\infty$-aug)
performs well when the Neumann eigenfields are excited.

Dirac (A$\infty$) is defined by the parameters
\begin{equation}  \label{eq:Ainfty}
\begin{split}
  & \begin{bmatrix} r & \beta & \gamma & \alpha' & \beta' & \gamma'
  \end{bmatrix} = 
  \begin{bmatrix} \frac 1{\hat k} &  \xi & a & 
  \frac 1{|\hat k|\hat k} & \rev a &  \rev a \end{bmatrix}, \\
  P &= \begin{bmatrix} \frac {\hat k^2}{(|\hat k|+(\hat k\xi)^{-1})\japsigma} & 
  \frac {\hat k}{(1+a)\japsigma} & \mv{\frac {\hat k}{2\japsigma}}  &
 \frac 1{1+\rev a} & \frac 1{1+\hat k^{-2}} &  \mv{\frac 1{1+\rev a}}
   \end{bmatrix},  \\
  P' &= \begin{bmatrix} \frac \japsigma{\hat k^2}  & \japsigma  & \mv{\japsigma} &
  1 & 1 & \mv 1
     \end{bmatrix},  \\
  N &=\begin{bmatrix} \frac {\hat k^2}{(|\hat k|\hat k\xi+1)\japsigma} & \frac 1{(1+a)\japsigma}  & \mv{\frac 1{2\japsigma}} &
 \frac {\rev a}{1+\rev a} & \frac 1{1+\hat k^{2}} & \mv{\frac 1{1+\rev a} }
   \end{bmatrix},  \\
  N' &= \begin{bmatrix} \frac \japsigma{a\hat k} & \frac \japsigma{|\hat k|}  & \mv{\frac \japsigma{\hat k}} &
  1 & 1 & \mv{\rev a}
     \end{bmatrix}.
\end{split}
\end{equation}
Here $a= \hat k/|\hat k|$, 
$\xi$ is as in \eqref{eq:xi}, and
\begin{equation}   \label{eq:sigmaparam}
\japsigma= 1+|k_+\hat k|.
\end{equation}
We have the following behaviour of Dirac (A$\infty$):
\begin{itemize}
\item
The coefficients in $G$ of Dirac (A$\infty$) are uniformly bounded and, similar to
$G$ of Dirac (A), the $G$ of Dirac (A$\infty$) is close to a nilpotent operator modulo compact operators
on smooth $\Gamma$ due to 
cancellations in blocks (1:2,3:4) and (7:8,5:6).
Furthermore, the norm of the (5:8,1:4) block is of order $k_-\japsigma$, whereas the norm of the 
(1:4,5:8) block is of order $|k_+\hat k/\japsigma|$.
\item
All  entries  in $N$ are bounded.
Inserting $N'$ and $P'$ into \eqref{eq:projdens}, it is seen that
\begin{equation}  \label{eq:Ascale}
  \|E^\pm|_\Gamma\|\lesssim \|h\| 
   \quad\text{and}\quad \|H^\pm|_\Gamma\|\lesssim \japsigma\|h\|.
\end{equation}
This possibility of having large transmitted and scattered fields, even
if the density $h$ is not large, explains why (A$\infty$) is able to 
accurately compute the fields when the Neumann eigenfields are excited. 
\item
For fixed $k_\pm\ne 0$, the choice of $\beta, \gamma, \alpha', \beta', \gamma'$
guarantees invertibility of $I+G$.
The limit operator $I+G_0$ as $k_\pm\to 0$ in the eddy current regime
is a Fredholm operator of index zero, and it has nullity $1$  regardless 
of the genus of $\Gamma$.
See~\ref{app:aug}.
At high conductivities, this analysis breaks down, but computations
suggest that the null space remains one-dimensional. 
\end{itemize}

Following Section~\ref{sec:chooseaug}, we  make one (L) augmentation
to remove the Dirichlet eigenfield. 
 Recall from Section~\ref{sec:physics} that this eigenfield cannot be
  excited by 
sources in $\Omega_-$, and therefore 
a homogeneous (L) augmentation is appropriate. 
For details we refer to~\ref{app:aug}. 
The field representation for Dirac (A$\infty$) fails to be a Fredholm map as 
$k_-\to 0$, as discussed in Section~\ref{sec:conditioning}, 
and (R) augmentations
are therefore not applicable. 

\subsection{Dirac (A$\infty$-aug)}   \label{sec:Ainf-aug}
To  remove  the Dirichlet eigenfield we make the homogeneous (L) augmentation 
of Dirac (A$\infty$)
\begin{equation} \label{eq:bc1}
\begin{split}
  c^1_D h &=\barint_\Gamma (E_{k_-}^- P'h)_6 d\Gamma,  \\
  b^1_D &= \begin{bmatrix} 0 & 0 & \mv 0 & 0 & 1 & \mv 0  \end{bmatrix}^T. 
\end{split}
\end{equation}
Here $\barint_\Gamma f d\Gamma$ denotes the average value of a function
$f$ on $\Gamma$.
We thus obtain a Dirac (A$\infty$-aug),
intended for $\Gamma$ of any genus. Given $f^0$, we solve 
the augmented system
\begin{equation}   \label{eq:Ainfty-augsystem}
  h+Gh+b^1_D(c^1_D h)= 2Nf^0,
\end{equation} 
with $G$, $c^1_D,b^1_D$ from \eqref{eq:G}, \eqref{eq:bc1}, and the parameters \eqref{eq:Ainfty} above.
This yields $h$, from which we compute
the fields using \eqref{eq:projdens}.

\section{Dirac (B)}   \label{sec:DirB}

 Recall  that we are using unit of length so that $L$ is of order~$1$.
In this section, we first formulate a Dirac BIE, referred to as (B), for  the
   MTP($k_-,k_+,\hat k^2$) that, after augmentation,  is intended to be
   used in the eddy current regime \eqref{eq:goodcond}.
Dirac (B) is radically different in the choices of parameters from both
Dirac (A) and Dirac (A$\infty$), 
and is designed to compute each of the fields $E^\pm, H^\pm$ 
accurately when the 
Neumann fields are not excited.
It is defined by the parameters
\begin{equation}   \label{eq:Bparam}
\begin{split}
  &\begin{bmatrix} r & \beta & \gamma & \alpha' & \beta' & \gamma'
  \end{bmatrix} = 
  \begin{bmatrix} \frac 1{\hat k} &  \frac{\hat k}{|\hat k|^2} & \frac{\hat k^2}\xi & 
  \frac 1\xi & \frac{1}{\hat k} &  \frac 1\xi \end{bmatrix}, \\
  P &= \begin{bmatrix} \frac 1{\xi\hat k^{-1}+ a^{-2}} & \frac {\hat k}{\xi+1} & \mv{\frac {\hat k}2} & 
   \frac{\hat k^2}\japsigma \frac 1{1+\xi\hat k^{-2}} &  
   \frac{\hat k^2}\japsigma \frac 1{\xi+\hat k^{-1}} & \mv{\frac 1{1+\xi\hat k^{-2}}} \end{bmatrix}, \\
  P' &= \begin{bmatrix} 1 & 1 & \mv 1 & \frac \japsigma{\hat k^2} & 
   \frac \japsigma{\hat k} & \mv 1
     \end{bmatrix},  \\
  N &= \begin{bmatrix} \frac 1{1+\xi a^2/\hat k} & \frac 1{\xi+1} & \mv{\frac 12} &
   \frac \xi\japsigma \frac 1{1+\xi\hat k^{-2}}  &  \frac 1\japsigma\frac 1{\xi+\hat k^{-1}} & \mv{\frac 1{1+\xi\hat k^{-2}}}   \end{bmatrix},  \\
  N' &= \begin{bmatrix} \frac \xi{\hat k}  &  \frac \xi{\hat k}  &  \mv{\frac 1{\hat k}}  &  \frac \japsigma{\hat k^2}  & \frac {\xi\japsigma}{\hat k^2}  &  \mv{\frac \xi{\hat k^2}}      
    \end{bmatrix}, 
\end{split}
\end{equation} 
  where $a=\hat k/|\hat k|$, $\xi$ is as in \eqref{eq:xi}, and
  $\japsigma$ is as in \eqref{eq:sigmaparam}. We have the following
  behaviour of Dirac~(B): 
\begin{itemize}
\item
The coefficients in $G$ are uniformly bounded for all $|\hat k|\gtrsim 1$.
The operator $G$ is close to a nilpotent operator modulo compact operators
on smooth $\Gamma$, but unlike Dirac (A) and (A$\infty$), this is due to 
cancellations mainly in blocks (3:4,1:2) and (5:6,7:8). 
Furthermore, the norm of the (5:8,1:4) block is of order $|k_+\hat k/\japsigma|$, whereas the norm of the 
(1:4,5:8) block is of order $k_-\japsigma$.

\item
All  entries  in $P', N, N'$ are uniformly bounded for all $|\hat k|\gtrsim 1$, and $N'$ is adapted to the generic sizes of the fields, as discussed in Section~\ref{sec:physics}, so that
\begin{equation}   \label{eq:Bscale}
  (|\hat k^2|/\japsigma)\|E^+|_\Gamma\|,  \quad \|E^-|_\Gamma\|, \quad
  \|H^+|_\Gamma\|, \quad \|H^-|_\Gamma\|
\end{equation}
are all $\lesssim\|h\|$.
\item
For fixed $k_\pm\ne 0$, the choice of $\beta, \gamma, \alpha', \beta', \gamma'$
guarantees invertibility of $I+G$.
The limit operator $I+G_0$ as $k_-\to 0$ in the eddy current regime
  is a Fredholm operator of index zero, but its nullity depends on
  $k_+$ as well as on the genus of $\Gamma$.
See~\ref{app:aug}.
\end{itemize}

Following Section~\ref{sec:chooseaug}, we make 
a number of augmentations.
First the field representation requires one or two (R) augmentations,
  depending on the genus of $\Gamma$, to become well-conditioned.
Once this is done, one or two (L) augmentations 
are needed to remove the Dirichlet eigenfield and the null space
   associated with the Neumann eigenfield.
It is only the Neumann eigenfield which can be excited by Maxwell sources
in $\Omega_-$, so that an inhomogeneous (L) augmentation is needed.
A final (L) augmentation is also needed to remove an
eigenfield which is present in the Dirac equation, but which is outside the
Maxwell equations.  
For details we refer to~\ref{app:aug}.

\subsection{Dirac (B-aug0)}  \label{sec:B-aug0}

   Consider $\Gamma$ of genus $0$.
   We first make an (R) augmentation 
\begin{equation}   \label{eq:bc2}
\begin{split} 
  c^R_D h &=\barint_\Gamma h_6 d\Gamma, \\
  b^R_D &= 2N E_{k_-}^- e_6,
\end{split}
\end{equation}
where $e_6= \begin{bmatrix} 0 & 0 & \mv 0 & 0 & 1 & \mv 0   \end{bmatrix}^T$,
which has the effect of adding the Dirichlet eigenfield to the
    field  representation. That field  is missing in Dirac (B).
This is needed to avoid a null space for the system. But since
the Dirichlet eigenfield cannot be excited by sources in $\Omega_-$, 
we also make a homogeneous (L) augmentation
\begin{equation}   \label{eq:bc3}
\begin{split} 
  c^2_D h &=\barint_\Gamma (E_{k_-}^- (P' h+e_6 (c^R_Dh)) )_6 d\Gamma, \\
  b^2_D &= e_1,
\end{split}
\end{equation}
where $e_1= \begin{bmatrix} 1 & 0 & \mv 0 & 0 & 0 & \mv 0   \end{bmatrix}^T$,
which again removes the Dirichlet eigenfield.
The two augmentations $b^R_Dc^R_D$ and $b^2_Dc^2_D$ together 
make the system stably solvable when $k_+\hat k\gtrsim 1$.
When $k_+\hat k \ll 1$ we also have a second eigenfield for the DTP.
   This field comes from one of the auxiliary HTPs, and we 
   make a homogeneous (L) augmentation 
\begin{equation} \label{eq:bc4}
  \begin{split}
  c^1_H h &=\barint_\Gamma (E_{k_+}^+\hat k N'h)_1 d\Gamma,  \\
  b^1_H &= \begin{bmatrix} 0 & 0 & \mv 0 & 0 & 1 & \mv 0  \end{bmatrix}^T.
\end{split}
\end{equation} 
We thus obtain a Dirac  BIE, referred to as (B-aug0),  intended for
   $\Gamma$ of genus $0$.
Given $f^0$, we solve 
the augmented system
\begin{equation}   \label{eq:B0system}
  h+Gh+(b^R_Dc^R_D+b^2_Dc^2_D+ \chi  b^1_Hc^1_H) h= 2Nf^0,
\end{equation} 
with $G$, $b^R_Dc^R_D$, $b^2_Dc^2_D$ and $b^1_Hc^1_H$ from \eqref{eq:G}, 
\eqref{eq:bc2}, \eqref{eq:bc3} and \eqref{eq:bc4} using the parameters \eqref{eq:Bparam}.
Here
\begin{equation}
\chi=
\begin{cases}
  0, & |k_+\hat k| \gtrsim 1, \\
  1, & |k_+\hat k| \ll 1,
\end{cases}
\end{equation}
but numerically $\chi=1$ seems to work in both cases. 
This yields $h$, from which we compute
the fields using the augmented field representation 
\begin{equation}     \label{eq:projdensB0}
\begin{split}   
  F^+|_\Gamma &= E_{k_+}^+ N'h, \\
   F^-|_\Gamma &= -E_{k_-}^- (P' h+e_6 (c^R_Dh)). 
\end{split}
\end{equation}

\subsection{Dirac (B-aug1)}  \label{sec:B-aug1}

 Consider an axially symmetric $\Gamma$ of genus $1$. 
  (The construction in this section generalizes to arbitrary $\Gamma$
   of genus $1$, if the $\tau$ and $\theta$ directions are suitably defined.) 
  As for genus $0$ we make the augmentations $b^R_Dc^R_D$ and $b^2_Dc^2_D$
  to remove the Dirichlet eigenfield.
  But for genus $1$, we also need to make an (R) augmentation 
  \begin{equation}  \label{eq:bc5}
\begin{split}
  c^R_N h &=\barint_\Gamma h_8 d\Gamma, \\
  b^R_N &=  2  \frac\japsigma{\hat k^2} P E_{k_+}^+ e_8,
\end{split}
\end{equation}
which has the effect of adding the Neumann eigenfield to the field
    representation. That field  is also missing in Dirac (B).
This gives the field representation
\begin{equation}   \label{eq:projdensB1}
\begin{split}   
  F^+|_\Gamma &= E_{k_+}^+(N'h + \frac\japsigma{\hat k^2}e_8 (c^R_N h)), \\
   F^-|_\Gamma &= -E_{k_-}^- (P' h+e_6 (c^R_Dh)). 
\end{split}
\end{equation} 
We use $e_8= \begin{bmatrix} 0 & 0 & 0 & 0 & 0 & 0 & 0 & 1  \end{bmatrix}^T$, which seems to work numerically, but the derivation in~\ref{app:aug} uses $e_8$ with the $\theta$ component on
$\Gamma$ of the interior Neumann eigenfield in its last component.  

For $k_+\hat k\ll 1$, we need to adjust the HTP augmentation
$b^1_Hc^1_H$ to \eqref{eq:projdensB1},
and we set
\begin{equation} \label{eq:bc6}
  \begin{split}
  c^2_H h &=\barint_\Gamma (E_{k_+}^+(\hat k N'h+\frac\japsigma{\hat k}e_8 (c^R_N h)))_1 d\Gamma,  \\
  b^2_H &= \begin{bmatrix} 0 & 0 & \mv 0 & 0 & 1 & \mv 0  \end{bmatrix}^T.
\end{split}
\end{equation} 

What remains is the more subtle inhomogeneous (L) augmentation of the Neumann eigenfield.
As discussed in Section~\ref{sec:physics}, the MTP is ill-posed in itself due to the
presence of the Neumann eigenfield, but we make an inhomogeneous (L) augmentation
 that  makes the system well-conditioned at the expense of adding an unbounded 
functional $d^1_N$ to the  right-hand  side. Let
\begin{equation}   \label{eq:bc7}
    \begin{split}
  c^1_N h &=\barint_\Gamma \Big( E_{k_+}^+((\hat k^2/\japsigma) N'+ e_8 c^R_N )h\Big)_8 w d\Gamma \\
  &+ \frac{\hat k^2}{\japsigma}\barint_\Gamma \Big( E^-_{k_-}P' h_{1:5}\Big)_8 w d\Gamma  
+ \frac{\hat k^2}{\japsigma}\barint_\Gamma \tfrac 12\Big((E_0-E_{k_-}) h_{7:8}\Big)_8 w d\Gamma, \\
  b^1_N &= \begin{bmatrix} 0 & 0 & 0 & 0 & 0 & 0 & 0 & 1  \end{bmatrix}^T= e_8, \\
  d^1_N f^0 &= \frac {\hat k^2}{\japsigma} \barint_\Gamma ( f^0 )_8 w d\Gamma,
\end{split} 
\end{equation}
where  $w=\tau\cdot H|_\Gamma$ is a weight function discussed below,
and $H$ denotes an exterior PEC Neumann eigenfield. 
Although the factor $\hat k^2/\japsigma$ makes $d^1_N$ unbounded, 
it turns out that $d^1_N f^0$ stays bounded  by $f^0$  unless the Neumann eigenfield
is excited. 
More precisely, $d^1_N f^0/\max_\Gamma |f^0|$, measures how 
much of the Neumann eigenfield will be excited. 
Indeed, an application of Stokes' theorem shows that
\begin{equation}
   d^1_N f^0 = \frac{ik_+\hat k}{\japsigma |\Gamma|}\int_{\Omega_-}
H^0\cdot H dx,
\end{equation}
where $|\Gamma|$ denotes the area of $\Gamma$ and $H^0$ is, as before, the incident 
magnetic field. 

The equation $c^1_N h= d^1_N f^0$ is seen to be equivalent to continuity
of the $\theta$ component of the electric field across $\Gamma$, and normalized 
 with  the generic size
of $E^+$.
See~\ref{app:aug}.
We denote by $h_{1:5}$ the density $h$ with components 6:8 set to zero,
and likewise $h_{7:8}$ denotes the density $h$ with components 
1:6 set to zero.
For this augmentation to work, it is essential to use a specific positive weight function $w$:
the $\tau$-component of the exterior PEC Neumann eigenfield, normalized so
that $\barint_\Gamma w d\Gamma=1$.
The weight function can be computed as $w=\theta\cdot f$,
where $f$ is the solution to the size $2\times 2$ block
   eigenfunction equation
\begin{equation}  \label{eq:1stweq}
  (\mv I+\mv M_0)f=0,
\end{equation}
and $\mv M_0$ is the static magnetic dipole operator \eqref{eq:magndipoleop}. 
Equivalently, and perhaps better from a numerical point of view, $w$ can be
computed as
\begin{equation}
  w= \tau\cdot H_0 + K^\tau_0 \psi,
\end{equation}
where $\psi$ is the solution to the single block Fredholm second
   kind integral equation
\begin{equation}  \label{eq:2ndweq}
  (I+K^\nu_0)\psi= -\nu\cdot H_0,
\end{equation}
which can be solved iteratively. 
Note also that the system in \eqref{eq:2ndweq} is half the size of the 
system~\eqref{eq:1stweq}. 
Here $K^\tau_0$ is defined as $K^\nu_0$, but with $\nu(x)$ replaced by $\tau(x)$,
and $H_0$ denotes the magnetic field produced by a steady current in a circular wire
around $\Omega_+$ (or any computable divergence- and curl-free vector field in $\Omega_-$ which 
is not a gradient field).

In total, we obtain a Dirac BIE referred to as (B-aug1)
and intended for $\Gamma$ of genus $1$. Given $f^0$, we solve 
the augmented system
\begin{equation}      \label{eq:B1system}
  h+Gh+(b^R_Dc^R_D+b^2_Dc^2_D+b^R_Nc^R_N+ b^1_Nc^1_N+\chi b^2_Hc^2_H)h= 2Nf^0+ b^1_N(d^1_Nf^0),
\end{equation} 
with $G$, $b^R_Dc^R_D$, $b^2_Dc^2_D$, $b^R_Nc^R_N$, $b^1_Nc^1_Nd^1_N$, $b^2_Hc^2_H$, from \eqref{eq:G}, \eqref{eq:bc2}, \eqref{eq:bc3}, \eqref{eq:bc5}, \eqref{eq:bc7}, \eqref{eq:bc6} using the parameters \eqref{eq:Bparam}.
This yields $h$, from which we compute
the fields using \eqref{eq:projdensB1}.

\section{Numerical examples}   \label{sec:numerics}

The properties of Dirac (A$\infty$-aug) from Section~\ref{sec:Ainf-aug} and (B-aug0/1) from Sections~\ref{sec:B-aug0} and~\ref{sec:B-aug1} are now illustrated in a series of numerical examples involving two objects $\Omega_+$ with axially symmetric surfaces:
\begin{itemize}
\item The ``rotated starfish'' has a surface $\Gamma$ of genus $0$, with
   generating curve
\begin{equation}
r(s)=(1+0.25\sin(5s))(\cos(s),\sin(s)), \qquad s\in [-\pi/2,\pi/2], \label{eq:rotstarfish} \end{equation} and generalized diameter $L\approx 2.4$.
\item The ``starfish torus'' has a surface $\Gamma$ of genus $1$, with
   generating curve
\begin{equation}
r(s)=1+0.5(1+0.25\sin(5s))(\cos(s),\sin(s)), \qquad s\in [-\pi,\pi], \label{eq:starfishtorus} \end{equation} and generalized diameter $L\approx 3.2$.
\end{itemize}
 The incident fields, when such are present, are either~\eqref{eq:partialwave} or~\eqref{eq:zcoil}. These are all of order $1$ on $\Gamma$, except $H^0$
in~\eqref{eq:zcoil}, which is of order $|k_-\log(k_-)|$. 

Our computations rely on Fourier--Nystr{\"o}m discretization~\cite{YouHaoMar12}, where a sequence of decoupled modal problems, with modal index $n$, are solved using a mix of $16$th- and $32$nd-order composite panel-based discretization and where linear systems are solved iteratively using GMRES with a stopping criterion threshold of machine epsilon in the estimated relative residual. The implementation of this numerical scheme is the same as that used in~\cite[Sec.~10]{HelsKarlRos20}. In particular, the scheme is thoroughly verified for Dirac (A) and genus $0$ in~\cite[Sec.~10.3]{HelsKarlRos20} both under mesh refinement and by comparison with semi-analytic results. In the present work, Dirac
(A$\infty$-aug) and (B-aug0/1) are verified against Dirac (A) to the extent possible. The codes are implemented in {\sc Matlab}, release 2020a, and executed on a workstation equipped with an Intel Core i7-3930K CPU and 64 GB of RAM.

 Both fields~\eqref{eq:partialwave} and~\eqref{eq:zcoil}  are axially symmetric and excite only the mode $n=0$, which is the only Fourier mode affected by our augmentations for any incident field. More precisely, in all our augmentations $bc$ for Dirac
(A$\infty$-aug) and (B-aug0/1), the vector $b$ is a mode-$0$ function, whereas $ch=0$ for all mode-$n$ functions $h$ with $n\ne 0$. This is straightforward to verify, given the fact that $E_k^\pm h$ is a mode-$n$ function, whenever $h$ is such a function. To see the latter, assume that $h_\alpha= e^{in\alpha} h$, where $h_\alpha$ denotes the function $h$ rotated an angle $\alpha$ around the $z$-axis. Write $h= h^++h^-$ in the splitting in Hardy subspaces from \cite[Thm.~9.3.9]{RosenGMA19}, where $h^\pm= E_k^\pm h$. Denote by $h_\alpha^\pm$ the function $h^\pm$ rotated as above. We have \begin{equation}
     h_\alpha^++ h_\alpha^-=h_\alpha=e^{in\alpha} h= e^{in\alpha} h^+ + e^{in\alpha}h^-.
\end{equation}
By the uniqueness in this splitting and the rotational invariance of $\dirac F= ikF$ and the Dirac radiation condition, it follows that we must have $h^\pm_\alpha= e^{in\alpha} h^\pm$. Thus $E_k^\pm h= h^\pm$ are also mode-$n$ functions as claimed.

Several of our experiments result in field images and error images.
When assessing the accuracy of computed fields, and in the absence of semi-analytic results, we adopt a procedure where to each numerical solution we also compute an overresolved reference solution, using roughly 50\% more points in the discretization of the system under study. The absolute difference between these two solutions is denoted the {\it estimated absolute error}. The fields are always computed at $90,\!000$ field points on a Cartesian grid in the computational domains shown.

\subsection{The number of accurate digits in field evaluations}
\label{sec:relerror}

Since the transmitted and scattered fields differ much in size, it is important to measure
their relative errors appropriately.
In the exterior $\Omega_-$, the measurable fields are $E^0+E^-$ and
$H^0+H^-$. Hence it is motivated to normalize the error in $E^-$ and $H^-$
by the maximum of $|E^0+E^-|$ and $|H^0+H^-|$ in $\Omega_-$, 
respectively.
In the interior $\Omega_+$, the measurable fields are $E^+$ and
$H^+$. Hence it is motivated to normalize the error in $E^+$ and $H^+$
by the maximum of $|E^+|$ and $|H^+|$ in $\Omega_+$. 
However, note that when $k_\pm\approx 0$, all field components are
   almost harmonic functions and the maximum principle for such
   functions motivates replacing the above maxima by the maxima of each
   field component on $\Gamma$. Summarizing, we use the relative errors 
\begin{equation}   \label{eq:relativeerrordef}
\left\{ \frac{\max_{\Omega_+}|E^+_\text{err}|}{\max_{\Gamma}|E^+|},
 \frac{\max_{\Omega_-\cap \mathcal{D}}|E^-_\text{err}|}{\max_{\Gamma}|E^0+E^-|}, \frac{\max_{\Omega_+}|H^+_\text{err}|}{\max_{\Gamma}|H^+|},
 \frac{\max_{\Omega_-\cap \mathcal{D}}|H^-_\text{err}|}{\max_{\Gamma}|H^0+H^-|} \right\}
\end{equation} 
in the four fields $\{E^+, E^-, H^+, H^-\}$, where $\mathcal{D}$ denotes the
   computational domain and $F_\text{err}$ denotes the estimated
   absolute error in a field $F$. The {\em number of accurate digits}
   in a field $F$ is 
\begin{equation}   \label{eq:Ydigits}
   Y= -\text{round}(\log_{10}\epsilon),
\end{equation}
where $\epsilon$ denotes the relative error defined in 
\eqref{eq:relativeerrordef}.

Note that the relative error  in  the equally measurable
eddy current $J$, will be the same as the relative error in $E^+$.

\subsection{The high conductivity ``rotated starfish''}   \label{sec:genus0}

\begin{figure}[t!]
\centering
\includegraphics[height=50mm]{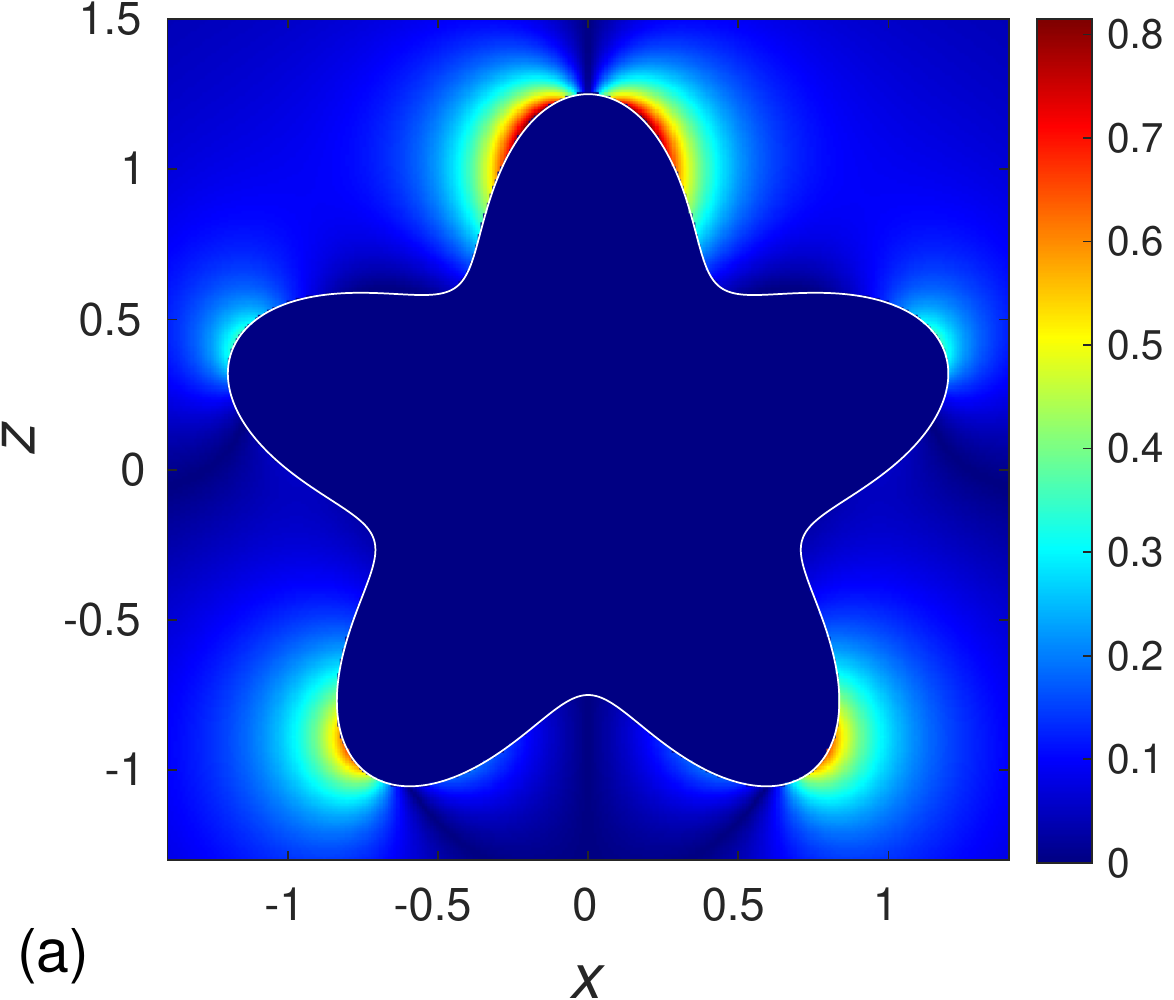}
\includegraphics[height=50mm]{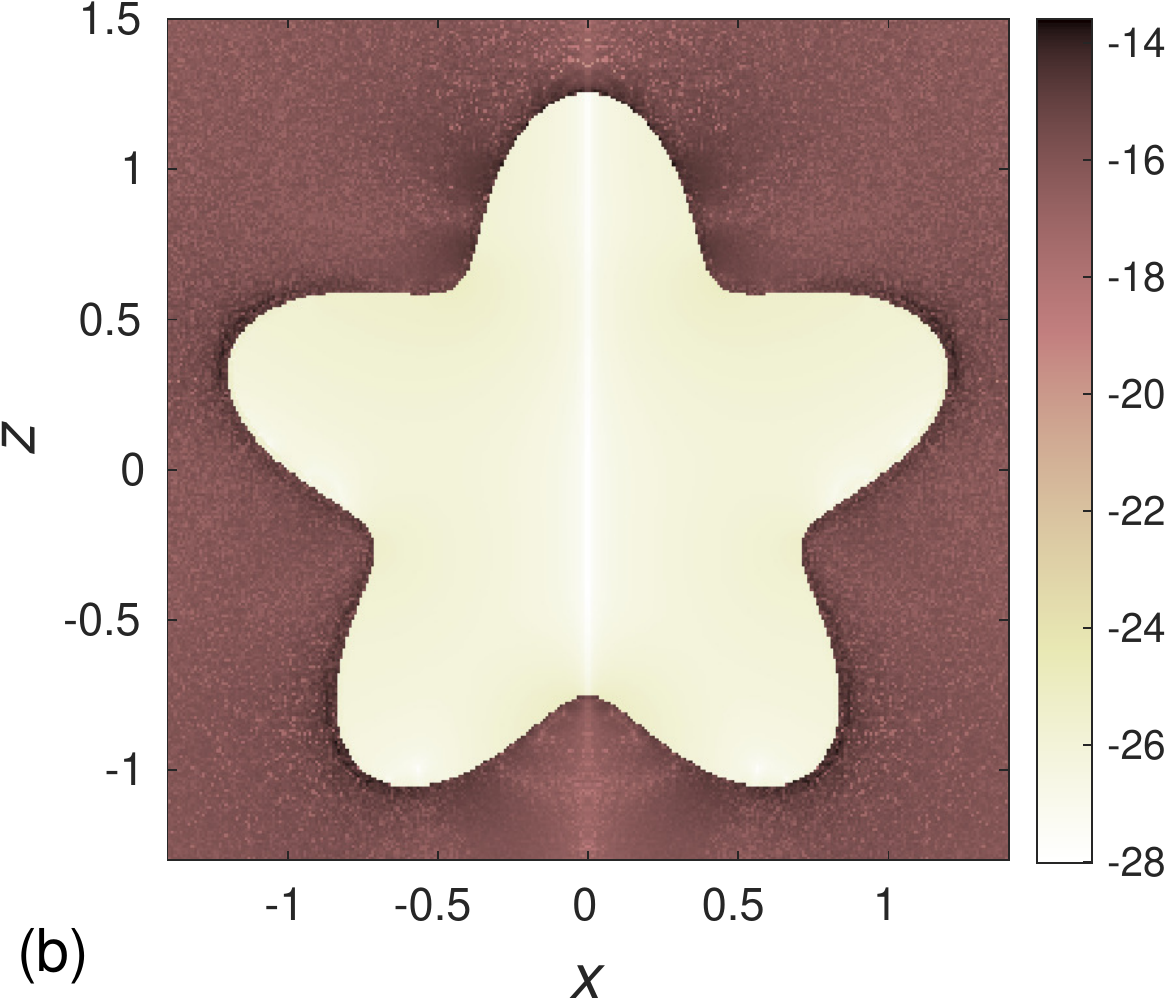}

\vspace{1mm}
\includegraphics[height=50mm]{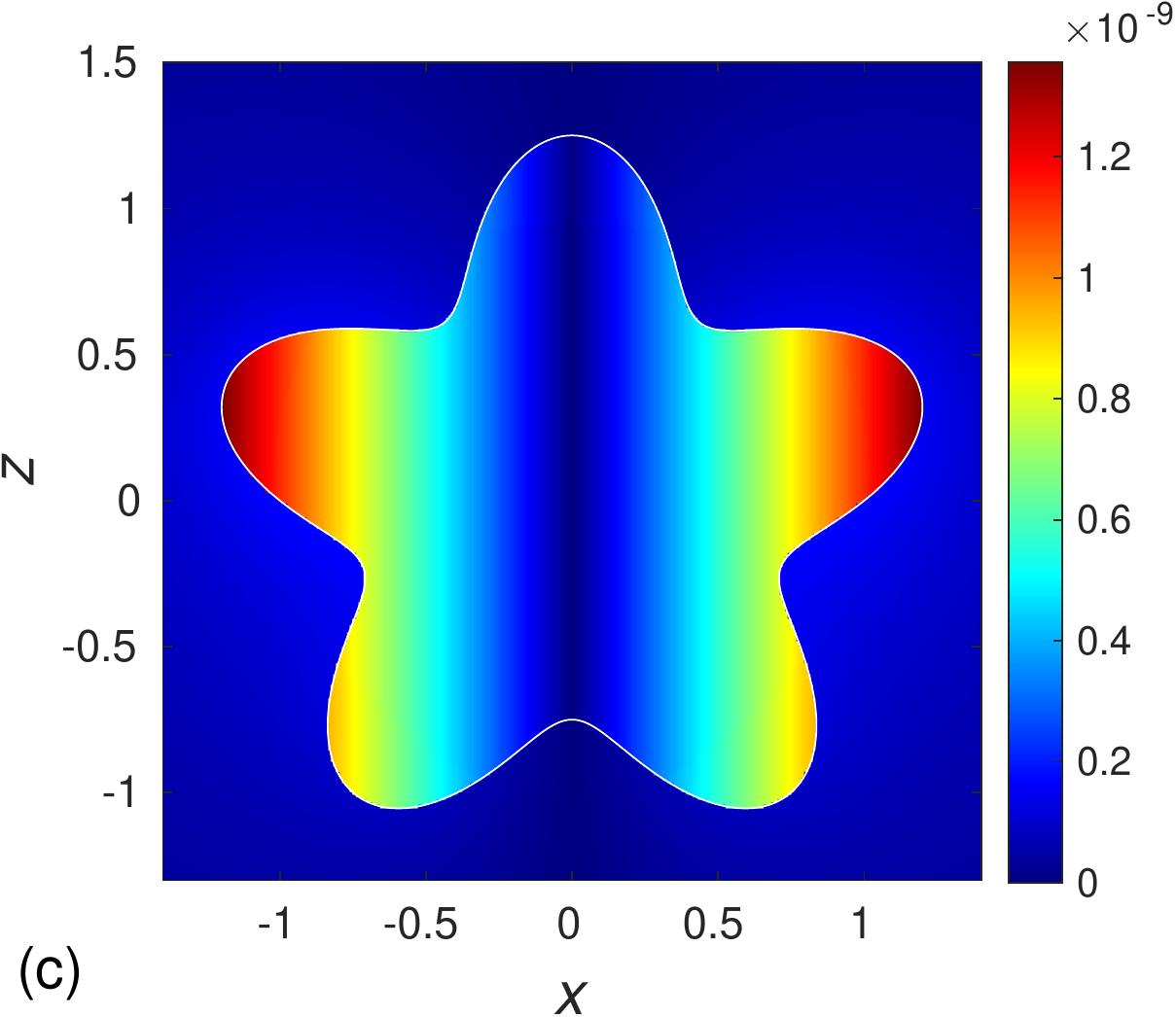}
\includegraphics[height=50mm]{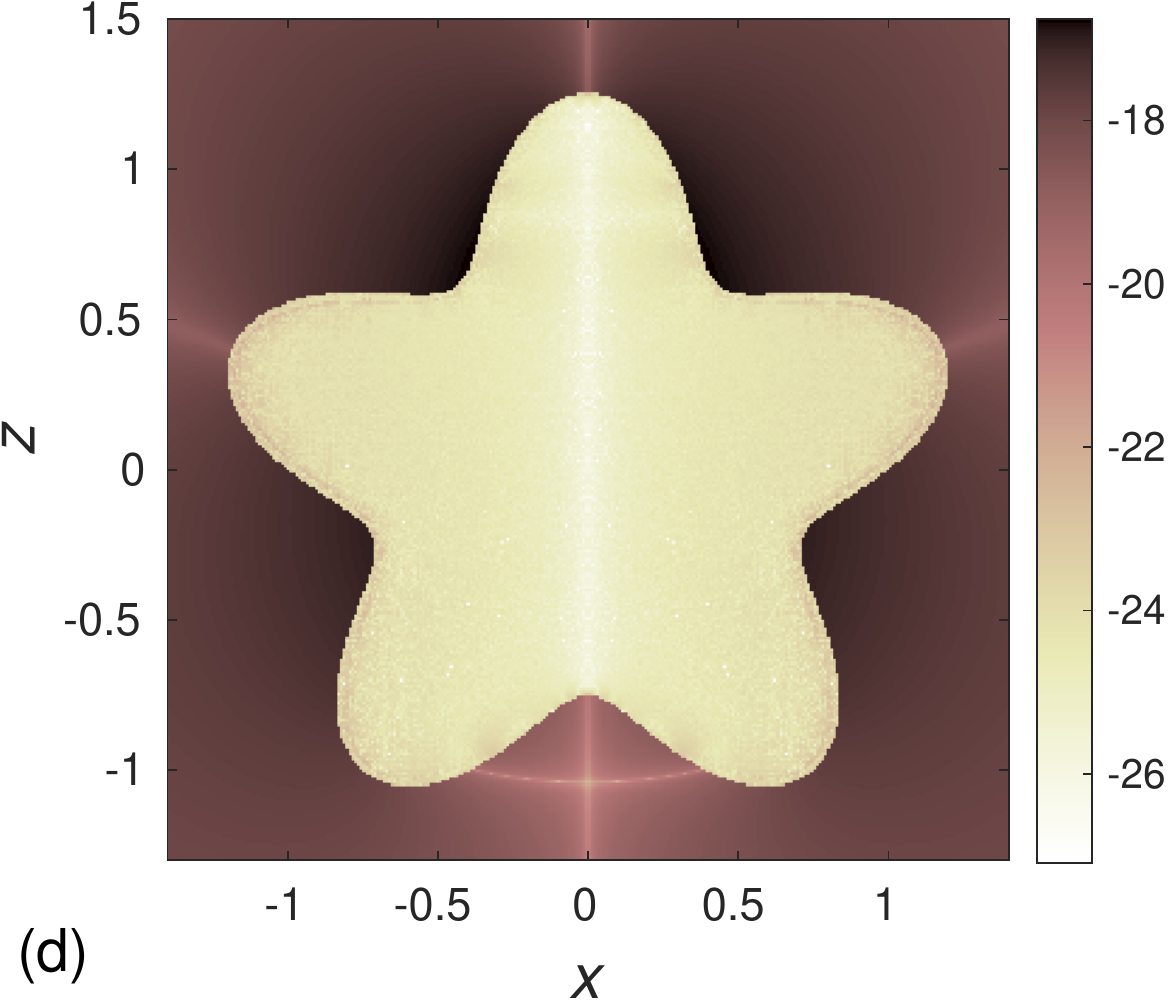}

\vspace{1mm}
\includegraphics[height=50mm]{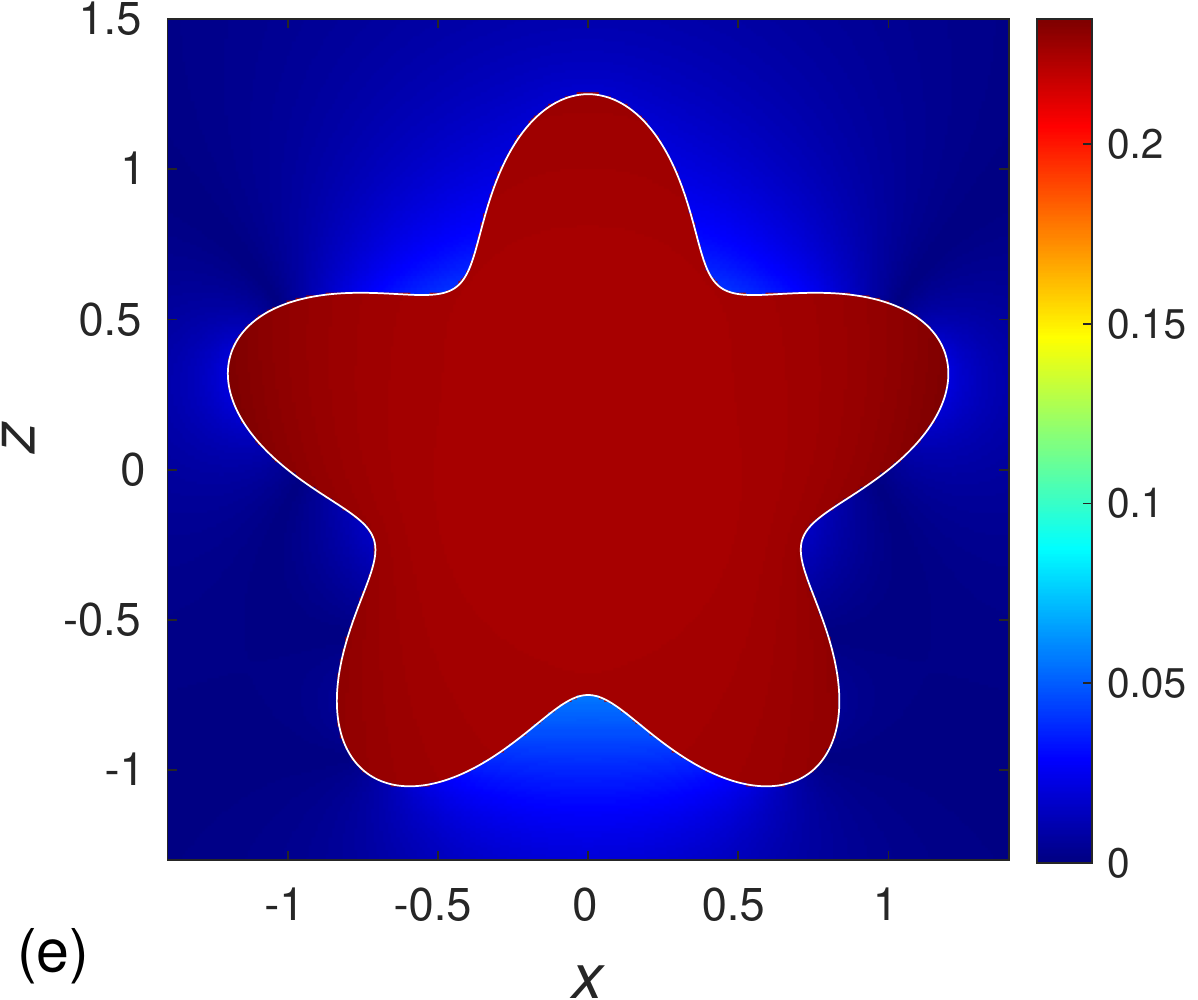}
\includegraphics[height=50mm]{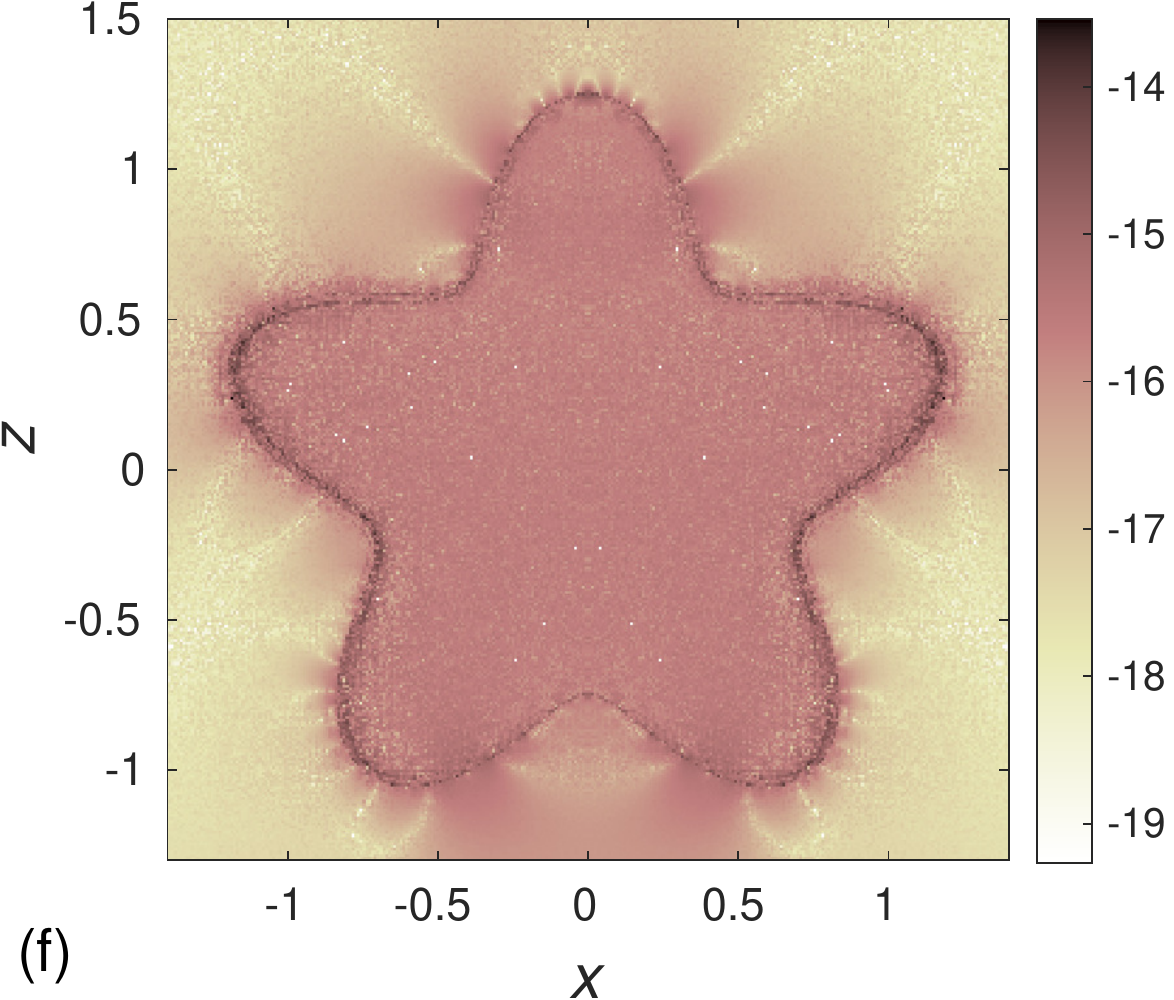}
\caption{\sf Field images for scattering of \eqref{eq:partialwave} by the
``rotated starfish''~\eqref{eq:rotstarfish} at $k_-=10^{-8}, k_+=1+i$;
(a) scattered/transmitted amplitude $|E_\rho|$;
(c) scattered/transmitted amplitude $|E_\theta|$;
(e) scattered/transmitted amplitude $|H_z|$.
(b,d,f) $\log_{10}$ of estimated absolute error of complex fields 
using (B-aug0), for (a,c,e) respectively.
}
\label{fig:genus0}
\end{figure}

We consider the field \eqref{eq:partialwave} incident on the 
``rotated starfish'' with  $L\approx 2.4$ cm  defined by \eqref{eq:rotstarfish}, 
at wavenumbers $k_-= 10^{-8}$ ${\rm cm}^{-1}$
and $k_+= 1+i$ ${\rm cm}^{-1}$. This corresponds to a frequency  $\omega\approx 300$ rad/s 
and conductivity $\sigma\approx 5.3\cdot 10^7$ S/m, which for example occurs for copper.
Since the surface $\Gamma$ has genus $0$, we compute the scattered
and transmitted fields using Dirac (B-aug0).
The number of  accurate digits~\eqref{eq:Ydigits}  obtained for the fields  $\{E^+, E^-, H^+, H^-\}$  are $\{13, 14, 13, 14\}$ and GMRES needs
   $33$ iterations.
The amplitudes of a selection of components of the fields $E^\pm$ and $H^\pm$
are shown in Figure~\ref{fig:genus0} along with field errors. 
Note that, unlike~\eqref{eq:Ydigits}, the errors shown in Figures~\ref{fig:genus0}--\ref{fig:finitegenus1-z} are estimated absolute
errors.  
The component $E_z$ is similar to $E_\rho$ and shows a scattered electric 
field $E^-$ which is close to normal on $\Gamma$. 
Figure~\ref{fig:genus0}(c) shows that the transmitted electric field $E^+$ is of order $10^{-9}$.
Figure~\ref{fig:genus0}(e) shows that most of the incident magnetic field $H^0$,
which is mainly in the  $z$ direction,  is transmitted into the non-magnetic object.
The components $H_\rho$ and $H_\theta$ (not shown)  are of order $10^{-2}$ and $10^{-9}$ respectively.   
   This agrees with the discussion beginning Section~\ref{sec:physics}
  which  predicts all the fields to be of order $1$ except $E^+$, which
  is of order $k_-= 10^{-8}$ ${\rm cm}^{-1}$.

For comparison, solving the same scattering problem with Dirac (A$\infty$-aug)
takes $36$ iterations and gives $\{5, 14, 6, 6\}$ accurate digits
   for  $\{E^+, E^-, H^+, H^-\}$. 
This gives numerical support for choosing (B-aug0) for eddy current
   scattering with surfaces of  genus $0$. This is so  since
   (A$\infty$-aug) computes the fields at the wrong scale
   \eqref{eq:Ascale}, with subsequent loss of accuracy.

\subsection{The high conductivity ``starfish torus''}  \label{sec:hightorus}

\begin{figure}[t!]
\centering
\includegraphics[height=50mm]{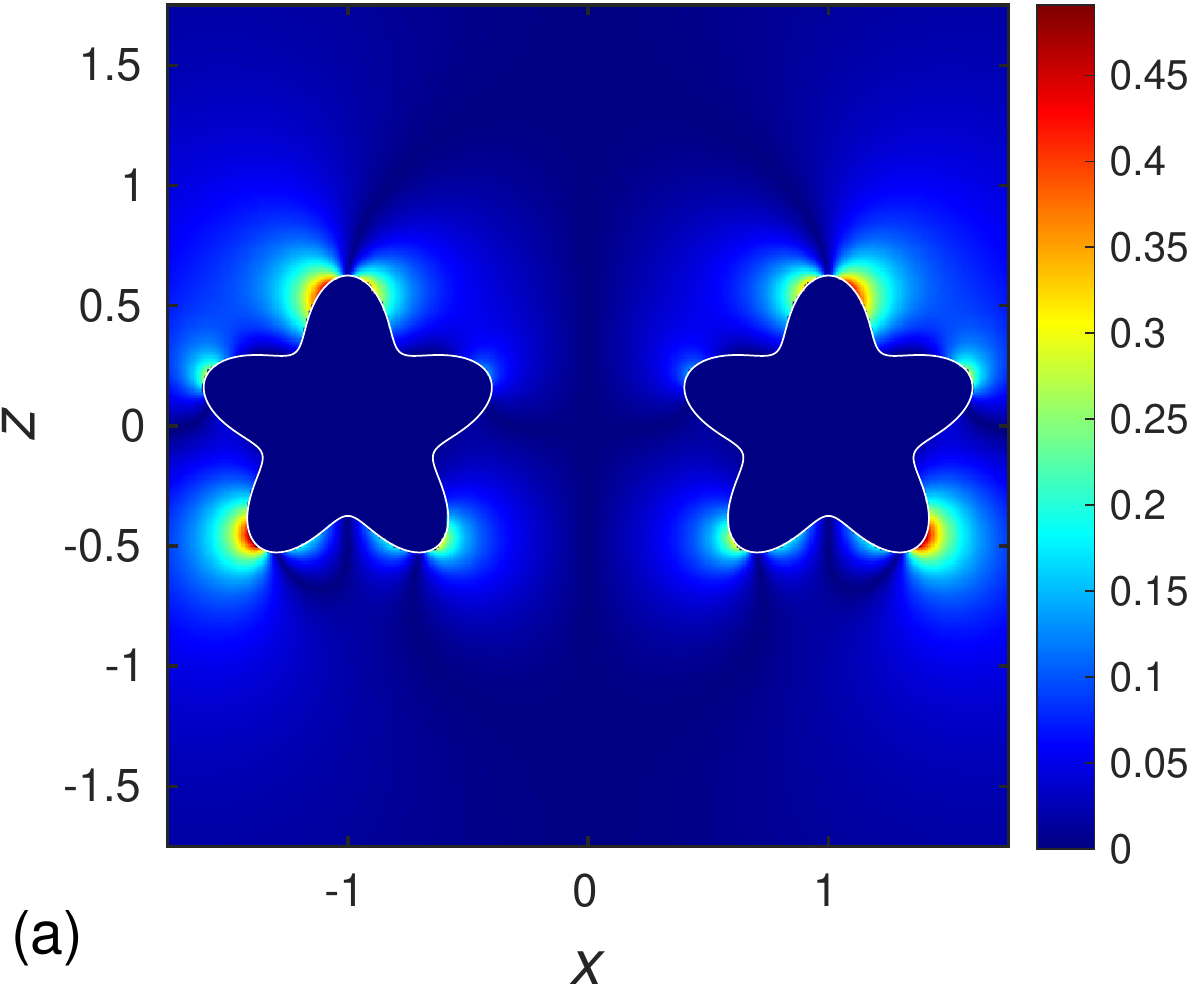}
\includegraphics[height=50mm]{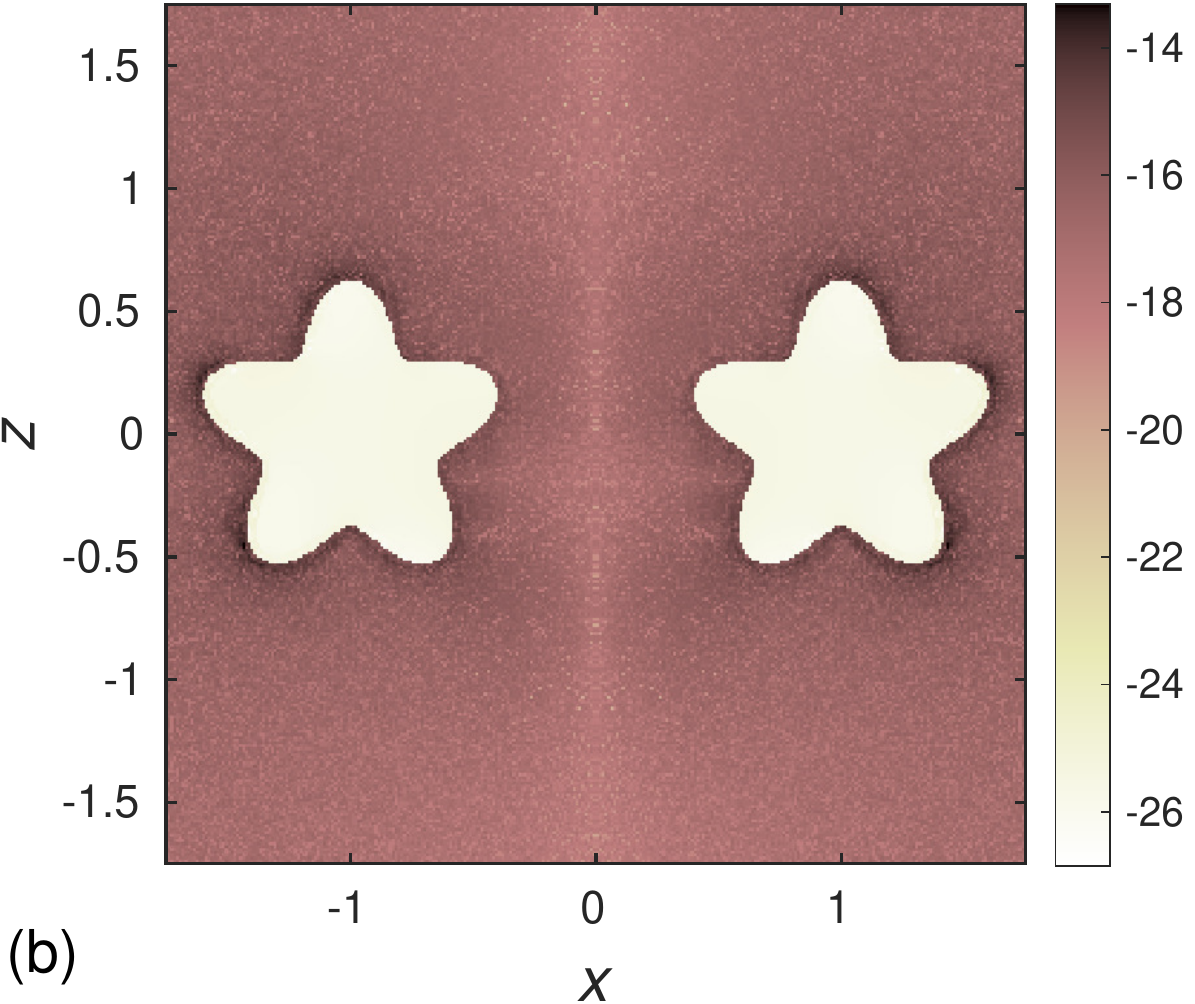}

\vspace{1mm}
\includegraphics[height=50mm]{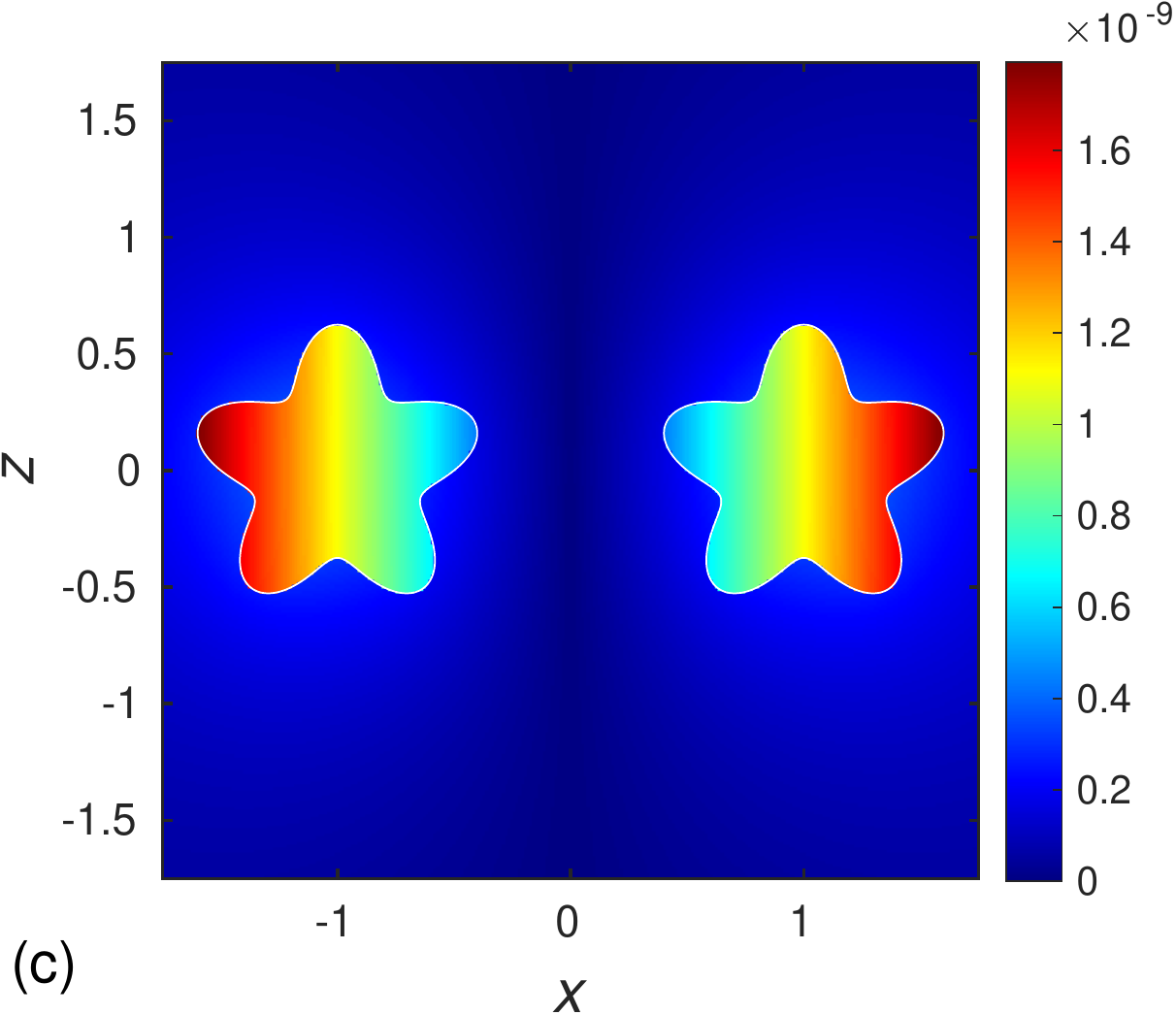}
\includegraphics[height=50mm]{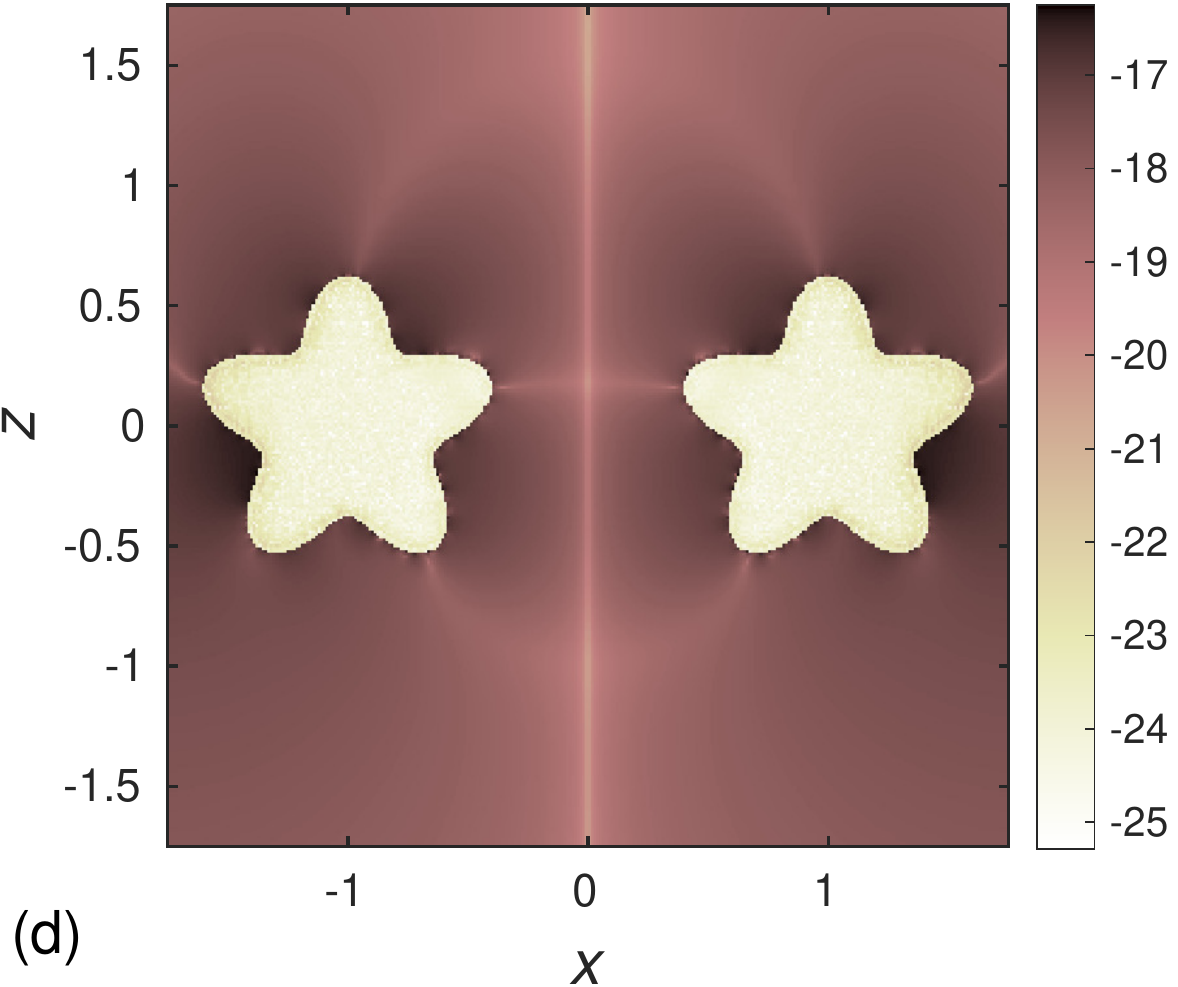}

\vspace{1mm}
\includegraphics[height=50mm]{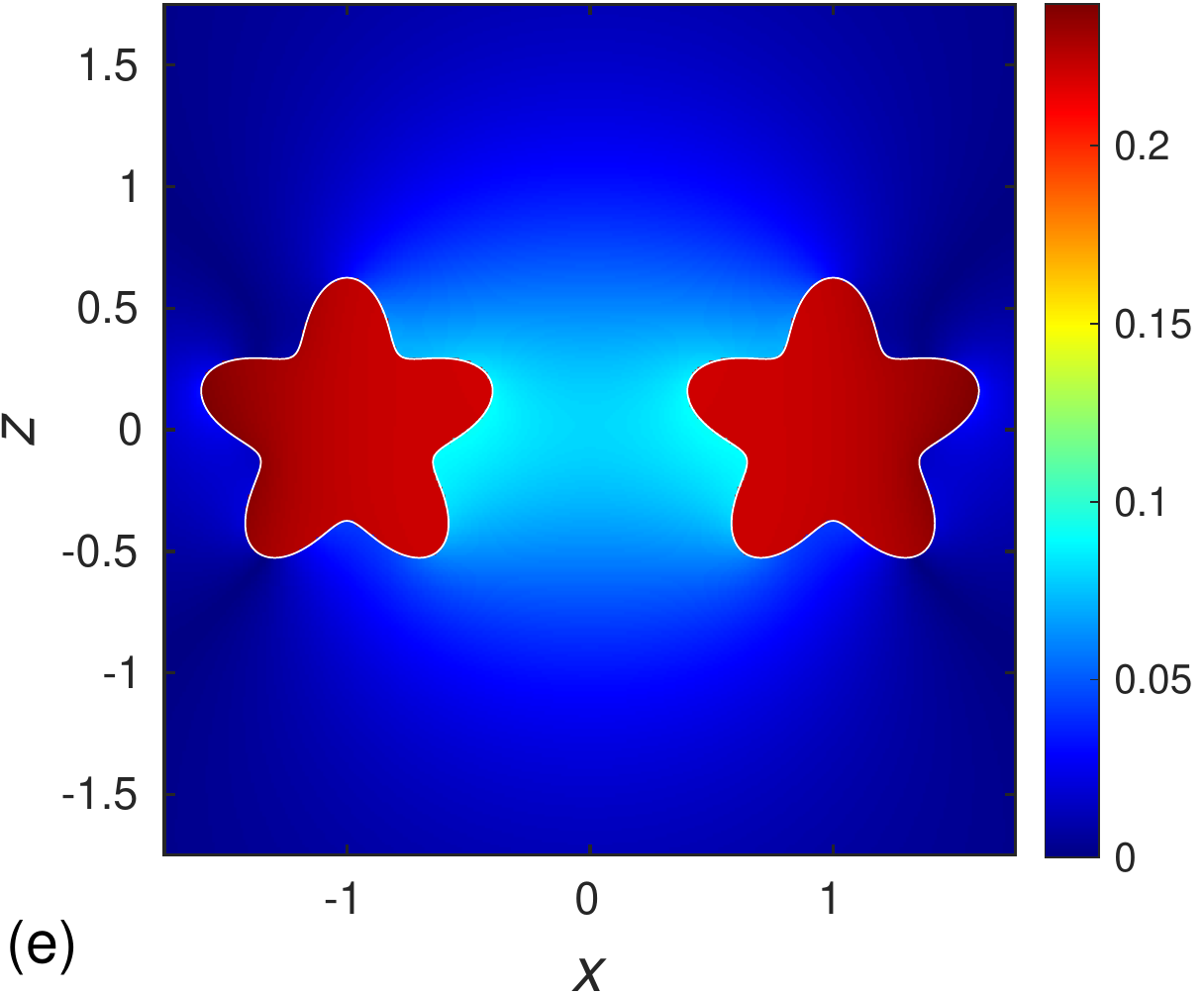}
\includegraphics[height=50mm]{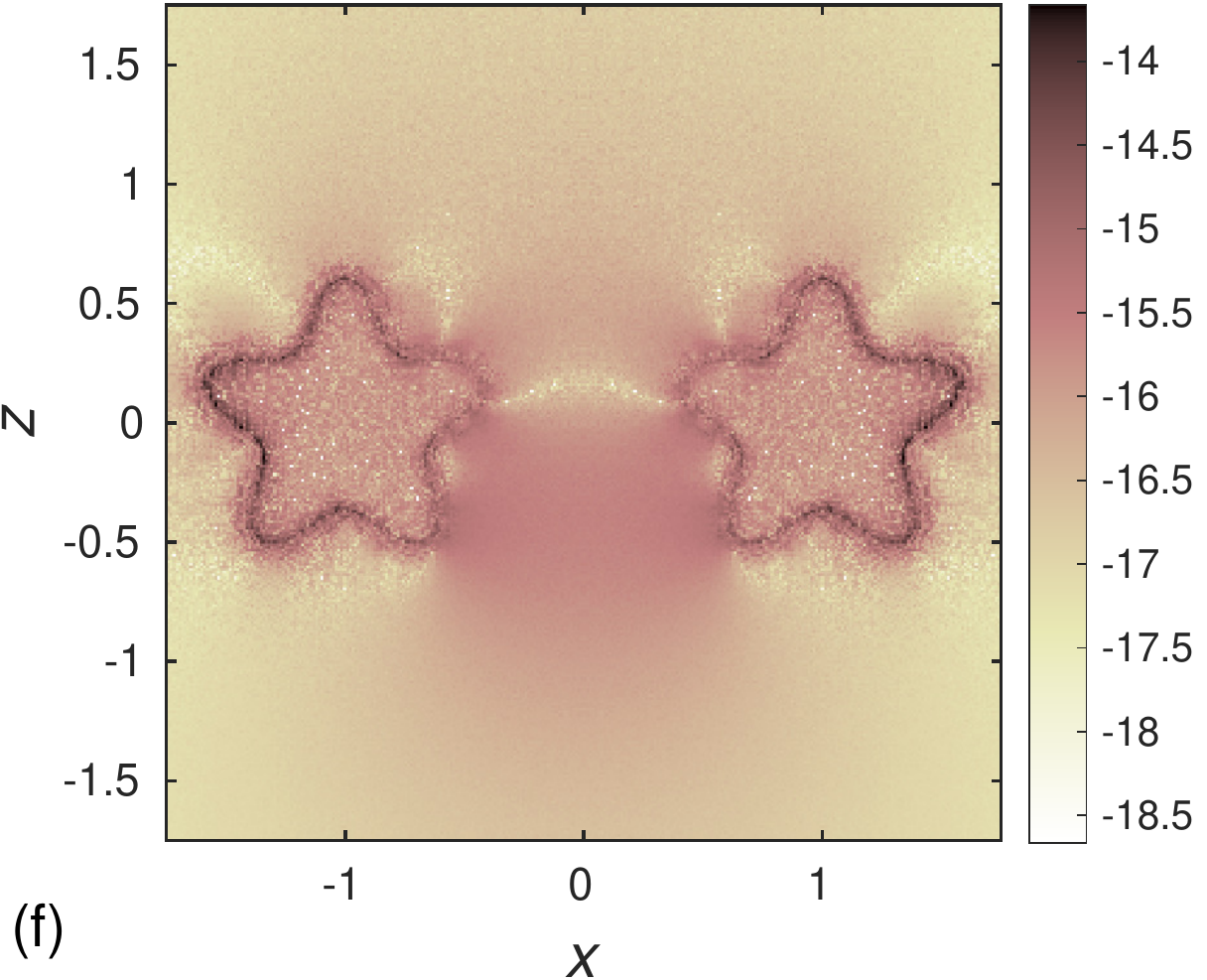}
\caption{\sf Field images for scattering of \eqref{eq:partialwave} by the
``starfish torus''~\eqref{eq:starfishtorus} at $k_-=10^{-8}, k_+=1+i$;
(a) scattered/transmitted amplitude $|E_\rho|$;
(c) scattered/transmitted amplitude $|E_\theta|$;
(e) scattered/transmitted amplitude $|H_z|$.
(b,d,f) $\log_{10}$ of estimated absolute error of complex fields 
using (B-aug1), for (a,c,e) respectively.
}
\label{fig:highgenus1-partial}
\end{figure}

\begin{figure}[t!]
\centering
\includegraphics[height=50mm]{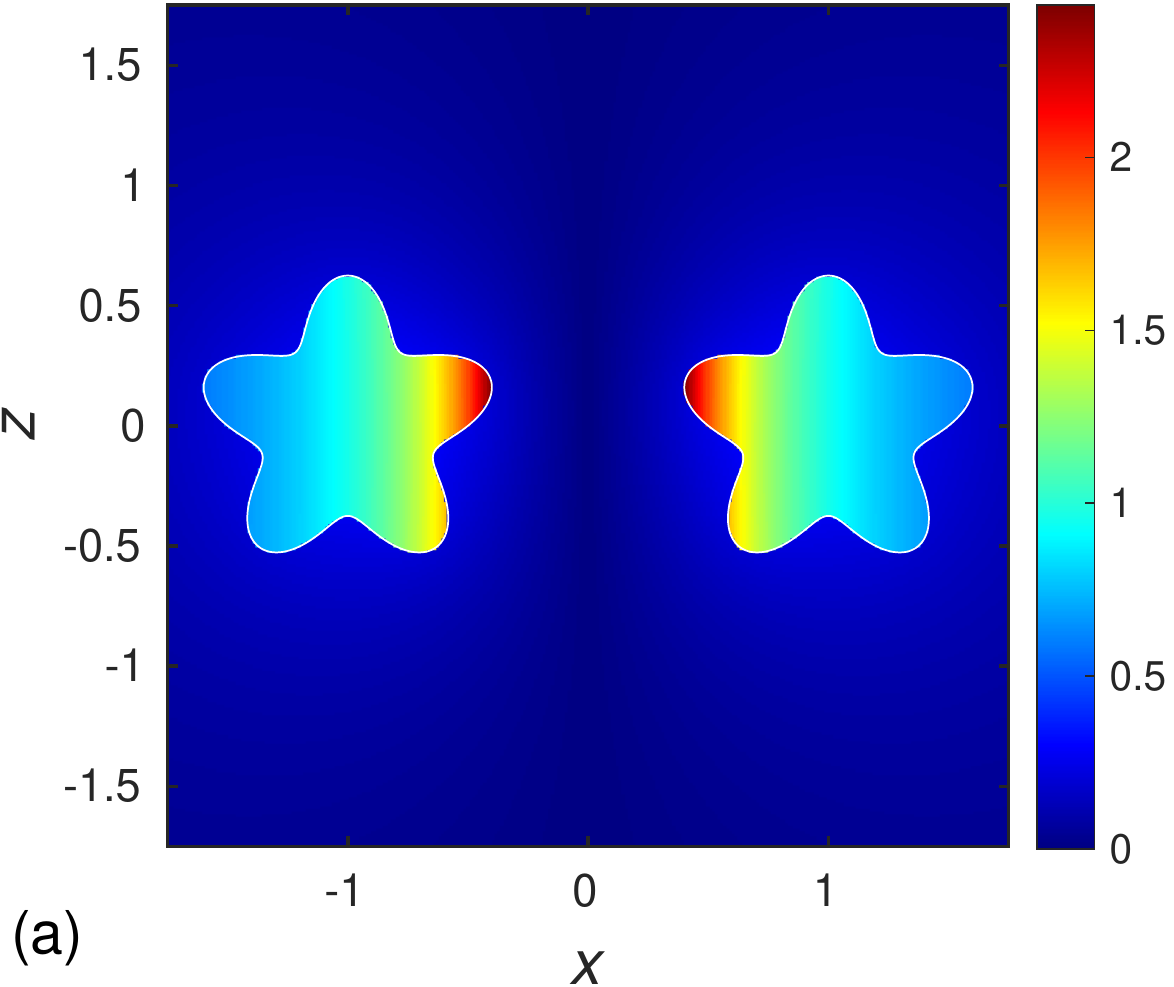}
\includegraphics[height=50mm]{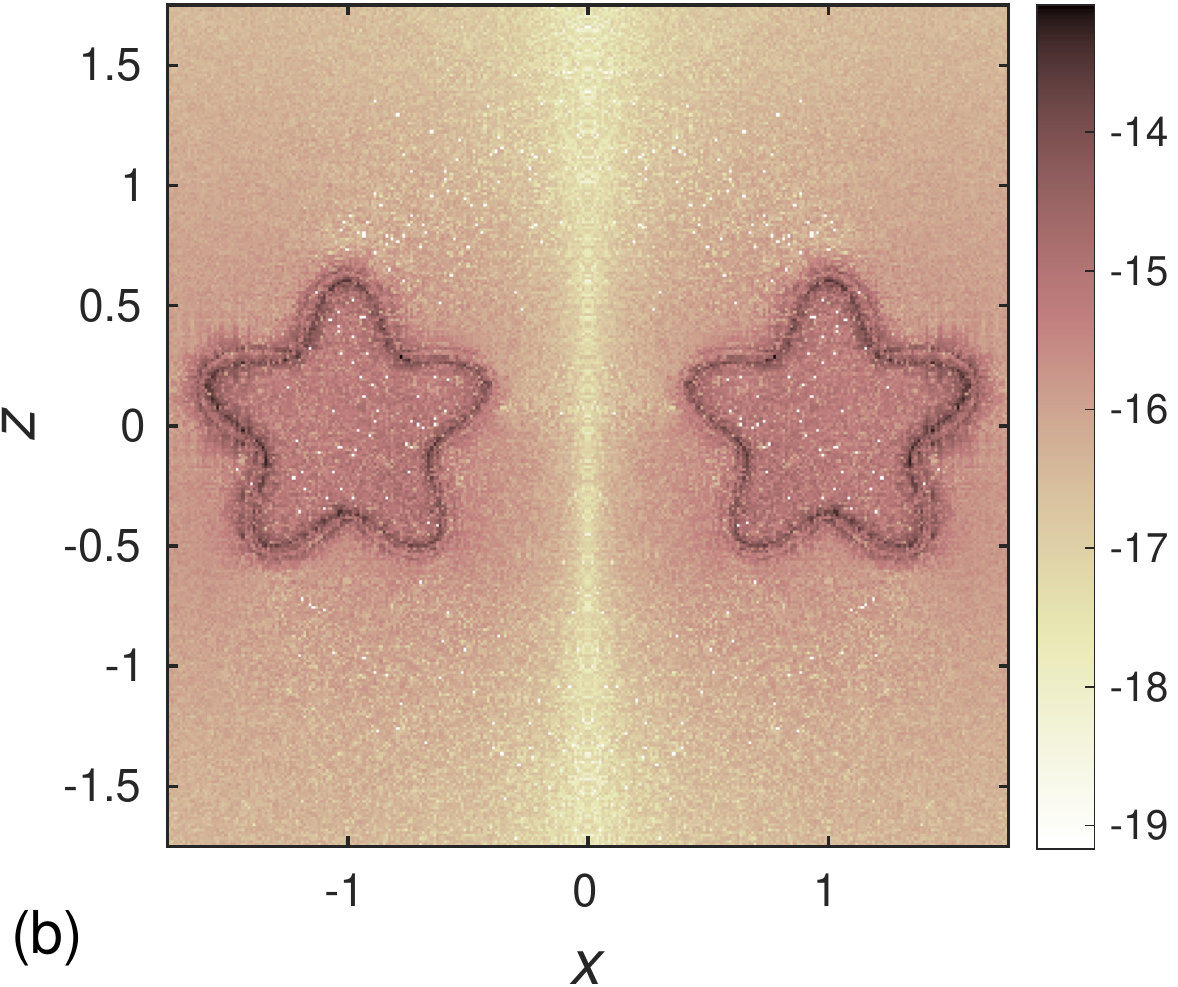}

\vspace{1mm}
\includegraphics[height=50mm]{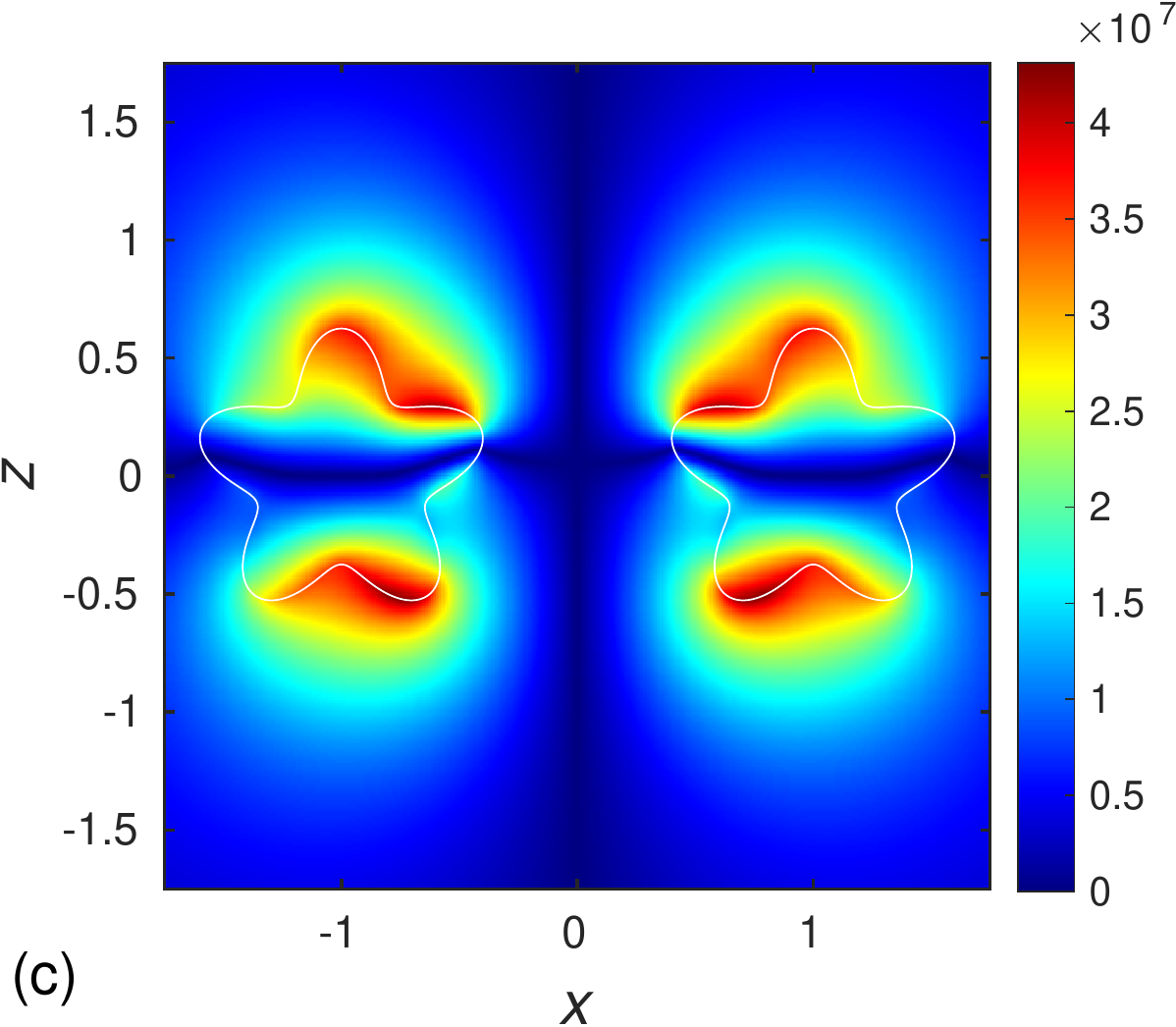}
\includegraphics[height=50mm]{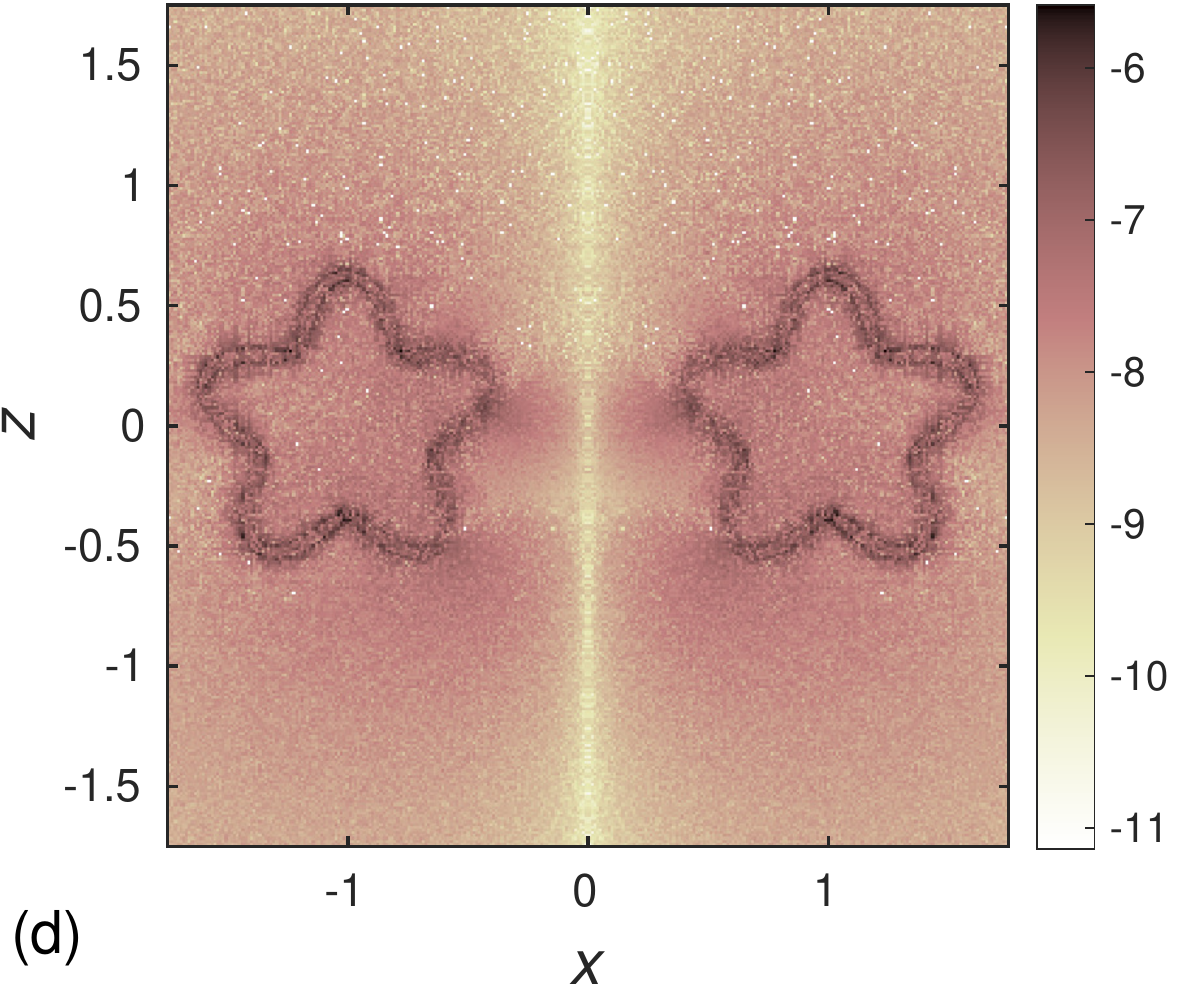}

\vspace{1mm}
\includegraphics[height=50mm]{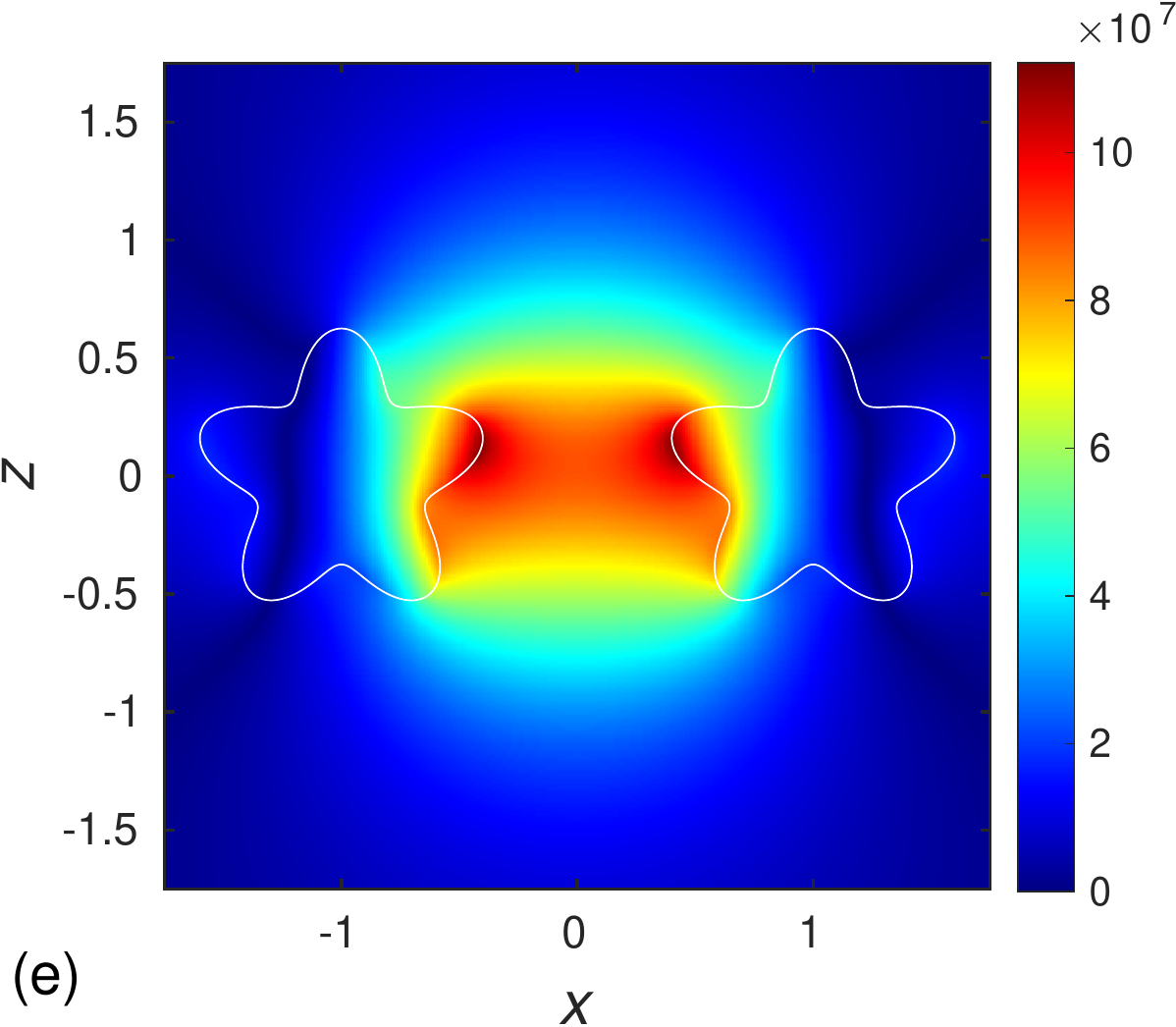}
\includegraphics[height=50mm]{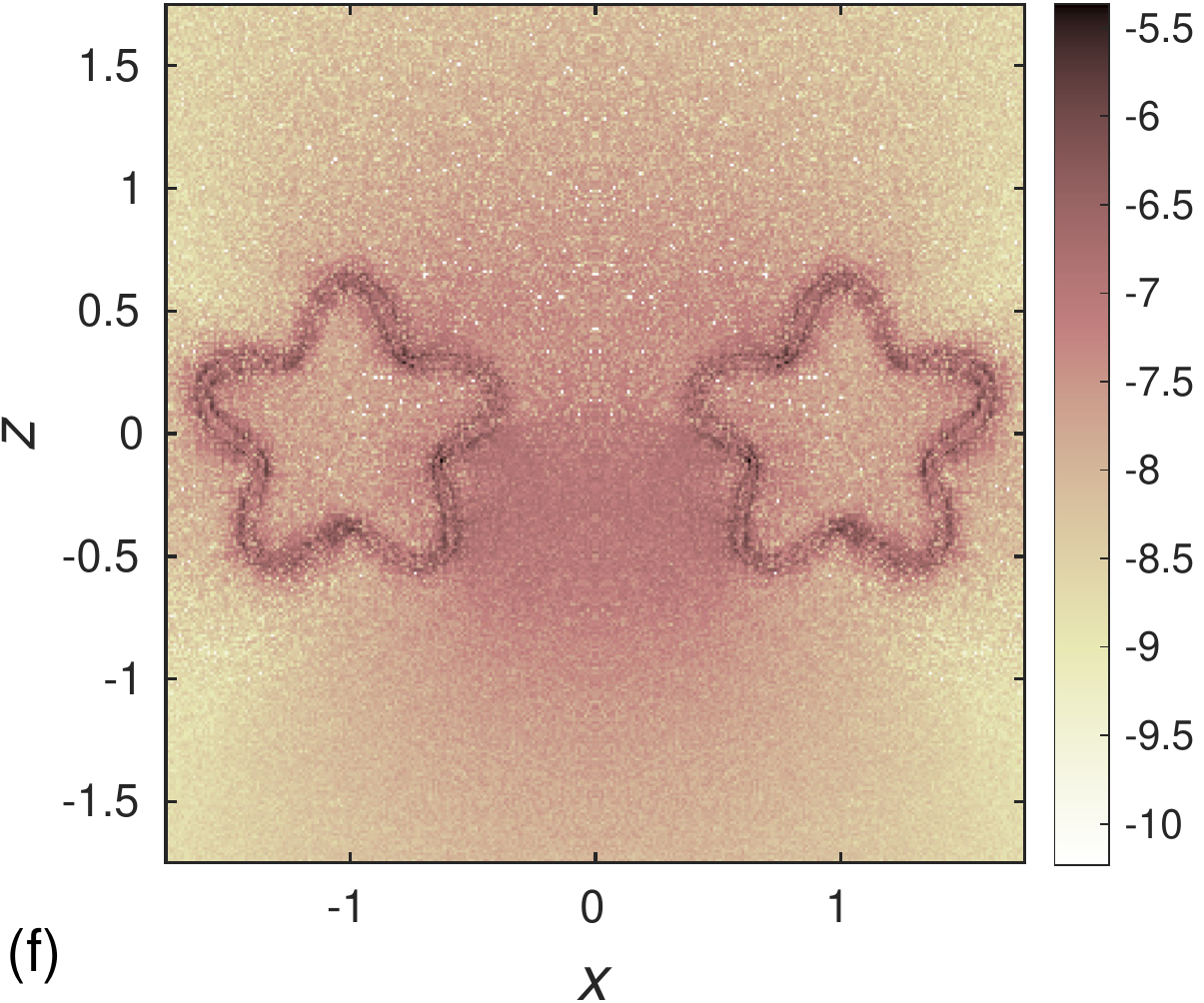}
\caption{\sf Field images for scattering of \eqref{eq:zcoil} by the
``starfish torus''~\eqref{eq:starfishtorus} at $k_-=10^{-8}, k_+=1+i$;
(a) scattered/transmitted amplitude $|E_\theta|$;
(c) scattered/transmitted amplitude $|H_\rho|$;
(e) scattered/transmitted amplitude $|H_z|$.
(b,d,f) $\log_{10}$ of estimated absolute error of complex fields 
using (A$\infty$-aug), for (a,c,e) respectively.
}
\label{fig:highgenus1-z}
\end{figure}

We consider the ``starfish torus''  with $L\approx 3.2$ cm  defined by
\eqref{eq:starfishtorus}, again 
at wavenumbers  $k_-= 10^{-8}$ ${\rm cm}^{-1}$ and $k_+= 1+i$ ${\rm cm}^{-1}$. 

   First, we use the incident field \eqref{eq:partialwave}, for which
   $|d^1_N f^0|/\max_\Gamma |f^0|\approx 0.4$. This indicates that
   \eqref{eq:partialwave} does not excite the Neumann eigenfield and
   motivates using Dirac (B-aug1). The number of accurate digits
   obtained for  $\{E^+, E^-, H^+, H^-\}$  are $\{13, 13, 13,
   14\}$, and GMRES needs $37$ iterations.  Solving the single block
   system~\eqref{eq:2ndweq}, which is needed only once for a given
   $\Gamma$, requires $20$ iterations.  The amplitudes of a
   selection of components of the fields $E^\pm$ and $H^\pm$ are shown
   in Figure~\ref{fig:highgenus1-partial} along with field errors.
   Qualitatively, the result is similar to that in
   Section~\ref{sec:genus0}.

For comparison, solving the same scattering problem with Dirac
   (A$\infty$-aug) also takes $37$ iterations but gives $\{6, 13, 7,
   7\}$ accurate digits for  $\{E^+, E^-, H^+, H^-\}$. 
This gives numerical support for choosing (B-aug1) for eddy current
   scattering with surfaces of genus $1$ when the Neumann eigenfield is
   not  excited. This is so  since (A$\infty$-aug) computes the fields at
   the wrong scale \eqref{eq:Ascale}, with subsequent loss of
   accuracy.

Second, we use the incident field \eqref{eq:zcoil}, for which
   $|d^1_N f^0|/\max_\Gamma |f^0|\approx 6\cdot 10^7$. This indicates
   that  \eqref{eq:zcoil}  does excite the Neumann eigenfield, and
   motivates using Dirac (A$\infty$-aug). The number of accurate digits
   obtained for  $\{E^+, E^-, H^+, H^-\}$  are $\{13, 14, 13,
   13\}$ and GMRES needs $24$ iterations. The amplitudes of a selection
   of components of the fields $E^\pm$ and $H^\pm$ are shown in
   Figure~\ref{fig:highgenus1-z} along with field errors. It is clearly
   seen that the incident field~\eqref{eq:zcoil} does excite the
   Neumann field. The scattered and transmitted fields are very similar
   to those shown in Figure~\ref{fig:Neumanneig}(g,h,i), except for the
   scale. Note that $k_-^{-1}E^+$ and $H^\pm$ are of order $10^8$,
   which is the scale \eqref{eq:Bscale} that (A$\infty$-aug) is adapted
   to, and a factor of $10^8$ larger than the scale of a generic field
   as in the discussion beginning Section~\ref{sec:physics}. The
   components $E_\rho$, $E_z$, and $H_\theta$ (not shown) are of order
   $10^{-16}$, $10^{-16}$, and $10^{-9}$.

For comparison, solving the same scattering problem with Dirac (B-aug1)
takes $32$ iterations and gives $\{13, 7, 13, 13\}$ accurate digits for 
 $\{E^+, E^-, H^+, H^-\}$. 
This gives numerical support for choosing (A$\infty$-aug) for eddy
   current scattering with surfaces of genus $1$ when the Neumann
   eigenfield is  excited. 
To see why (B-aug1) gives loss of accuracy in $E^-$ in
this example, note that the right-hand side in \eqref{eq:B1system}
is of order $10^8$ due to $d^1_Nf^0$ in the second term.
Since the system in \eqref{eq:B1system} is  well-conditioned 
and has norm of order $1$,
the density $h$ will also be of order $10^8$.
When finally computing the fields with \eqref{eq:projdensB1} then,
according to \eqref{eq:Bscale}, if no cancellation occurs 
$E^-$ will be of order $10^8$.
But as we see in Figure~\ref{fig:highgenus1-z}(a) $E^-$ is of order $1$,
and this is due to cancellation in 
\eqref{eq:projdensB1} which leads to the loss of accuracy.

\begin{figure}[t!]
\centering
\includegraphics[height=50mm]{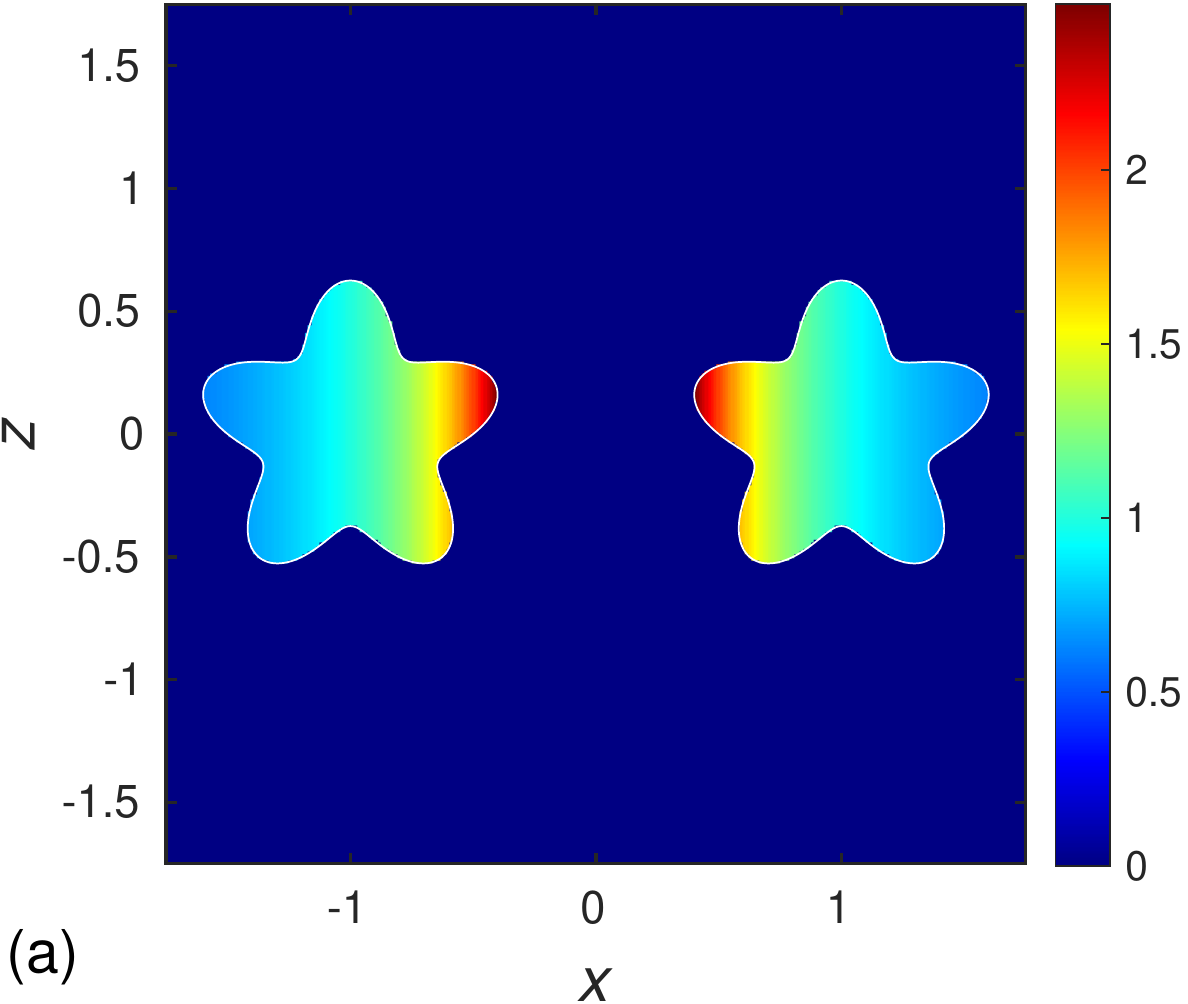}
\includegraphics[height=50mm]{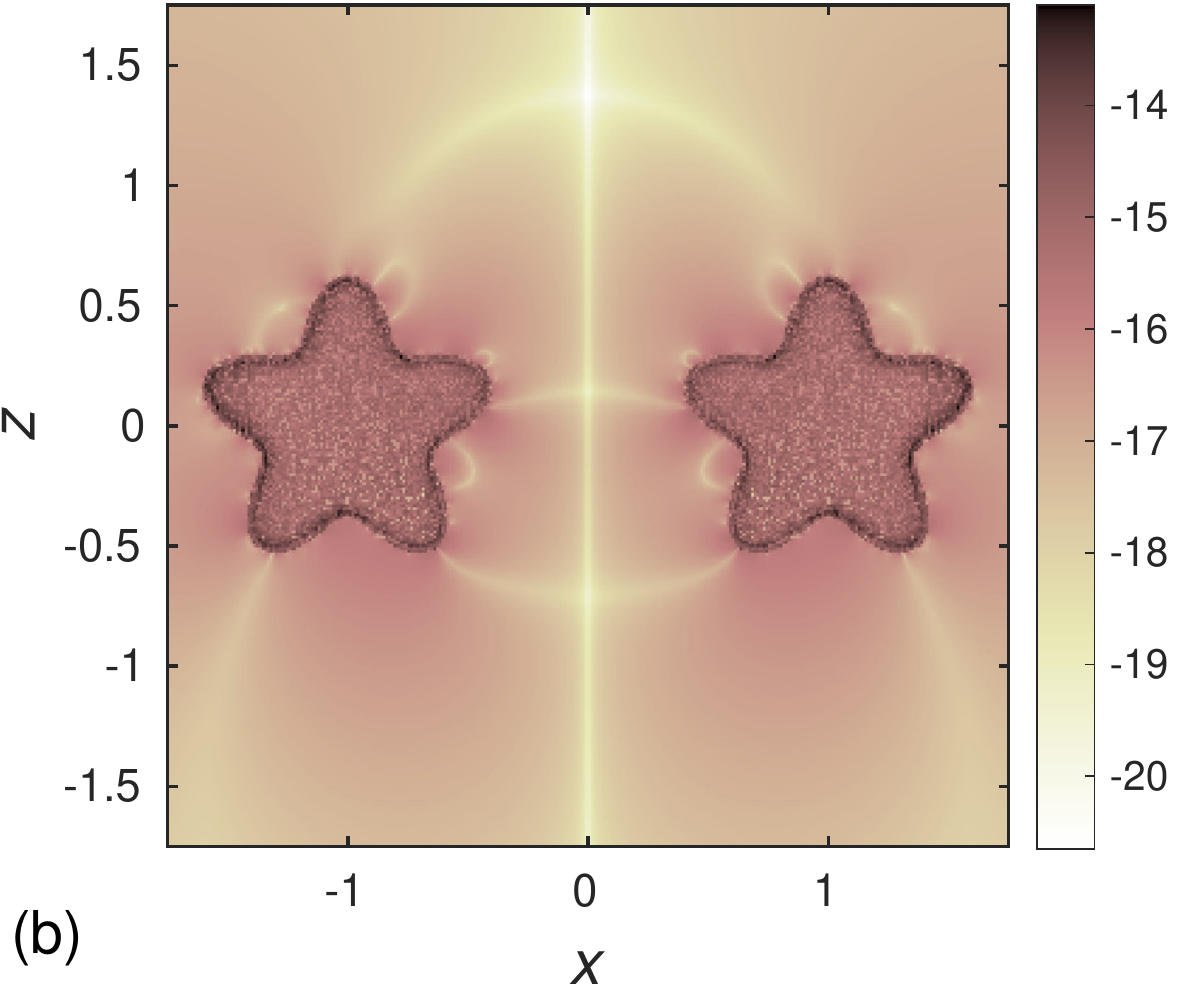}

\vspace{1mm}
\includegraphics[height=50mm]{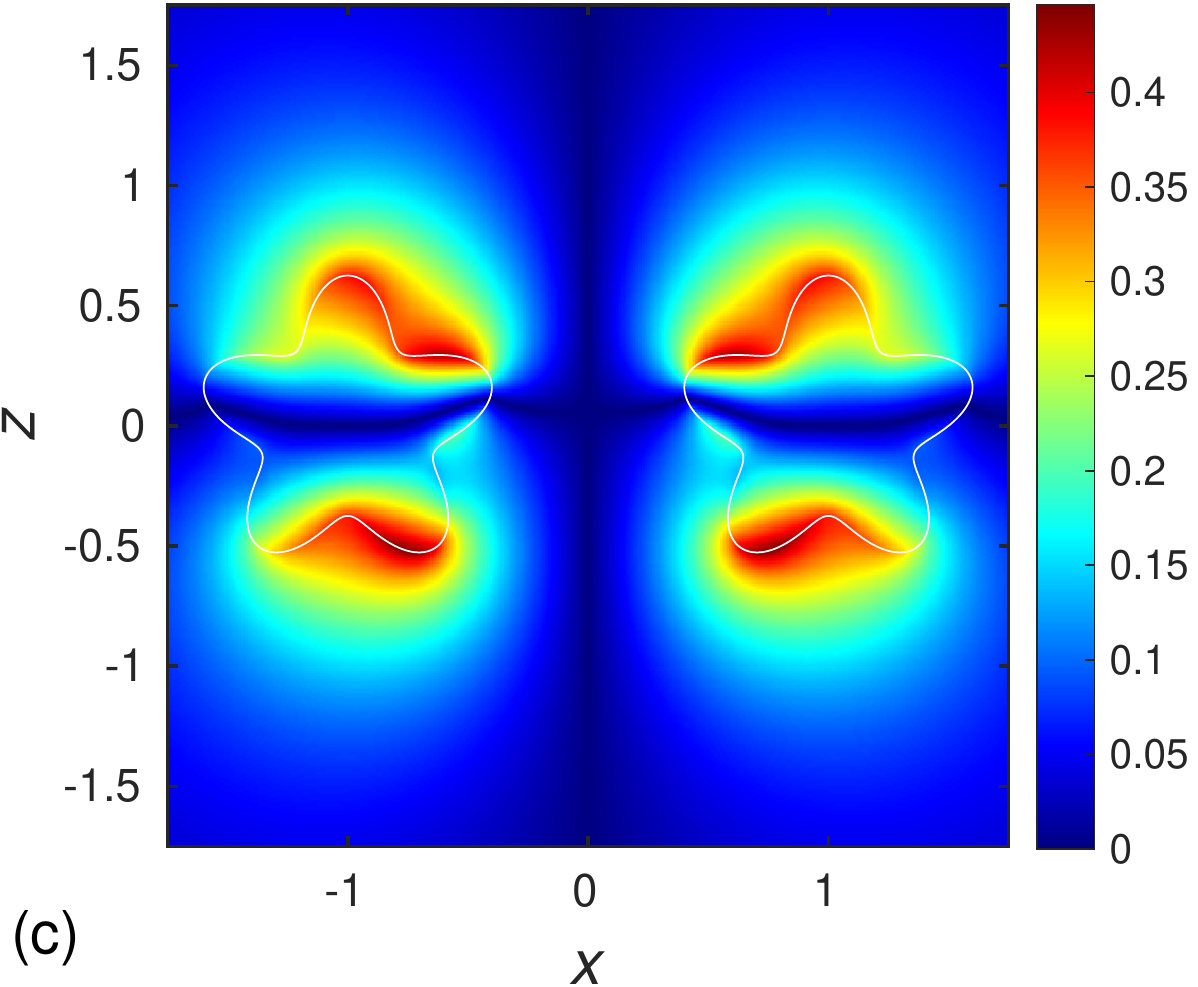}
\includegraphics[height=50mm]{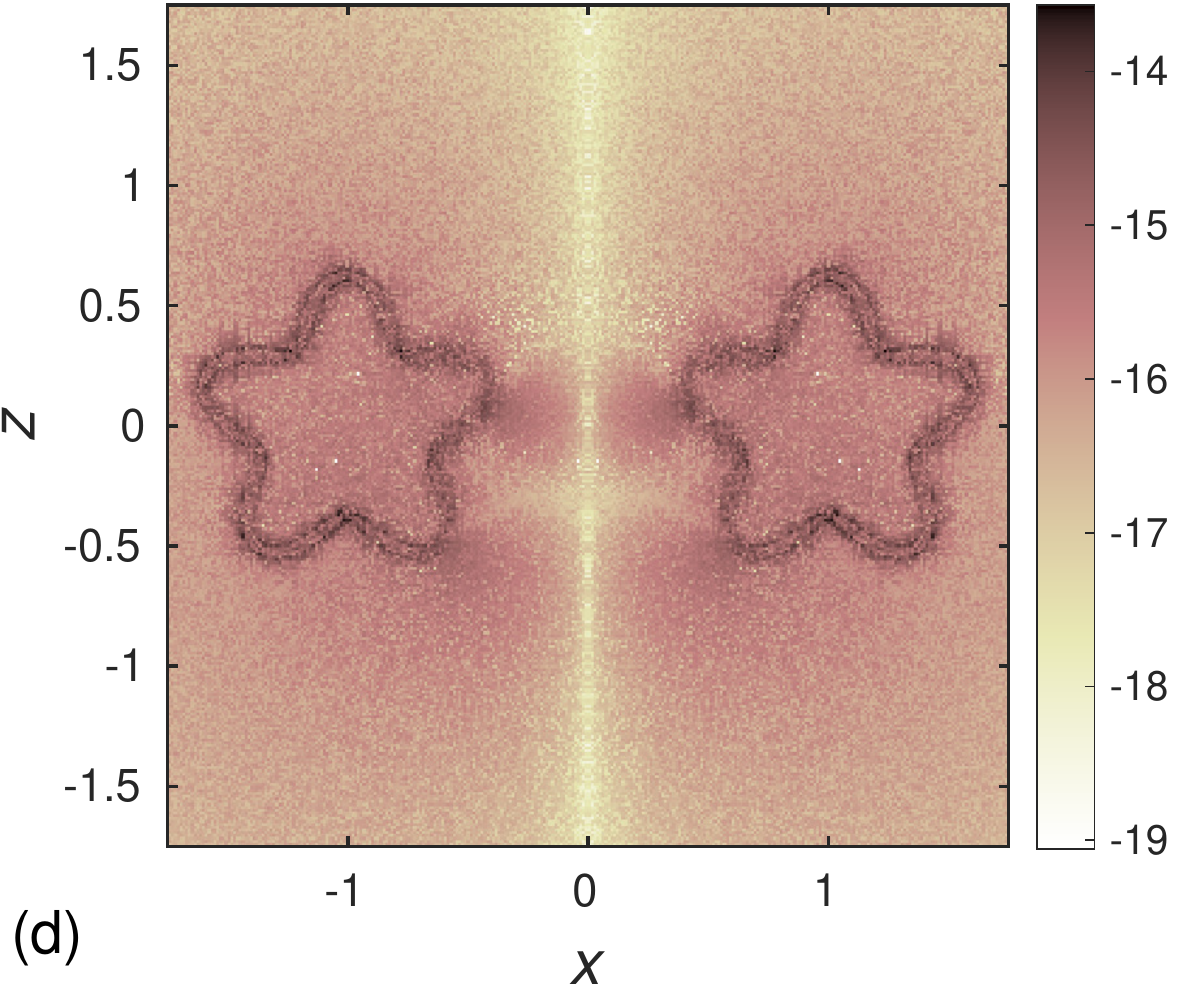}

\vspace{1mm}
\includegraphics[height=50mm]{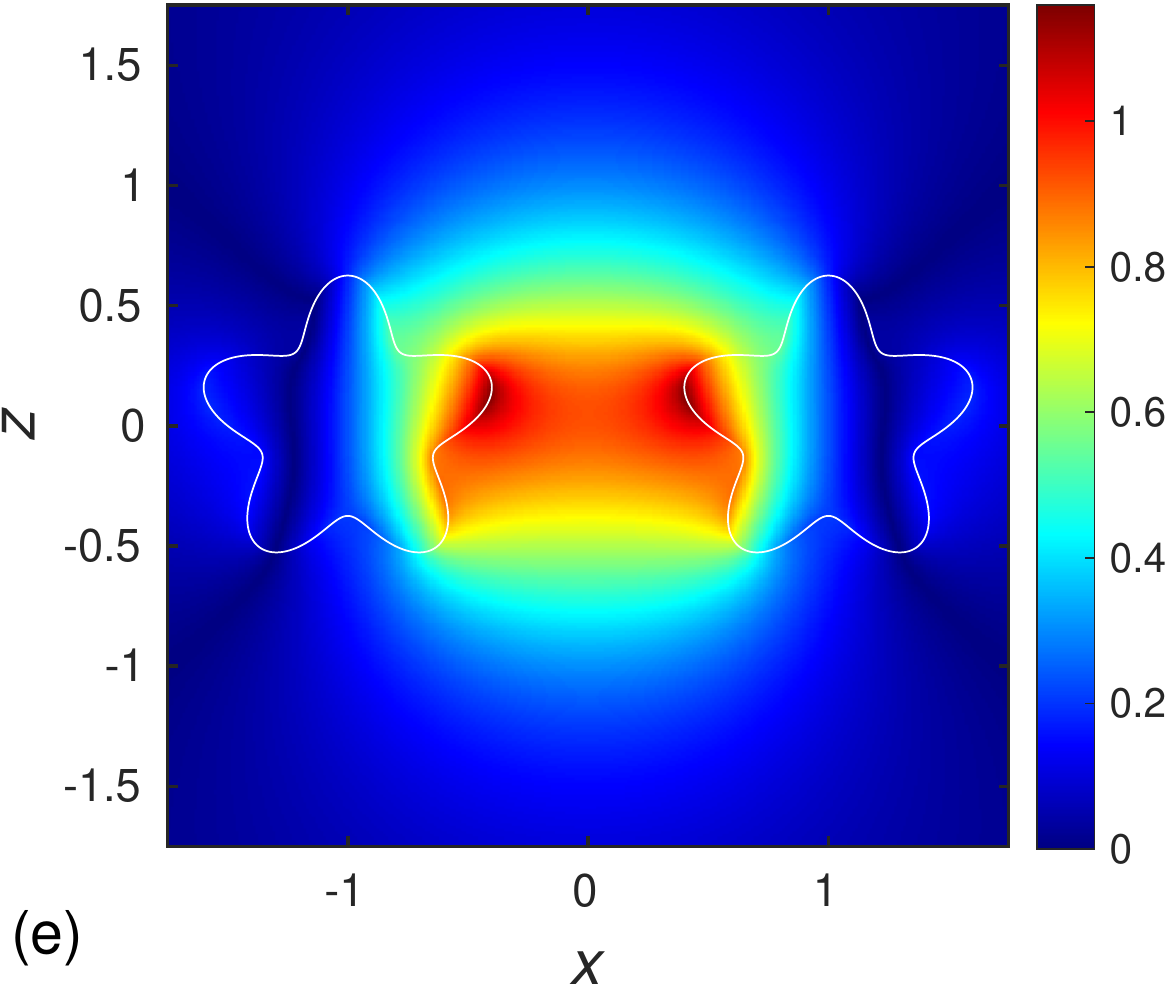}
\includegraphics[height=50mm]{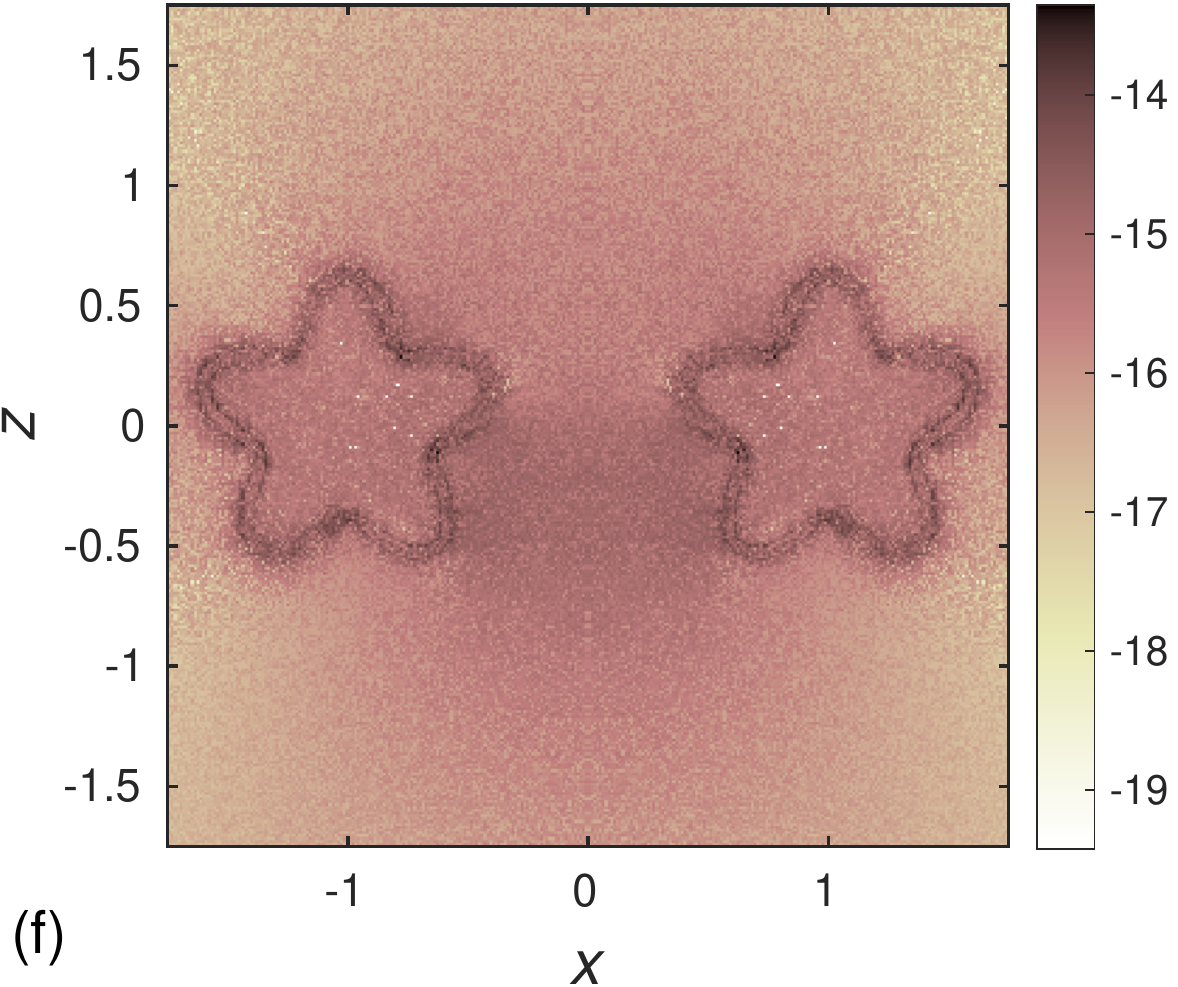}
\caption{\sf 
Field images for scattering of \eqref{eq:zcoil} by the
``starfish torus''~\eqref{eq:starfishtorus} at $k_-=10^{-8}, k_+=10^{-4}(1+i)$;
(a) scattered/transmitted amplitude $|E_\theta|$;
(c) scattered/transmitted amplitude $|H_\rho|$;
(e) scattered/transmitted amplitude $|H_z|$.
(b,d,f) $\log_{10}$ of estimated absolute error of complex fields 
using (A$\infty$-aug), for (a,c,e) respectively.
}
\label{fig:finitegenus1-z}
\end{figure}

\subsection{The  medium  conductivity ``starfish torus''}  \label{sec:finitetorus}

   We consider the field \eqref{eq:zcoil} incident on the ``starfish
   torus''~\eqref{eq:starfishtorus}  with $L\approx 3.2$~cm,  now
   at wavenumbers  $k_-= 10^{-8}$ ${\rm cm}^{-1}$ and $k_+=
   10^{-4}(1+i)$ ${\rm cm}^{-1}$.  This corresponds to a frequency
    $\omega\approx 300$ rad/s  and a conductivity $\sigma \approx
   0.53$~S/m, which for example occurs for seawater. Furthermore,
   $|d^1_N f^0|/\max_\Gamma |f^0|\approx 4\cdot 10^7$, which indicates
   that \eqref{eq:partialwave} does excite the Neumann eigenfield and
   motivates using Dirac (A$\infty$-aug). The number of accurate digits
   obtained for  $\{E^+, E^-, H^+, H^-\}$  are $\{13, 15, 13,
   13\}$ and GMRES needs $16$ iterations. The amplitudes of a selection
   of components of the fields $E^\pm$ and $H^\pm$ are shown in
   Figure~\ref{fig:finitegenus1-z} along with field errors. This result
   is very similar to Figure~\ref{fig:highgenus1-z}, with the
   difference that $J$ and $H$ are a factor of $10^8$ smaller now.
   However, since the conductivity is also a factor of $10^8$ smaller,
   the transmitted electric field $E^+$ has barely changed and the
   fields $E^+$ and $H^\pm$ are still a factor of $10^8$ larger than
   what is expected in a generic scattering situation, according to the
   discussion beginning Section~\ref{sec:physics}. The important point
   that we want to make is that, although none of the
   transmitted and scattered fields $E^\pm$, $H^\pm$ are significantly
   larger than the incident fields $E^0$, $H^0$, we are looking at an
   excited Neumann eigenfield,  according to \eqref{eq:defneigfield}. 
   Indeed, in the generic scattering
   situation when the eigenfield is not excited, the magnitude of $E^+$
   would be of order $10^{-8}$. 

   For comparison, solving the same scattering problem with Dirac
   (B-aug1) takes $27$ iterations and gives $\{13, 8, 13, 13\}$
   accurate digits for  $\{E^+, E^-, H^+, H^-\}$.  Again this
   gives numerical support for choosing (A$\infty$-aug) for eddy
   current scattering with surfaces of genus $1$ when the Neumann
   eigenfield is excited.

\subsection{Condition numbers for systems and field representations}
\label{sec:conditioning}

\begin{figure}[t!]
\centering
\includegraphics[height=50mm]{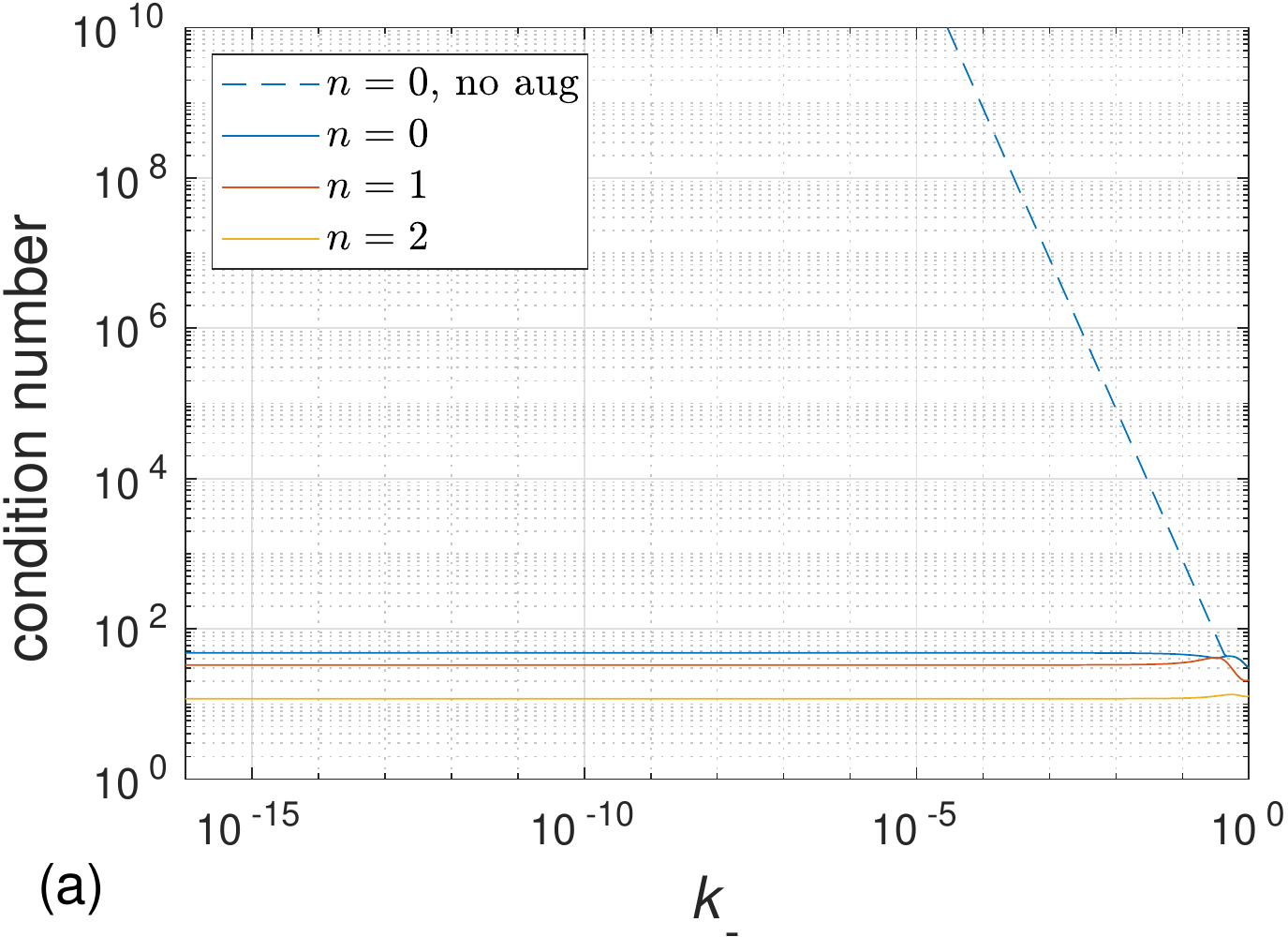}
\includegraphics[height=50mm]{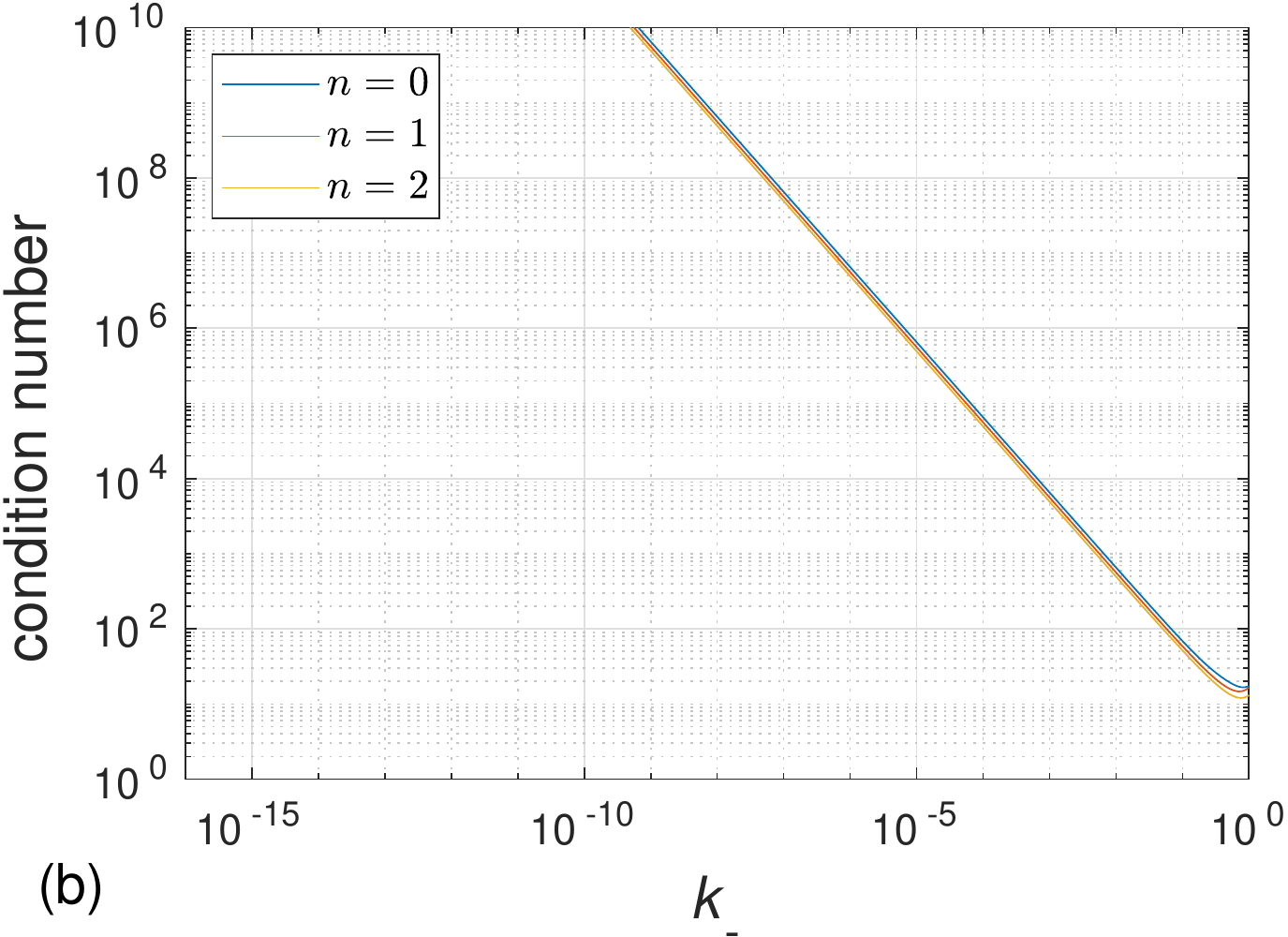}

\vspace{1mm}
\includegraphics[height=50mm]{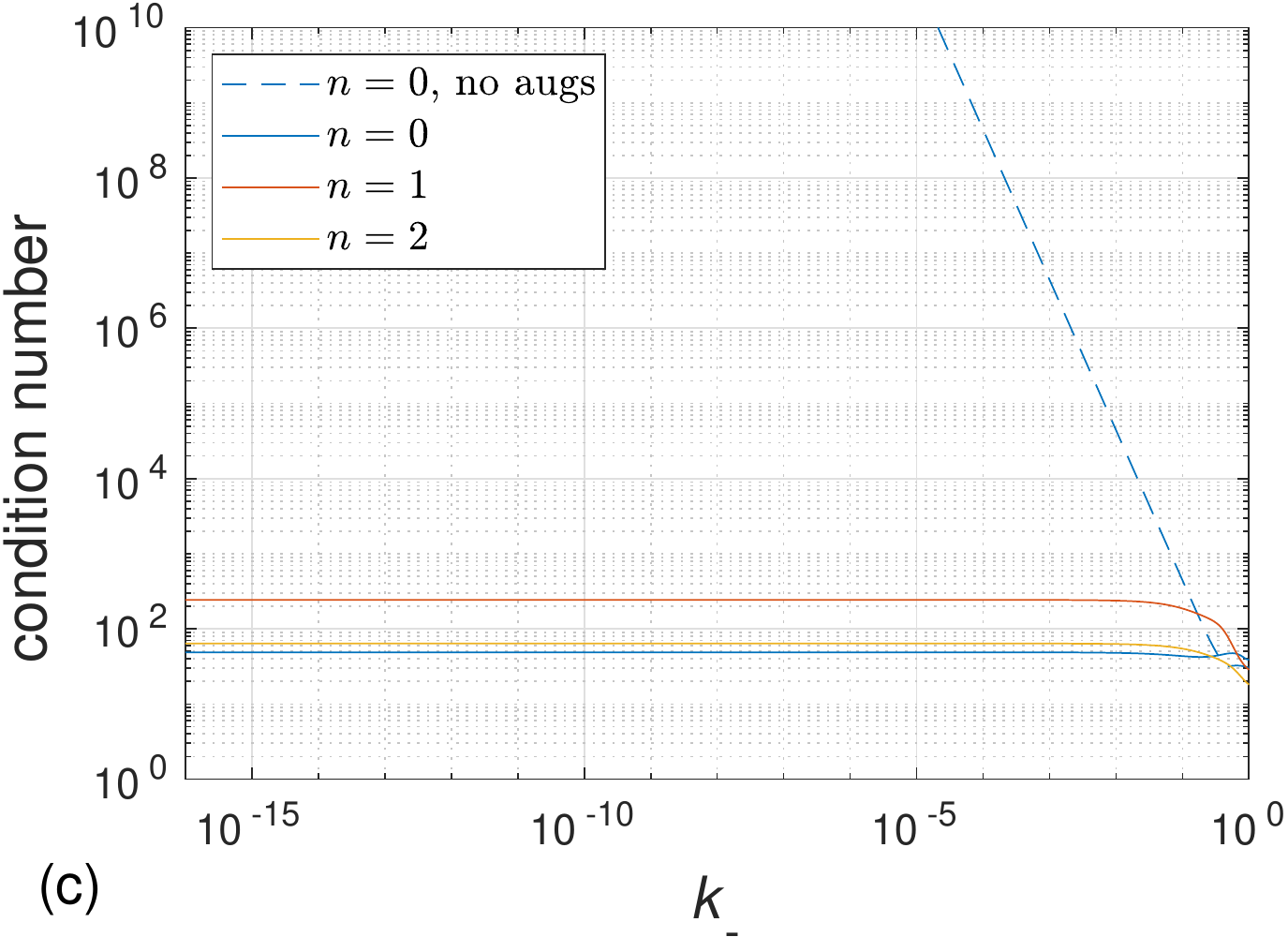}
\includegraphics[height=50mm]{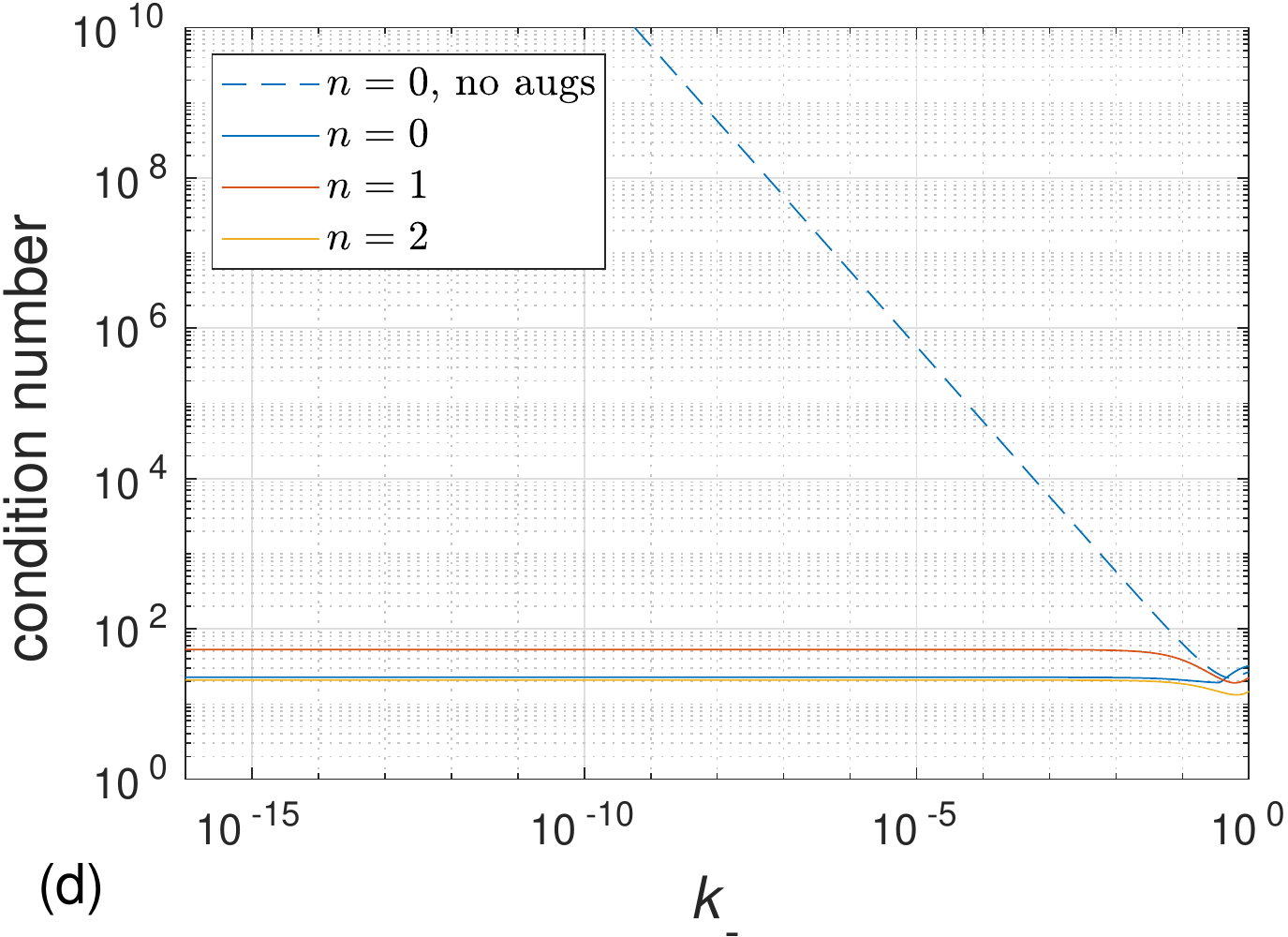}

\caption{\sf The ``starfish torus''~\eqref{eq:starfishtorus} at
high conductivities $k_-\in [10^{-16}, 1]$, $k_+=1+i$.
First row: condition numbers for (A$\infty$-aug);
(a) system \eqref{eq:Ainfty-augsystem};
(b) field representation \eqref{eq:projdens}.
Second row: condition numbers for (B-aug1);
(c) system \eqref{eq:B1system};
(d) field representation \eqref{eq:projdensB1}.
}
\label{fig:conditioning}
\end{figure}
 
We  here examine Dirac (A$\infty$-aug) and Dirac (B-aug1) from a
   condition number point of view. 
Recall from the discussion in the
Introduction that the computation of the transmitted and 
scattered fields involves (a) solving a linear system for the 
density $h$, followed by (b) applying the field formulas
to $h$.
The important point that we want to stress is that condition numbers
for the system (a) alone give  insufficient information for assessing 
a BIE. 
Indeed, Figure~\ref{fig:conditioning}(a) shows that Dirac (A$\infty$-aug),
after augmentation, has a  well-conditioned  system.
But we have seen that (A$\infty$-aug) in general only computes the fields
accurately when the Neumann eigenfield is excited.
Figure~\ref{fig:conditioning}(b) reveals that the important missing
information is that the field representation \eqref{eq:projdens} is ill-conditioned
for (A$\infty$-aug), for all  modes $n$.  
We see a low-frequency breakdown in this field
representation since it fails to be a Fredholm map as $k_-\to 0$. 
This is unavoidable since 
the only way that (A$\infty$-aug) can compute fields much smaller
   than the Neumann eigenfield is by cancellation in the field
   evaluations. 
To be precise, the condition numbers in Figure~\ref{fig:conditioning}(b,d)
refer to the map
\begin{equation}   \label{eq:genericscaling}
   h\mapsto (\hat k^2/\japsigma E^+|_\Gamma, E^-|_\Gamma, 
H^+|_\Gamma, H^-|_\Gamma),
\end{equation}
where we have scaled the fields by their generic size in the eddy current regime,
as discussed in Section~\ref{sec:physics}.

Figure~\ref{fig:conditioning}(c,d) shows that for (B-aug1), we have succeeded in
constructing a BIE where both the system and the field representation
are  well-conditioned,  after augmentation.
That this is possible for (B-aug0) and genus $0$ is perhaps less surprising, since
the MTP in this case is well-conditioned.
But for (B-aug1) and genus $1$, we recall that the MTP itself is ill-conditioned.
Our design of (B-aug1) is such that the Neumann eigenfield is hiding in the
preprocessing, the computation of the  right-hand  side in 
\eqref{eq:B1system}, and in particular in computing $d^1_Nf^0$. 
To be precise, in order to assess the efficiency of a BIE one must
  take into account  three computations: 
the preprocessing involved in computing the right-hand side $g$ in
   \eqref{eq:abstractIER}, the solution of the main linear system that
   produces the density $h$, and finally the postprocessing involved in
   computing the fields $F^\pm$. 
For our ill-conditioned MTP, it is clearly the best option to let 
the Neumann eigenfield appear only in the preprocessing, as
in (B-aug1). 
It should be noted that $d^1_Nf^0$ in general requires careful computation, since 
for general incident fields $f^0$ this integral involves cancellations.
However, for mode-$0$ fields like \eqref{eq:partialwave} and \eqref{eq:zcoil}
there are no such cancellations.

\begin{figure}[t]
\centering
\includegraphics[height=90mm]{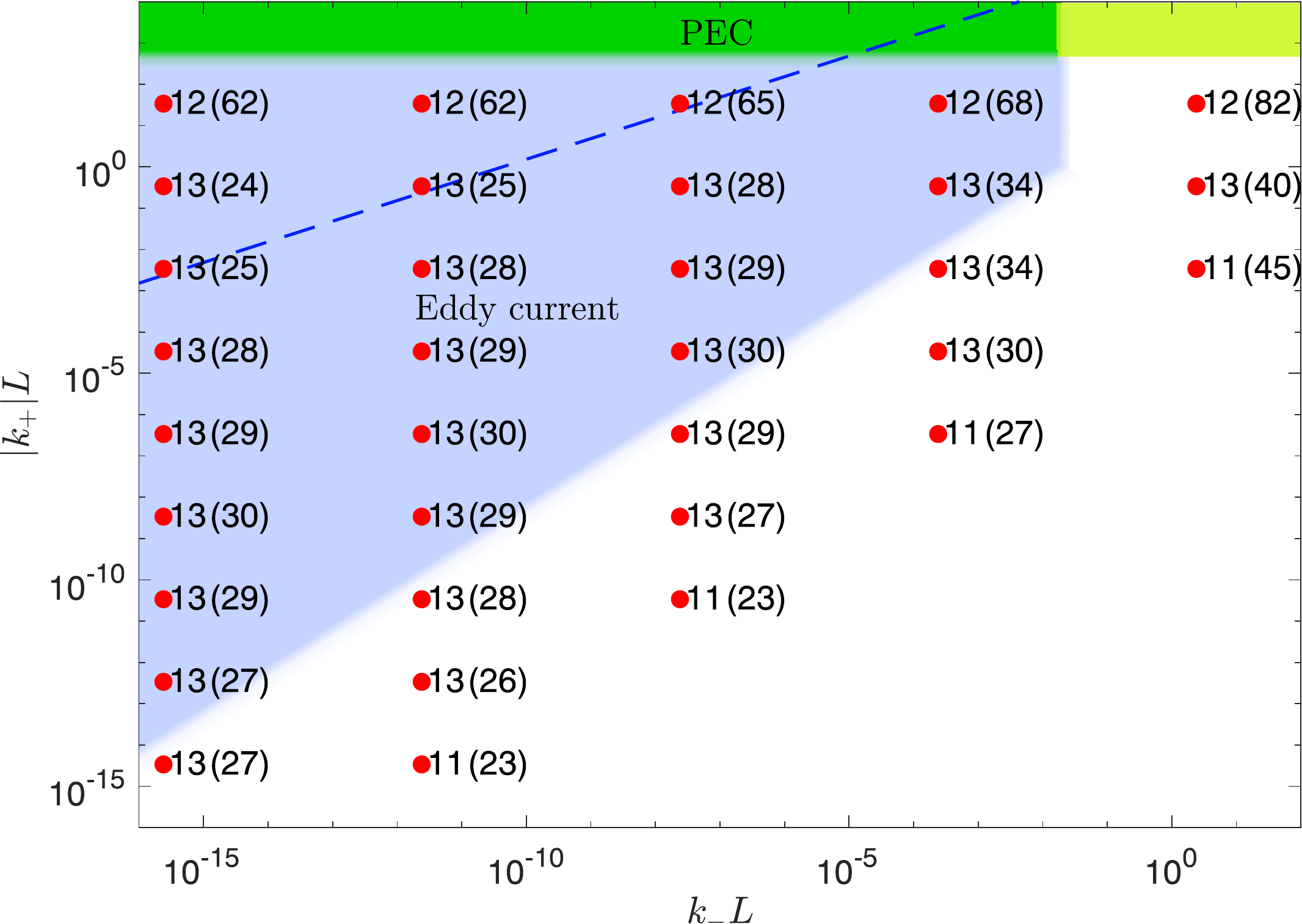}
\caption{\sf 
Performance in the eddy current regime of Dirac (B-aug0) on the
   ``rotated starfish''~\eqref{eq:rotstarfish} with incident partial waves~\eqref{eq:partialwave}.
Notation as in Figure~\ref{fig:quadrant0genus1}.
}
\label{fig:quadrant0genus0}
\end{figure}

\subsection{Performance of (B-aug0/1) in the eddy current regime}\label{sec:perform}

We  conclude this section by surveying  the accuracy and speed 
of  Dirac  (B-aug0/1) across the regime~\eqref{eq:goodcond}, with the incident field~ \eqref{eq:partialwave}, which we have seen 
does not excite the Neumann eigenfield.
We compute the minimum number of accurate digits, as defined in Section~\ref{sec:relerror},
in the four fields  $\{E^+, E^-, H^+, H^-\}$  at all
   $90,\!000$ field points in the computational domain. 
This minimum $Y$, at pairs of wavenumbers across the
   regime~\eqref{eq:goodcond}, is reported 
in Figure~\ref{fig:quadrant0genus1} for Dirac (B-aug1)
and the ``starfish torus''~\eqref{eq:starfishtorus},
and in Figure~\ref{fig:quadrant0genus0} for Dirac (B-aug0)
and the ``rotated starfish''~\eqref{eq:rotstarfish}.
Within parentheses is also reported in these figures, the
number of iterations $X$ that it takes GMRES to compute
the density $h$.
We conclude that there is no low-frequency breakdown for  Dirac 
(B-aug0/1) in the regime~\eqref{eq:goodcond}.

\begin{figure}[t]
\centering
\includegraphics[height=90mm]{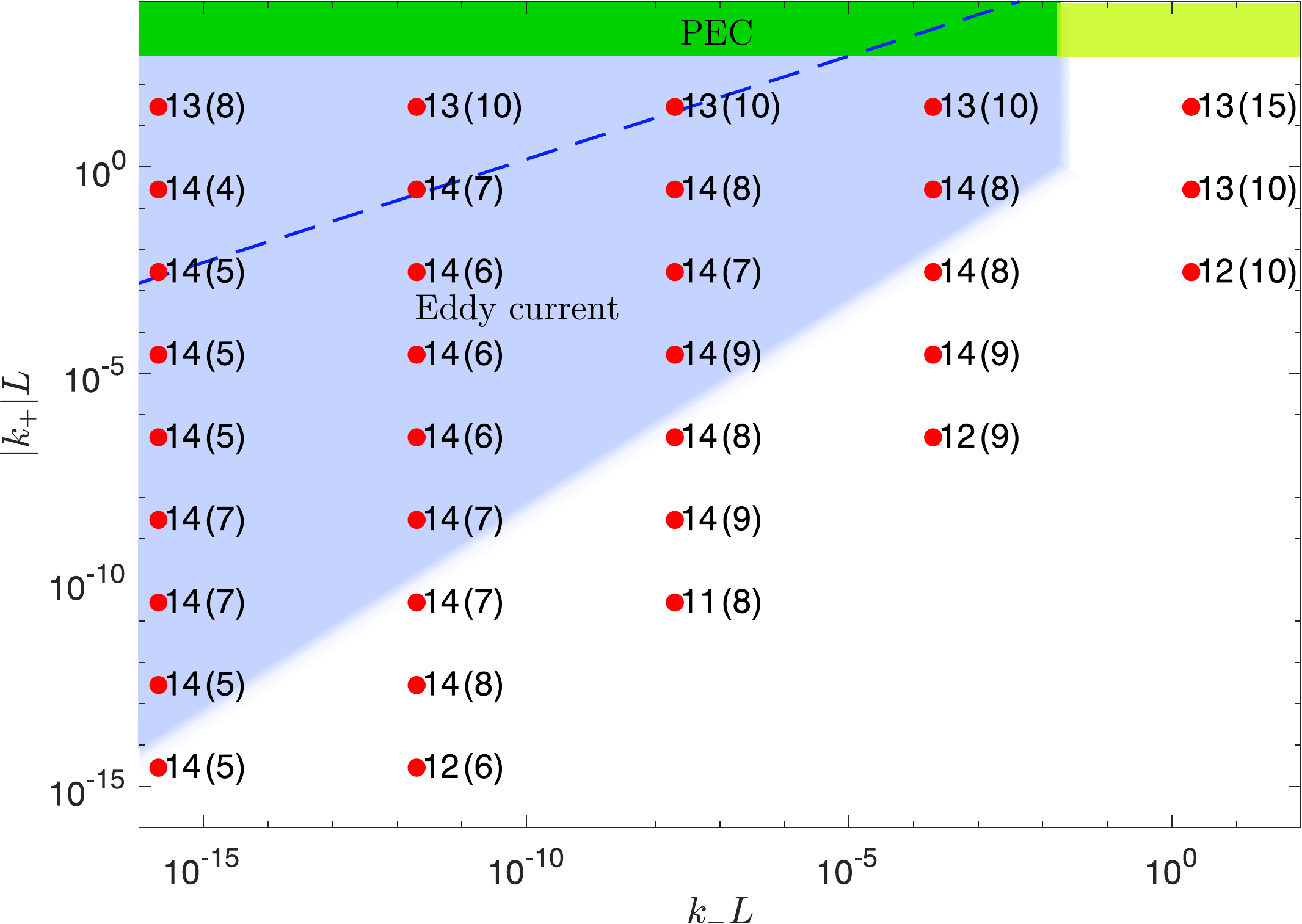}
\caption{\sf 
 Performance in the eddy current regime of Dirac (B-aug0)
on the unit sphere with incident partial waves~\eqref{eq:partialwave} and with reference
solutions from Mie theory. Notation as in Figure~\ref{fig:quadrant0genus1}.
}
\label{fig:quadrant0genus0Mie}
\end{figure}

As customary, we also show, in Figure~\ref{fig:quadrant0genus0Mie}, 
results analogous to those in Figure~\ref{fig:quadrant0genus0}
but for the unit sphere and where the reference solutions are semi-analytic solutions
given by Mie theory -- rather than overresolved, purely numerical, solutions. The
computational domain is 
$\mathcal{D}=\{-2\le x\le 2, -2\le z\le 2\}$. Figure~\ref{fig:quadrant0genus0Mie}
shows that the
solutions obtained with Dirac (B-aug0) in the eddy current regime agree with the
Mie solutions at least as well as they agree with the overresolved reference solutions
in Figure~\ref{fig:quadrant0genus0}. The number of GMRES iterations required is lower, however, because
the scattering problem on the sphere is simpler.

\section{The Maxwell essential spectrum}   \label{sec:Maxwess}

We prove in this section the following result, announced in \cite[Sec. 6]{HelsKarlRos20}.

\begin{thm}   \label{thm:essspectrum}
  Let $\Gamma\subset \R^3$ be a bounded Lipschitz surface, and let 
  $k_\pm\in\C\setminus\{0\}$, $\im(k_\pm)\ge 0$.
  Assume that $\hat k$ is not negative real.
  Then the non-magnetic Maxwell transmission problem 
  MTP($k_-$,$k_+$,$\hat k^2$) defines a Fredholm map,
  in $L^2_{\text{loc}}$ norm of the fields up to $\Gamma$, if and only if
\begin{equation}  \label{eq:dlpcond}
  (1+\hat k^2)/(1-\hat k^2)\notin \sigma_{\text{ess}}(K^{\nu'}_0; H^{1/2}(\Gamma)),
\end{equation}
where  $K^{\nu'}_0$ is the Neumann--Poincar\'e operator,
that is \eqref{eq:acousticNP} with $k=0$. 
\end{thm}

Note that for all passive non-magnetic materials, the technical condition
that $\hat k$ should not be negative real, is always satisfied.
The essential spectrum appears when $\hat \epsilon= \hat k^2$
is negative real. 
In three dimensions, $\sigma_{\text{ess}}(K^{\nu'}_0; H^{1/2}(\Gamma))$
may be non-symmetric with respect to $0$.
Theorem~\ref{thm:essspectrum}
shows in particular that the Fredholm property of MTP($k_-$,$k_+$,$\hat k^2$)
only depends on $\hat \epsilon$, and is not in general symmetric
when replacing $\hat \epsilon$ by $\hat \epsilon^{-1}$.

\begin{proof}
The idea is to use an auxiliary Dirac BIE with parameters
\begin{equation}   \label{eq:essparams}
  \begin{bmatrix} r & \beta & \gamma & \alpha' & \beta' & \gamma'
  \end{bmatrix} = 
  \begin{bmatrix} 1/\hat k &  1 & 1 & 
  1/\hat k & 1/\hat k &  1 \end{bmatrix}.
\end{equation}
Here we are not concerned with false eigenwavenumbers, and 
have tuned the free Dirac parameters only so that we avoid false essential
spectrum, assuming that $\hat k$ is not negative real.
Preconditioning is also irrelevant for this proof, and we let $P=P'=I$.
Consider the Dirac integral operator
\begin{equation}   \label{eq:essDiracBIE}
  E_{k_+}^+(rM')+ ME_{k_-}^-=
  \tfrac 12(rM'+M+E_{k_+}(rM')-ME_{k_-})
\end{equation}
from \eqref{eq:DiracBIEfactors}.
With our choice of parameters~\eqref{eq:essparams}, we have
\begin{equation}
rM'+M=\diag\begin{bmatrix} 
1+\tfrac 1{\hat k} & \tfrac 1{\hat k}+\tfrac 1{\hat k}
  & \mv{\tfrac 1{\hat k}}+\mv{\tfrac 1{\hat k}}  & 1+1 & 1+\tfrac 1{\hat k^2} 
  & \mv{\tfrac 1{\hat k}}+\mv 1
  \end{bmatrix}.
\end{equation}
Modulo compact operators, the operator $E_{k_+}(rM')-ME_{k_-}$
equals the entry-wise product of
\begin{equation}
\begin{bmatrix}
 1-\tfrac 1{\hat k}  & 0 & {\bf  \tfrac 1{\hat k}-\tfrac 1{\hat k}} &
 0 & \hat k-\tfrac 1{\hat k} & {\bf 0}  \\
 1-\tfrac 1{\hat k} & \tfrac 1{\hat k}-\tfrac 1{\hat k} & 
  {\bf \tfrac 1{\hat k}-\tfrac 1{\hat k}} & \hat k-\tfrac 1{\hat k} & 0 & {\bf 1-\tfrac 1{\hat k}} \\
 {\bf 1-\tfrac 1{\hat k}} & {\bf \tfrac 1{\hat k}-\tfrac 1{\hat k}} & 
  {\bf \tfrac 1{\hat k}-\tfrac 1{\hat k}} & {\bf \hat k-\tfrac 1{\hat k}} & {\bf 0} & {\bf 1-\tfrac 1{\hat k}} \\
  0 & 1-1 & {\myvec 0} & 1-1 & 0 & {\bf \tfrac 1{\hat k}-1} \\
 \hat k-\tfrac 1{\hat k^2} & 0 & {\bf 1-\tfrac 1{\hat k^2}} & 1-\tfrac 1{\hat k^2} &
 1-\tfrac 1{\hat k^2} & {\bf \tfrac 1{\hat k}-\tfrac 1{\hat k^2}} \\
 {\bf \hat k-1} & {\myvec 0} & {\bf 1-1} & {\bf 1-1} & 
{\bf 1-1} & {\bf \tfrac 1{\hat k}-1}
\end{bmatrix}
\end{equation} 
(not simplifying some entries to zero in order to show what cancellations
occur) 
and $E_k$ from \eqref{eq:EkCauchyinte}, with $\Phi_{k}$ replaced
by  $\Phi_0$  in the operators, and the factor $ik_-$ in front of the
single layer operators. 
See \cite[Eq.~(132)]{HelsRose20}. 
The compactness of the approximation of $E_k$ was proved
  in~\cite[Lem.~3.20]{AxThesisPub4:06}. 

We now exploit a number of block triangular structures, modulo
   compact operators, in the matrix~\eqref{eq:essDiracBIE}
   which~\eqref{eq:essparams} entails. Recall from \cite[Sec.~5]{HelsRose20}
that the Dirac integral operator acts in the function space
\begin{equation}
  \mH_3 = 
  H^{1/2}(\Gamma)\oplus
  H^{-1/2}(\Gamma)\oplus
  H^{-1/2}(\curl,\Gamma) 
\oplus 
   H^{1/2}(\Gamma)\oplus
  H^{-1/2}(\Gamma)\oplus
 H^{-1/2}(\curl,\Gamma),
\end{equation}
which coincides, up to equivalence of norms, with the function
space $\mE$ from \cite{AxThesisPub4:06}.
Note that we here have omitted the Hodge star present in 
\cite[Eq.~(64)]{HelsRose20}, since we in the present paper 
conform to the standard vector representation of the
magnetic field. 
The function space $H^{-1/2}(\curl,\Gamma)$ consists, roughly
speaking, of tangential vector fields in $H^{-1/2}$ with tangential
curl also in  $H^{-1/2}$.
The precise definition of $H^{-1/2}(\curl,\Gamma)$ on Lipschitz $\Gamma$
is in \cite[Eq.~(65)]{HelsRose20}.

First, we claim that the (5:8,1:4) block is compact.
Indeed, the only possibly non-compact operators are in blocks (5,2)
and (7:8,3:4), and here we have cancellation $1-1$.
This shows that \eqref{eq:essDiracBIE} is a Fredholm operator
on $\mH_3$ if and only if 
its diagonal (1:4,1:4) and (5:8,5:8) blocks are so.   
Second, within these diagonal blocks we have cancellation 
in blocks (1:2,3:4) and
(7:8,5:6). This gives a block triangular structure inside
the diagonal (1:4,1:4) and (5:8,5:8) blocks, which shows 
that \eqref{eq:essDiracBIE} is Fredholm if and only if
its diagonal (1,1), (2,2), (3:4,3:4), (5,5), (6,6) and (7:8,7:8)
blocks are Fredholm operators. 
This is true for the (2,2), (3:4,3:4) and (5,5) blocks,
since these operators are $(1/\hat k)I_{H^{-1/2}(\Gamma)}$, 
$(1/\hat k)I_{H^{-1/2}(\curl,\Gamma)}$ and $I_{H^{1/2}(\Gamma)}$, 
respectively.
Also the (1,1) and (7:8,7:8) blocks in \eqref{eq:essDiracBIE} are 
Fredholm operators since by assumption $(\hat k+1)/(\hat k-1)\notin (-1,1)$,
and the essential spectra of 
$K_0^{\nu'}: H^{1/2}(\Gamma)\to H^{1/2}(\Gamma)$
and $\mv{M}_k^*: H^{-1/2}(\curl,\Gamma)\to H^{-1/2}(\curl,\Gamma)$
are contained in $(-1,1)$. 
This is a  well-known  consequence of the Plemelj symmetrization principle and boundary 
Hodge decompositions.
A precise reference is \cite[Cor.~4.7(i)]{AxThesisPub4:06}, 
which contains the spectral estimates of both $K_0^{\nu'}$ and $\mv{M}_0^*$ upon 
letting $k=0$ and writing $E_k$ in matrix form as in \eqref{eq:EkCauchyinte}.

We conclude that \eqref{eq:essDiracBIE} is  a Fredholm operator
if and only if the (6,6) block is so, which by duality is equivalent to
\eqref{eq:dlpcond}. Moreover, assuming that $\hat k$ is not negative real, 
then the auxiliary 
Maxwell and Helmholtz problems corresponding to the 
parameters $\beta, \gamma, \alpha', \beta', \gamma'$,
all define Fredholm maps, and it follows from the proofs of 
\cite[Props.~8.4, 8.5]{HelsRose20} that \eqref{eq:essDiracBIE} is a Fredholm
operator if and only if MTP($k_-$,$k_+$,$\hat k^2$) defines a Fredholm map.
Combining these two equivalences completes the proof.
\end{proof}

\section{Concluding remarks}   \label{sec:concl}

We conclude with some guiding remarks for the reader chiefly interested in coding the BIEs proposed in this paper.
Eddy current computations are difficult since it is difficult to achieve both
(a) a well-conditioned BIE for computing the density $h$
and (b) a well-conditioned representation of the fields.
Dirac (B-aug0), presented in Section~\ref{sec:B-aug0}, always achieves
(a) and (b) for boundary surfaces of genus $0$.
Dirac (B-aug1), presented in Section~\ref{sec:B-aug1}, almost always achieves
(a) and (b) for boundary surfaces of genus $1$.
There are no null spaces for the systems and there are no low-frequency breakdowns
in these two integral equation reformulations 
for the Maxwell transmission problem.
The only incident field which leads to loss of accuracy with Dirac (B-aug1)
for genus $1$ is \eqref{eq:zcoil}, which is a pathological field which typically 
does not appear in applications.
For \eqref{eq:zcoil} we can accurately compute the fields with Dirac (A$\infty$-aug), 
as presented in Section~\ref{sec:Ainf-aug}.

\appendix
\section{Augmentations}   \label{app:aug}

In this appendix, we derive the augmentations proposed in 
Sections~\ref{sec:finitecond} and \ref{sec:DirB}, following the general principles explained in Section~\ref{sec:chooseaug}.
Throughout this section, we consider the limit as $k_-\to 0$ in the eddy current regime \eqref{eq:goodcond}. To be able to perform the analysis below, we set
the tuning factor $\xi=1$, and also assume that $k_+\to 0$. 
Note that this holds for all ordinary conductors. 
For non-physical limits when $\sigma\to\infty$ so that $k_+\not\to 0$,
the analysis below needs to be adjusted and becomes more complicated. However, the heuristics is that the larger $|k_+|$ is, 
the fewer augmentations are needed. 
We denote by $K$ and $S$, operators of the form
\begin{equation}
  K f(x) = 
  \pv\int_\Gamma v(x,y)\cdot\nabla\Phi_0(y-x) f(y) d\Gamma(y),\qquad x\in \Gamma,
\end{equation}
and
\begin{equation}
  S f(x) = k_+\hat k/\japsigma
  \int_\Gamma u(x,y) \Phi_0(y-x) f(y) d\Gamma(y),\qquad x\in \Gamma,
\end{equation}
for some given vector fields $v(x,y)$ and scalar functions $u(x,y)$.
Note that $|k_+\hat k/\japsigma|\le 1$.

We recall that $\nul(I+ K_0^{\nu'})$ is spanned by the constant function $1$,
that $(K^{\nu'}_0)^*= -K_0^\nu$,
and that $\nul(I- K_0^{\nu})$ is spanned by $f=N_{\Omega_-}(1)$, where $N_{\Omega_-}$ is the 
Dirichlet-to-Neumann map for $\Omega_-$.
Further recall that the spectrum of all operators 
$K^{\nu'}_0, K_0^\nu$ and $\mv M^*_0$ are
contained in $[-1,1]$, with $\sigma(K_0^{\nu'})\cap\{-1,+1\}=\{-1\}$
and $\sigma(K_0^{\nu})\cap\{-1,+1\}=\{+1\}$.
Further $\sigma(\mv M^*_0)\cap\{-1,+1\}=\emptyset$ if the genus of $\Gamma$
is zero and otherwise equals $\{-1,+1\}$.

\subsection*{Dirac (A$\infty$-aug)}

The limit of $I+G$ is seen to be
\begin{equation}    \label{eq:Ainftylimit}
I+G_0^A=
\begin{bmatrix}
I-K^{\nu'}_0  & 0 & \mv 0 &
  0 & S & \mv{ 0}  \\
 K & I-\frac{a-1}{a+1} K_0^\nu & 
  \mv 0 & S & 0 & \mv {S} \\
 \mv K & \mv K & 
 \mv I & \mv S & \mv 0 & \mv {S} \\
 0 & 0 & \mv{ 0}
  & I-\frac{a-1}{a+1} K_0^{\nu'} & 0 & \mv 0 \\
  0 & 0 & \mv 0 & K &
 I-K_0^\nu & \mv K \\
 \mv 0 & \mv 0 & \mv 0 & \mv 0 & 
\mv 0 & \mv I+\frac{a-1}{a+1}\mv M_0^*
\end{bmatrix},
\end{equation}
where we assume that $\hat k\notin (0,\infty)$ so that
$(a-1)/(a+1)\notin [-1,1]$, where $a=\hat k/|\hat k|$.
We note the block triangular structures and that the only non-invertible 
diagonal block is the (6,6) block.
The null space is seen to be spanned by a density of the form
$h_D^A=\begin{bmatrix} h_1 & h_2 & \mv {h_{3:4}} & 0 & f & \mv 0\end{bmatrix}$, and we have $c^1_D h_D^A\ne 0$ at $k_-=0$.
Moreover, the adjoint $(I+G_0^A)^*$ is seen to have null space
spanned by a density of the form
$\begin{bmatrix} 0 & 0 & \mv {0} & h_5 & 1 & \mv{h_{7:8}}\end{bmatrix}$, which is not orthogonal to $b_D^1$.
Therefore the homogeneous (L) augmentation $b_D^1c_D^1$ will  remove  the Dirichlet eigenfield, since $c^1_D h= 0$, that is \eqref{eq:DirEminusaug},
holds for all incident fields under consideration.

\subsection*{Dirac (B-aug0)}
The limit of $I+G$ is seen to be
\begin{equation}  \label{eq:Blimit}
I+G_0^B=
\begin{bmatrix}
I+K^{\nu'}_0  & 0 & \mv K &
 0 & 0 & \mv{ 0}  \\
 0 & I & 
  \mv 0 & 0 & 0 & \mv 0 \\
 \mv 0 & \mv 0 & 
 \mv I & \mv{ 0} & \mv 0 & \mv 0 \\
  0 & S & \mv{ 0} & I-K_0^{\nu'} & 0 & \mv 0 \\
 S & 0 & \mv {S} & K &
 I-K_0^\nu & \mv 0 \\
 \mv 0 & \mv 0 & \mv 0 & \mv 0 & 
\mv 0 & \mv I+\mv M_0^*
\end{bmatrix}.
\end{equation}
In this section, we assume that $\Gamma$ has genus $0$, so that the
(7:8, 7:8) block is invertible.
We note again block triangular structures, and now the non-invertible 
diagonal blocks are the (1,1) and (6,6) blocks.
If $(I+G_0^B)h=0$, then
$h=\begin{bmatrix} c1 & 0 & \mv {0} & 0 & h_6 & \mv{0} \end{bmatrix}$,
where $c\in\C$.
Here $h_6$ satisfies $cS1+ (I-K_0^\nu)h_6=0$, which forces $c=0$
and $h_6=f$ unless $\int_\Gamma S1 d\Gamma$ is zero.
Inspection of the operator $S$ appearing in the (6,1) block shows
that $\int_\Gamma S1 d\Gamma= 2ik_+\hat k/\japsigma |\Omega_+|$, 
where $|\Omega_+|$ denotes the volume of $\Omega_+$. 
Therefore the null space includes
$h_D=\begin{bmatrix} 0 & 0 & \mv {0} & 0 & f & \mv{0}\end{bmatrix}$
and, as $k_+\hat k\to 0$, also 
$h_H=\begin{bmatrix} 1 & 0 & \mv {0} & 0 & 0 & \mv{0}\end{bmatrix}$.
A similar argument applied to the adjoint $(I+G_0^B)^*$ shows that
its null space always includes  
$h^*_D=\begin{bmatrix} 1 & 0 & \mv {h_{3:4}} & 0 & 0 & \mv{0}\end{bmatrix}$, but as $k_+\hat k\to 0$, also
$h^*_H=\begin{bmatrix} 0 & h_2 & \mv {h_{3:4}} & h_5 & f & \mv{0}\end{bmatrix}$.

(1)
The homogeneous (L) augmentation $b^1_Hc^1_H$ will  remove  the additional null space appearing
as $k_+\hat k\to 0$ since $b^1_H$ is not orthogonal to $h^*_H$ and 
since $c^1_H(h_H)\ne 0$. Note that $F_0=0$ for any electromagnetic 
field $F$ satisfying \eqref{eq:hodgedirac}, and consequently 
$c^1_H(h)= 0$ holds for all incident fields under consideration.

(2)
The augmentation for $h_D$ is less straightforward since the (6,6) 
 entry  in
$P'$ vanishes for Dirac (B) at $k_-=0$, causing the corresponding 
fields to vanish.
This means that the field representation \eqref{eq:projdens}
needs to be augmented.
Recalling \eqref{eq:DiracBIEfactors}, we factorize 
\begin{equation}
I+G^B
=  \begin{bmatrix} PE_{k_+}^+D^{-1} & -N E_{k_-}^-\end{bmatrix}
  \begin{bmatrix} DE_{k_+}^+N' \\- E_{k_-}^-P' \end{bmatrix}
= L^B R^B.
\end{equation}
Taking into account \eqref{eq:genericscaling} and Remark~\ref{rem:BH},
we have rescaled the interior fields using the  size  $8\times 8$ matrix 
$D= \diag\begin{bmatrix} \hat k & \hat k & \mv{\hat k} & \hat k^2/\japsigma
& \hat k^2/\japsigma & \mv{\hat k^2/\japsigma} \end{bmatrix}$. 
In what follows, we write $\widetilde E^\pm= \tfrac 12(I\pm \widetilde E)$, where $\widetilde E$ denotes
the (1:4,1:4) = the (5:8,5:8) block in $E_0$, and $\widetilde S$ denotes the (5:8,1:4)
block in $E_0$, but with $ik$ replaced by $ik_+\hat k/\japsigma$ in all operators $S$.
For the right factor $R^B$, corresponding to the field representation, we have
\begin{equation}   \label{eq:rightfactorlimit}
DE^+_{k_+}N'\to 
\begin{bmatrix} \widetilde E^+ & 0 \\ \widetilde S & \widetilde E^+D_1 \end{bmatrix}
\quad\text{and}\quad
E^-_{k_-}P'\to
\begin{bmatrix} \widetilde E^- & 0 \\ 0 & \widetilde E^-D_2 \end{bmatrix},
\end{equation}
with  size $4\times 4$ matrices 
$D_1=\diag\begin{bmatrix} 1 & 1 & \mv{\japsigma^{-1}} \end{bmatrix}$
and $D_2=\diag\begin{bmatrix} 0 & 0 & \mv 1 \end{bmatrix}$.
For the left factor $L^B$, corresponding to the DTP, we have
\begin{equation}
PE^+_{k_+}D^{-1}\to 
\begin{bmatrix} D_3\widetilde E^+ & 0 \\ S & D_4\widetilde E^+ \end{bmatrix}
\quad\text{and}\quad
NE^-_{k_-} \to
\begin{bmatrix} D_5\widetilde E^- & 0 \\ 0 & D_6\widetilde E^- \end{bmatrix},
\end{equation}
with  size $4\times 4$ matrices 
$D_3=\diag\begin{bmatrix} 0 & 1/2 & \mv {1/2} \end{bmatrix}$,
$D_4=\diag\begin{bmatrix} 1 & 1 & \mv 0 \end{bmatrix}$,
$D_5=\diag\begin{bmatrix} 1 & 1/2 & \mv {1/2} \end{bmatrix}$
and $D_6=\diag\begin{bmatrix} \japsigma^{-1} & \japsigma^{-1} & \mv 1 \end{bmatrix}$.

(3)
If $h$ is in the null space for both limit operators in \eqref{eq:rightfactorlimit}, then $h_{1:4}=
\widetilde E^+h_{1:4}+ \widetilde E^-h_{1:4}=0$.
From $\widetilde E^-D_2h=0$ it follows that $h_{7:8}$ is the boundary trace
of  a harmonic  vector field for $\Omega_+$ with tangential boundary
conditions. 
As in \cite[Exc. 10.6.12]{RosenGMA19}, it follows that $h_{7:8}=0$,
as we assume $\Gamma$ to have genus $0$.
Moreover, with $h_{7:8}=0$, we again see $h_D$ appearing from
$\widetilde E^+D_1h=0$.

The (R) augmentation $b^R_Dc^R_D$ has $c^R_D(h_D)\ne 0$, since
$\int_\Gamma f d\Gamma\ne 0$.
We note that $b^R_D$ appears from applying $L^B$ 
to the fields $F^+=0$ and $F^-=E_{k_-}^- e_6$, and it remains to check that
$\widetilde E^- e_6$ is not in the range of $\widetilde E^-D_2$. This in turn follows from 
the divergence theorem, since no divergence-free vector field in 
$\Omega_+$ can have normal component $\nu$. 

(4)
With the (R) augmented field representation \eqref{eq:projdensB0},
we can now make the homogeneous (L) augmentation $b^2_Dc^2_D$
to  remove  the Dirichlet eigenfield.
As for $b^1_Dc^1_D$, we note that $c^2_D h= 0$, that is \eqref{eq:DirEminusaug}, holds for all incident fields under consideration.
However, now we use $b_D^2= e_1$, since this is not orthogonal to
$h_D^*$.

\subsection*{Dirac (B-aug1)}
(1)
We repeat the augmentations for (B-aug0), but now assume
that $\Gamma$ has genus $1$.
To be able to perform the analysis, we also assume that $\japsigma\to \infty$.
Now $\mv I+\mv M_0^*$ has a one-dimensional null space, leading 
to additional null vectors
$h_N$ and 
$h_N^*$, non-zero only in the 7:8  components,  for $I+G_0^B$ and $(I+G_0^B)^*$ respectively.

(2)
Since the Neumann eigenfield can be excited by sources in $\Omega_-$,
we aim to remove the above null space by an inhomogeneous (L)
augmentation. For the abstract equation \eqref{eq:ceqd}, we use
\begin{equation}   \label{eq:NeuLaugbase}
  (\hat k^2/\japsigma) \theta\cdot E^+ - (\hat k^2/\japsigma)\theta\cdot E^- =
 (\hat k^2/\japsigma) \theta\cdot E^0,
\end{equation}
which follows by rescaling the first equation in 
\eqref{eq:maxwtranspr}.
On the  left-hand  side it is the $\theta\cdot E^+$ term which will be dominant,
and the generic size of $E^+$, as discussed in Section~\ref{sec:physics},
motivates the factor $\hat k^2/\japsigma$.
However, in computing the fields with \eqref{eq:projdensB0},
we note that the (7:8,7:8)  entries  in $N'$ vanishes at $k_-=0$.
This indicates that further (R) augmentation is needed for $\Gamma$
of genus $1$.

(3)
Inspecting the limit operators in \eqref{eq:rightfactorlimit},
we see that the new null vector $h_N$ is the boundary trace of 
a harmonic  vector field in $\Omega_+$ with tangential boundary 
conditions,
and hence $h_N$ is in the $\theta$ direction for a torus.
Therefore $c_N^R(h_N)\ne 0$ at $k_-=0$.
We note that $b^R_N$ appears from applying $L^B$ 
to the fields $F^+=E_{k_+}^+e_8$ and $F^-=0$, and it remains to prove that
$\widetilde E^+ e_8$ is not in the range of $\widetilde E^+D_1$. 
For this, we assume that $e_8$ is the boundary trace of the 
interior Neumann eigenfield, and in particular has zero surface curl. 
The assumption that  
$E^+e_8=e_8$ is in the range of $\widetilde E^+D_1$ is seen to 
be equivalent to the existence of a vector field $F$ and a scalar function
$U$ in $\Omega_-$, decaying at $\infty$, with
$\nabla\times F= \nabla U$, $\nabla\cdot F=0$ in $\Omega_-$,
and $F$ having tangential part $e_8$ on $\Gamma$.
It follows that $G= \nabla\times F$ is an exterior Neumann eigenfield.
But since $G$ is a gradient vector field, this forces $G=0$.
It follows that $F$ is curl-free in $\Omega_-$, which contradicts Stokes'
Theorem.
Therefore $e_8$ cannot be in the range of $\widetilde E^+D_1$.

(4)
With the doubly (R) augmented field representation \eqref{eq:projdensB1}, 
we adjust the augmentation $b^1_Hc^1_H$ to $b^2_Hc^2_H$, and
proceed to the inhomogeneous (L) augmentation $b^1_Nc^1_N$
based on \eqref{eq:NeuLaugbase}.
We use \eqref{eq:projdensB1} to write $E^\pm$ in terms of $h$, and
integrate \eqref{eq:NeuLaugbase}  with respect to $wd\Gamma$.
To see that the obtained  left-hand  side defines a bounded functional 
$c^1_Nh$ as in \eqref{eq:bc7}, we rewrite as follows.
That the term $(\hat k^2/\japsigma) \theta\cdot E^+$ depends boundedly
on $h$ is readily seen by inspecting $N'$.
For the term $(\hat k^2/\japsigma) \theta\cdot E^-$, we see from $P'$
that it only depends boundedly on $h_{1:5}$.
Inspection of the (8,6) block of $E_{k_-}$ in \eqref{eq:EkCauchyinte}
reveals that this yields a gradient vector field, with zero circulation around
the torus. Therefore there is no dependence on $h_6$ and the 
Dirichlet (R) augmentation term.
Finally we consider the crucial dependence on $h_{7:8}$, which as it stands
is unbounded.
Note from \eqref{eq:EkCauchyinte} that the (7:8,7:8) block in 
$E^-_{k_-}$ is $(\mv I+\mv M_{k_-}^*)/2$.
The reason for the specific choice of weight $w$, is that 
$w\theta$ is orthogonal to $\ran(\mv I+\mv M_0^*)$.
This allows us to subtract the zero term 
$\int_\Gamma ((\mv I+\mv M_0^*)h_{7:8})_8 wd\Gamma/2$ to obtain 
the third regularized term for $c^1_N h$ in \eqref{eq:bc7}, 
which depends boundedly on $h$.

To motivate the choice of $b^1_N$, we note that $(\mv I+\mv M_0)(h_N^*)_{7:8}=0$.
Hence 
$h^*_N=\begin{bmatrix} 0 & 0 & \mv {0} & 0 & 0 & 0 & w\end{bmatrix}$,
and it follows that $b_N^1=e_8$ is not orthogonal to $h^*_N$.

\end{document}